\documentclass{amsart}

\usepackage{parskip}

\usepackage{layout} 
\usepackage[a4paper,top=1in,bottom=1in,left=1.3in,right=1.3in,marginparwidth=1.75cm]{geometry}
\usepackage[english]{babel}
\usepackage[utf8]{inputenc}
\usepackage[T1]{fontenc}
\usepackage{amsmath}
\usepackage{amssymb}
\usepackage{mathrsfs}
\usepackage{amsthm}
\usepackage{enumerate} 
\usepackage[colorlinks=true,urlcolor=blue, citecolor=red,linkcolor=blue]{hyperref}
\usepackage{enumitem}
\usepackage{comment}
\usepackage{mathtools}
\usepackage{esint}
\usepackage{color}
\usepackage[capitalize]{cleveref}
\usepackage[noadjust]{cite}
\usepackage{cases}

\DeclareFontEncoding{LS1}{}{}
\DeclareFontSubstitution{LS1}{stix2frak}{m}{n}
\DeclareSymbolFont{stixoperp}{LS1}{stix2frak}{m}{n}
\DeclareMathSymbol{\operp}{\mathbin}{stixoperp}{"A8}

\usepackage{tikz}

\newtheorem{theorem}{Theorem}[section]
\newtheorem{proposition}[theorem]{Proposition}

\newtheorem{lemma}[theorem]{Lemma}
\newtheorem{remark}[theorem]{Remark}
\newtheorem{defi}[theorem]{Definition}
\newtheorem{claim}[theorem]{Claim}

\newcommand{\scal}[2]{\left\langle #1,#2 \right\rangle}

\newcommand{\g}{\nabla}

\DeclareMathOperator{\pv}{\mathrm{p.v.}}
\DeclareMathOperator{\lap}{\Delta}
\newcommand{\dr}{\partial}
\newcommand{\bz}{\bar{z}}

\DeclareMathOperator{\Span}{\mathrm{Span}}

\newcommand{\id}{\mathrm{Id}}

\DeclareMathOperator{\sgn}{\mathrm{sign}}

\DeclareMathOperator{\VMO}{\mathrm{VMO}}
\DeclareMathOperator{\BMO}{\mathrm{BMO}}
\newcommand{\dd}{\mathrm{d}}
\newcommand{\ii}{\mathrm{i}}
\newcommand{\ee}{\mathrm{e}}
\newcommand{\loc}{\mathrm{loc}}

\DeclareMathOperator{\supp}{\mathrm{supp}}
\DeclareMathOperator{\dist}{\mathrm{dist}}
\DeclareMathOperator{\Diff}{\mathrm{Diff}}

\DeclareMathOperator{\tr}{\mathrm{tr}}
\DeclareMathOperator{\proj}{\mathrm{proj}}

\DeclareMathOperator{\Hh}{\mathrm{H}}

\newcommand{\R}{\mathbb{R}}
\newcommand{\C}{\mathbb{C}}
\newcommand{\N}{\mathbb{N}}

\newcommand{\s}{\mathbb{S}}
\newcommand{\T}{\mathbb{T}}
\newcommand{\Z}{\mathbb{Z}}
\newcommand{\D}{\mathbb{D}}

\newcommand{\B}{\mathbb{B}}

\newcommand{\Mr}{\mathcal{M}}

\newcommand{\Fr}{\mathcal{F}}

\newcommand{\Dr}{\mathcal{D}}
\newcommand{\Ar}{\mathcal{A}}
\newcommand{\Cr}{\mathcal{C}}

\newcommand{\Sr}{\mathcal{S}}

\newcommand{\Pc}{\mathscr{P}}
\newcommand{\Cc}{\mathscr{C}}

\newcommand{\Ll}{\mathrm{L}}
\newcommand{\W}{\mathrm{W}}
\newcommand{\Cd}{\mathrm{C}}

\DeclareMathOperator{\Mob}{\mathrm{M \ddot{o} b}}
\newcommand{\Chat}{\widehat{\mathbb C}}

\DeclareMathOperator{\Confhat}{\widehat{\mathrm{Conf}}}

\newcommand{\IL}{\mathrm{I^L}}

\newcommand{\ust}{\underset}

\newcommand{\vp}{\varphi}
\newcommand{\eps}{\varepsilon}

\newcommand{\bb}{\mathrm{b}}
\newcommand{\Jj}{\mathrm{J}}

\DeclareMathOperator{\aleq}{\lesssim}

\numberwithin{equation}{section}

\title{Comparing the Möbius and the Loewner energy of a Jordan curve}
\author{Dorian Martino}
\address{Department of Mathematics, ETH Zürich, Rämistrasse 101, 8092 Zürich, Switzerland}
\email{dorian.martino@math.ethz.ch}
\author{Yilin Wang}
\address{Department of Mathematics, ETH Zürich, Rämistrasse 101, 8092 Zürich, Switzerland}
\email{yilin.wang@math.ethz.ch}
\date{\today}

\begin{document}
	
	\begin{abstract} 
		There are two Möbius-invariant energies naturally associated with a Jordan curve, namely, the Möbius energy which arises from knot theory, and the Loewner energy which arises from random curves SLEs and Teichmüller theory. Each of them is finite if and only if the curve is Weil--Petersson. 
    
		We compare these two energies from two perspectives. We first relate both to the fermionic Grunsky operator introduced in \cite{KSS2025}, and express the Möbius energy as a Hilbert--Schmidt norm and the Loewner energy as a log-determinant. 
        From this, we show that these two energies are mutually bounded. 
        Second, we study two optimization problems for the Möbius energy and show that they exhibit similar analytical properties as those for the Loewner energy, although the minimizers are different. In the appendix, we prove regularity up to the boundary for the $1/2$-Laplacian with Dirichlet boundary conditions in dimension $1$, which is instrumental in the study of the optimization problems.
	\end{abstract}
	
	\maketitle

	\tableofcontents
	
	\section{Introduction}


\subsection{Background}

    The goal of this paper is to compare two M\"obius-invariant energies associated with Jordan curves in $\Chat = \mathbb C \cup \{\infty\}$: the M\"obius energy and the Loewner energy. They are different and arise from different contexts, but are known to be finite for the exact same class of curves, the ``Weil--Petersson class'' of quasicircles. Let us give some background.
    
	Although the name ``Weil--Petersson'' originates from the work of Takhtajan--Teo \cite{TT2006} in 2006 referring to the Weil--Petersson metric defined on the space of those curves, the notion itself or closely related notions such as BMO domains have already been studied in many works \cite{AZ1991,C2000,G2000,NS1995,PP2009,PR2009}. A curve $\gamma\colon \T\to \C$, with $\T \coloneqq \R/2\pi\Z$, is called \emph{Weil--Petersson} if it is rectifiable and if its unit tangent vector lies in $\W^{1/2,2}$, see also Section \ref{sec:WP} for a precise definition. Bishop \cite{B2025} has recently proved and gathered a large number of equivalent characterizations. This notion can also be seen as the one-dimensional version of the notion of weak immersions in a critical Sobolev space introduced in \cite{R2014} in dimension 2 and \cite{MR2026} in higher dimensions in the context of Willmore surfaces and their generalizations.
	
	The Loewner energy is one of the most important Möbius invariant functionals for Jordan curves.
	It was introduced in probability theory independently by Rohde and the second author \cite{W2019,RW2021}, who interpreted it as the large deviation rate function of ${\rm SLE}_{\kappa}$ as $\kappa\to0$ \cite{peltola_wang}, the Onsager--Machlup functional of ${\rm SLE}_{\kappa}$ loops for all $\kappa \in (0,4]$ \cite{CarfagniniWang2024}, and by Friz--Shekhar \cite{FS2017} to study finite-energy Loewner drivers. Since then, several equivalent expressions of this energy have been discovered: It can be expressed using the log-derivative of uniformizing maps or $\zeta$-regularized determinants of Neumann jump operators \cite{W2019-2}; as the Fredholm determinant of Neumann--Poincar\'e operator or Grunsky operator \cite{TT2006}; as a K\"ahler potential of Weil--Petersson metric on universal Teichm\"uller space \cite{TT2006,W2019-2}; as the renormalized volume between Epstein surfaces \cite{BBVW2025}; as the renormalized moving frame energy \cite{MW2024}; using the welding \cite{FVW2026}; or as the free energy of Coulomb gas \cite{JV2026,johansson2021strong,wiegmann_zabrodin_22}. 
These identities provide a rich interplay between various subjects.

We use the following definition of the \emph{Loewner energy} of a Jordan curve $\gamma\colon \T\to \C$ (which is also known as the \emph{universal Liouville action} from \cite{TT2006} and identified to the Loewner energy in  \cite[Theorem 1.4]{W2019-2}). We denote by $\Omega$ and $\Omega^*$ the bounded and the unbounded connected components of $\C\setminus \gamma(\T)$,  respectively. Let $f\colon \D\to \Omega$ and $g\colon \D^*\to\Omega^*$, with $\D \coloneqq \{z \in \C \mid |z| < 1\}$,  $\D^*\coloneqq \C\setminus \overline{\D}$, be two biholomorphic maps with $g(\infty)=\infty$. The Loewner energy of $\gamma$ is then given by
	\begin{equation}\label{eq:def_IL}
		\IL(\gamma) \coloneqq \frac{1}{\pi}\int_{\D} \left| \frac{f''}{f'} (z)\right|^2  |\dd z|^2+ \frac{1}{\pi}\int_{\D^*} \left|\frac{g''}{g'} (z)\right|^2 |\dd z|^2 + 4\, \log\frac{ |f'(0)| }{|g'(\infty)|}.
	\end{equation}
See Section~\ref{sec:Loewner} for more details and the expression when the curve passes through $\infty$.

	On the other hand, the \emph{Möbius energy} is a M\"obius invariant functional for knots in $\mathbb R^3$ introduced by O'Hara \cite{O1991,O1992,O1992-2}  through a length-renormalization process, in order to provide a notion of canonical representation of knots. It is defined for an embedded curve $\gamma\colon I\to \C$, with $I=\T$ for closed curve and $I=\R$ for noncompact curves, as follows: 
	\begin{equation}\label{eq:intro_def_Mob}
		\Mob(\gamma) \coloneqq \int_{I\times I} \left(\frac{1}{|\gamma(s)-\gamma(t)|^2} - \frac{1}{\dist_{\gamma}(s,t)^2} \right)\, |\gamma'(s)|\, |\gamma'(t)|\, \dd s\, \dd t,
	\end{equation}
	where $\dist_{\gamma}$ is the distance on $I$ induced by $\gamma$ (i.e., the arclength of the shortest arc between $\gamma(s)$ and $\gamma (t)$). Since the work of Freedman--He--Wang \cite{FHW1994} on minimization of $\Mob$ among knots of a given class and the book of O'Hara \cite{ohara2003} on the geometry of this functional, most of the efforts focused on the analytical properties of this functional, for which we refer to \cite{BRS2018}. Blatt \cite{B2012} proved that the Möbius energy of a knot is finite if and only if the knot has a parameterization by arc-length and its tangent vector lies in $\W^{1/2,2}$. In the case of planar curves, this is exactly the definition of Weil--Petersson curves \cite{B2025}.

	In this manuscript, we compare the Möbius energy of Weil--Petersson curves with their Loewner energy from different viewpoints, which provides partial answers to the questions raised by Bishop in \cite[Section 23]{B2025}.
	Indeed, in addition of sharing a common regularity setting, both functionals are also invariant by Möbius transformations, see the paragraph following Theorem 1.1 in \cite{RW2021} for $\IL$ and Theorem 2.1 in \cite{FHW1994} for $\Mob$. Moreover, both energies depend only on the image $\Gamma = \gamma (\mathbb T)$ and not on its parameterization.

\subsection{Comparison via fermionic Grunsky operators}

	Although $\IL$ can also be expressed as the Dirichlet energy of the harmonic extension of its tangent vector \cite{VW2020}, a direct comparison between the Loewner and M\"obius energies for a given curve seems complicated and not optimal, since both energies are non-local and satisfy different additive properties \cite[Corollary 4.4]{W2019} or \cite[Lemma 2.6]{BJRW2025}. 
    
    Instead, we express both functionals in terms of well-chosen operators. On the one hand, it is known that $\IL$ can be expressed in terms of (classical) Grunsky operators, see \cite{B2025,FVW2026,JV2026,TT2006}. 
    We prove that the Möbius energy can also be expressed in terms of well-chosen operators. 
    
    Without loss of generality, we assume that the Jordan curve $\Gamma \subset \C$.
	\begin{theorem}\label{thm:intro_cauchy}
		Let $\gamma\colon \T\to\C$ be a Weil--Petersson curve with $\Gamma = \gamma(\T)$. 
        The Cauchy operator $\Cc_{\Gamma}\colon \Ll^2(\Gamma)\to \Ll^2(\Gamma)$ is defined as
		\begin{equation*}
			\forall f\in \Ll^2(\Gamma),\qquad \Cc_{\Gamma}(f)\colon z\mapsto \frac{1}{2\ii \pi}\, \pv \int_{\Gamma} \frac{f(w)}{z-w}\, \dd w.
		\end{equation*}
		Let $\Cc^*$ denote its adjoint operator and $\|\cdot\|_{\rm HS}$ the Hilbert--Schmidt norm of an operator, then we have 
		\begin{equation*}
			\Mob(\gamma) = 4 + 2\pi^2\, \|\Cc_{\Gamma}-\Cc_{\Gamma}^*\|^2_{\mathrm{HS}\left(\Ll^2(\Gamma)\right)}.
		\end{equation*}
	\end{theorem}
    This result is a direct consequence of \cite{BB2007,KS1978}, which relate the Cauchy operator to conformal angles. We record its proof in Section \ref{sec:Mobius-Cauchy} for completeness.
	
	In order to relate to the Loewner energy, 
    we will use \textit{fermionic Grunsky operator} $\mathrm{Gr}^{1/2}_{f}$ on half-differential forms introduced in \cite{KSS2025}.  
  This operator has a kernel representation (by \cite[Corollary 4.6]{KSS2025}, and see \eqref{eq:def-Grunsky}) given by $(2\ii \pi)^{-1}$ times
	\begin{equation*}
	K(z,\zeta)\coloneqq \frac{1}{\zeta-z} - \frac{\sqrt{f'(z)}\, \sqrt{f'(\zeta)}}{f(\zeta)-f(z)}, \qquad z,\zeta \in \s^1
	\end{equation*}  
    where $f$ is a conformal map from the unit disk $\D$ onto the bounded connected component $\Omega$ of $\C \setminus \Gamma$.
    We recall its background and definition in Section \ref{sec:Grunksy}, where we also establish many formulas. For instance, using Theorem~\ref{thm:intro_cauchy}, we prove that the Möbius energy is described using the Hilbert--Schmidt norm of this fermionic Grunsky operator (Proposition~\ref{pr:Mob-Grunsky}), with $\s^1=\dr\D$:
		\begin{equation*}
			\Mob(\Gamma)-4 = 4\pi^2\, \left\| \mathrm{Gr}^{1/2}_f \right\|^2_{\rm HS\left(\overline{\Ar^{1/2}(\D) } \right) } = \int_{\s^1\times\s^1} |K(z,\zeta)|^2\, |\dd z|\, |\dd \zeta|.
		\end{equation*}
    This result also answers a question of \cite{KSS2025} about characterization of Weil--Petersson curves. Our result implies that as long as $\mathrm{Gr}^{1/2}_f$ is defined, a Jordan curve $\Gamma$ is Weil--Petersson if and only if its fermionic Grunsky operator is Hilbert--Schmidt (Remark~\ref{rem:SP}). 
    
    Then we study, analogously to \cite{FVW2026}, the action of conformal welding in this space of half-differential forms. By decomposing the space into two orthogonal subspaces, we derive algebraic identities leading to a new expression of the Loewner energy as a log-determinant involving the fermionic Grunsky operator (Theorem~\ref{th:Loewner-Grunsky}):
    \begin{equation*}
			\IL(\Gamma) = 24\log\det\left(I+\left[\mathrm{Gr}^{1/2}_{f}\right]^*\, \mathrm{Gr}^{1/2}_{f}\right).
		\end{equation*}
    
    Comparing the values of $\Mob$ and $\IL$ for each curve boils down to comparing the Hilbert--Schmidt norm and the log-determinant. We prove the following inequalities in Section \ref{sec:Mobius-Loewner}.
	\begin{theorem}[See Theorem~\ref{thm:compare}]
		Let $\Gamma\subset \C$ be a $K$-quasicircle. Then there exists a constant $C_K>0$ that depends only on $K$ such that 
		\begin{equation*}
			\IL(\Gamma) \leqslant \frac{6}{\pi^2}\left(\Mob(\Gamma)-4\right) \leqslant C_K\, \IL(\Gamma).
		\end{equation*}
	\end{theorem}

    There is also a non-explicit upper bound of the quasiconformal constant $K$ in terms of $\IL$ in Proposition 2.9 of \cite{RW2021}, so a priori, the second inequality can be written solely in terms of $\IL$.
	

    \subsection{Comparison via optimization problems}
    Another approach to compare both energies is to compare the solutions to similar variational problems. 
    The following minimization problems for Loewner energy were studied in \cite{BJRW2025,MRW2026,RW2021,W2025}.
    
    Given $n\geqslant 3$ distinct points $Z\coloneqq (z_1,\ldots,z_n)\in \C^n$ and a relative 
    isotopy class $\tau$ of curves passing through $Z$, 
    we consider the following optimization problem
	\begin{equation}\label{eq:Pc1-IL}
		\Pc^1_{\IL}(Z,\tau) \coloneqq \inf \left\{ \IL(\gamma) \mid \gamma\colon \T\to \Chat\text{ is Weil--Petersson},\ z_1,\ldots,z_n\in \gamma(\T), \gamma \in \tau \right\}.
	\end{equation}
	
	A very closely related problem consists in filling an arc by an optimal curve. Given $\gamma_0\colon \T\to \Chat$ and $u,v\in \T$, we consider the following problem  
	\begin{equation}\label{eq:Pc2-IL}
		\Pc^2_{\IL}(\gamma_0,u,v) \coloneqq \inf \left\{ \IL(\gamma) \mid \gamma\colon \T\to \Chat \text{ is Weil--Petersson},\ \gamma|_{(u,v)} = \gamma_0|_{(u,v)} \right\}.
	\end{equation}
    
	The solutions to both optimization problems are unique and well understood.  
    
    \begin{remark}
        Since the global minimizers of both Möbius and Loewner energies are given by circles, the minimizers when $n = 2, 3$ are trivially given by circles passing through the given points, which is not unique when $n = 2$. This is the reason we assume $n \geqslant 3$. 
    \end{remark}
    
    \begin{remark}\label{rem:two_op}
        We note that the solution 
    $\gamma_{Z,\tau}$ of \eqref{eq:Pc1-IL} is also the solution of \eqref{eq:Pc2-IL} when we take $\gamma_0 = \gamma_{Z,\tau}$ and $u,v$ be such that $\gamma_0|_{(u,v)}$ contains all points in $Z$.
    \end{remark}

Regarding $\Pc^1_{\IL}$, the uniqueness is rather nontrivial and is proved in \cite{BJRW2025}. The regularity follows from explicit computations in \cite{MRW2026}. 
\begin{theorem}[\cite{BJRW2025,MRW2026}]
		Let $n \geqslant 3$. For any distinct points $Z= (z_1,\ldots,z_n)\in\C^n$ and relative isotopy class $\tau$,  there exists a unique minimizer $\gamma = \gamma_{Z,\tau}$ of the Loewner energy. 
        Moreover, there exists $c>0$ such that for any $x,y\in\s^1$ there holds
            \begin{equation*}
                |\gamma'(x)-\gamma'(y)| \leqslant c\, |x-y|\, \left( 1+ \big| \log|x-y| \big| \right).
            \end{equation*}
	\end{theorem}
See Section~\ref{sec:loewner_arc_completion} for more details.

    For $\Pc^2_{\IL}$, we will use an equivalent definition of the Loewner energy as the Dirichlet energy of its driving function (see Section~\ref{sec:Loewner} for the background). From this definition, the solution is easily identified and the regularity follows from the correspondence between the H\"older regularity of the curve and of its Loewner driving function \cite{W2014,RW2021,Lind_Tran}. More precisely, we show the following.
	\begin{theorem}[See Theorem \ref{thm:loewner_arc_completion}]
		Given a $\Cd^{1,\alpha}$ Jordan curve $\gamma_0\colon \T\to \Chat$, with $\alpha > 1/2$ and $u,v\in\T$, there is a unique solution to $\Pc^2_{\IL}(\gamma_0,u,v)$. The solution is given by completing $ \gamma_0|_{(u,v)}$ with the unique hyperbolic geodesic in the slit complement $\Chat \setminus \gamma_0([u,v])$ connecting the tips of the slit.  It is analytic away from the slit tips, while the completed Jordan curve is $\Cd^{1,1/2}$ at the junctions.  
	\end{theorem}
The reason that we start with a $\Cd^{1,\alpha}$ curve with $\alpha > 1/2$, is that this is the largest H\"older class contained in the Weil--Petersson class (there are examples of $\Cd^{1,1/2}$ curves that are not Weil--Petersson).
    
	These results inspired us to study in this paper the analogous problems for M\"obius energy:
	\begin{equation}\label{eq:P1Mob}
		\Pc^1_{\Mob}(Z,\tau) \coloneqq \inf \left\{ \Mob(\gamma) \mid \gamma\colon \T\to \Chat\text{ is Weil--Petersson},\ z_1,\ldots,z_n\in \gamma(\T), \gamma\in\tau \right\},
	\end{equation}
	and
	\begin{equation}\label{eq:P2Mob}
		\Pc^2_{\Mob}(\gamma_0,u,v) \coloneqq \inf \left\{ \Mob(\gamma) \mid \gamma\colon \T\to \Chat\text{ is Weil--Petersson},\ \gamma|_{(u,v)} = \gamma_0|_{(u,v)} \right\}.
	\end{equation}
	We prove the existence of a solution to $\Pc^1_{\Mob}$ and $\Pc^2_{\Mob}$ in Section \ref{sec:wp_attainment}. We study the regularity of solutions to $\Pc^1_{\Mob}$ and $\Pc^2_{\Mob}$ in Sections \ref{sec:reg-P1mob} and \ref{sec:reg-P2mob} respectively. We study the uniqueness of critical points in Section \ref{sec:uniqueness}.
    We gather these results as follows.
	\begin{theorem}\label{th:main-Mob}
		\begin{enumerate}
            \item There is a unique solution to $\Pc^1_{\Mob} (Z,\tau)$. Each solution $\gamma\colon \T\to\Chat$ is analytic outside of $Z$ and has second derivative globally in $\BMO$. In particular, the tangent vector is log-Lipschitz across the marked points: there exists $c>0$ such that for any $x,y\in\T$ there holds
            \begin{equation*}
                |\gamma'(x)-\gamma'(y)| \leqslant c\,  |x-y|\, \left( 1+ \big| \log|x-y| \big| \right).
            \end{equation*}

			\item There is a unique solution to $\Pc^2_{\Mob}$ which is analytic inside the completed part and globally $\Cd^{1,1/2}$.
		\end{enumerate}
	\end{theorem}

These statements about the regularity and uniqueness are the same as for the optimizations of the Loewner energy. 

Uniqueness and classification results concerning critical points of the Möbius energy are few in the literature. Blatt \cite{B2020} proved that round circles are the only Jordan curves that are critical points of the Möbius energy. Blatt--Gilsbach--Reiter--von der Mosel \cite{BGRM2025} recently proved the existence of at least two critical points for knots in certain torus knot classes. 
In our setting, we benefit from the planar setting, which allows us to use the holomorphic motion of Jordan curves \cite{F1988,AM2001,BF2016}. We also use the fact that the Möbius energy has an explicit expression in terms of the parameterization of the curve, unlike the Loewner energy. 
Assuming that we have two distinct minimizers $\gamma_{-1}$ and $\gamma_1$, to either $\Pc^1_{\Mob}$ or $\Pc^2_{\Mob}$, lying in the same relative isotopy class, then by \cite{F1988,AM2001,BF2016}, we can find, for any $\theta>0$ small enough, a continuous family $(\gamma_z)_{z\in\overline{\D}}$ of embeddings of $\T$ in $\C$ such that 
\begin{enumerate}
    \item For $\zeta\in\s^1$, we have $\gamma_{\zeta}(\T)=\gamma_1(\T)$ as unparametrized Jordan curves, if $\Re(\zeta)>\theta$ and $\gamma_{\zeta}(\T)=\gamma_{-1}(\T)$ if $\Re(\zeta)<-\theta$,
    \item For each $s\in\s^1$, the map $(z\in\D\mapsto\gamma_z(s))$ is holomorphic,
    \item For each $z\in\overline{\D}$, the curve $\gamma_z(\T)$ is admissible for either $\Pc^1_{\Mob}$ or $\Pc^2_{\Mob}$.
\end{enumerate}
Using this, we compare the kernels of the fermionic Grunsky operator associated with the curves $\gamma_{-1}(\T)$ with that of $\gamma_1(\T)$. We prove that their distance in $\Ll^2$ is controlled by $\theta$. This controls the distance between the cross-ratios of any four points on the two curves. Since these curves coincide at least three points, we deduce that they are equal by letting $\theta\to 0$.

 The logarithmic singularities of the derivative of the solution to $\Pc^1_{\Mob}$ correspond to logarithmic singularities of the curvature, see Remark \ref{rk:log-singularities}. These are reminiscent of the logarithmic singularities of the uniformizing maps for $\Pc^1_{\IL}$ proved in \cite[Theorem 4.3]{MRW2026}.
    
	The main difficulty in the regularity part of \Cref{th:main-Mob} is studying boundary regularity for fractional degenerate elliptic operators. Already for the fractional Laplacian, it has been proved in \cite{AFLY2025,B2016,CS1998,G2014,KW2025,RoS2014} that the boundary regularity is very different from that of the local elliptic operators. Moreover, we work here in dimension $n=1$, which is the only dimension for which $n/2=1/2 \in(0,1)$. This implies additional difficulties when working with the $1/2$-Laplacian. We could not find the boundary regularity for fractional Laplace systems in dimension 1 in the literature even though we believe that it is well-known to experts; therefore, we wrote it in Section \ref{sec:CZ}. 
	
	The two types of constraints have different impacts on the regularity of solutions. Concerning $\Pc^2_{\Mob}$, the two junctions are genuine nonlocal boundary points, reminiscent of the Dirichlet boundary condition. As proved in \cite{RoS2014}, the expected tangent regularity is $\Cd^{0,1/2}$. Concerning $\Pc^1_{\Mob}$, there is no boundary, the right-hand side of the Euler--Lagrange equation is given by Dirac masses which lie in $\W^{-1,\infty}$ in dimension 1. Hence, the curve is analytic on the open arcs between the points $z_i$, but the singularities at the points are of $\log$-type for the second derivatives. Hence, we get $\gamma''\in\BMO$, which implies that $\gamma'$ is log-Lipschitz, see Lemma \ref{lm:logLip}.
	
	
	Describing the geometry of solutions to $\Pc^1_{\Mob}$ and $\Pc^2_{\Mob}$ is a challenging problem.
	Concerning $\Pc^1_{\Mob}$, we prove that the Möbius invariance implies balancing identities showing that the singularities are completely described by a polynomial of degree at most $n-4$ in Section \ref{sec:balance}. We also prove the following unique continuation property in Section \ref{sec:unique_mob}.
	\begin{theorem}
		Let $\eta$ be a solution to $\Pc^2_{\Mob}$ in \eqref{eq:P2Mob}. If there exists an interval $I\subset \T\setminus (u,v)$ such that the optimized part $\eta|_I$ parametrizes a round (circular or straight) arc, then $\gamma_0|_{(u,v)}$ is also a round arc and $\eta(\T)$ is a circle (or a line).
	\end{theorem}
By Remark~\ref{rem:two_op}, the same statement is also true for $\Pc^1_{\Mob}$.
 We also prove the same unique continuation property for the Loewner energy, which is much simpler, see Section~\ref{sec:unique_cont_IL}.
    
    Despite all these similarities, we prove by a careful analysis of the fermionic Grunsky operator in Section \ref{sec:distinct} that solutions to $\Pc^1_{\IL}(Z,\tau)$ and $\Pc^1_{\Mob}(Z,\tau)$ are different, except when they are round circles.

	\textbf{Structure of the paper:} In \Cref{sec:Sobolev}, we record all the notations and analytic properties surrounding Sobolev spaces that will be used in this article. In Section \ref{sec:WP}, we record all the regularity estimates implied by the Möbius energy. In \Cref{sec:Loewner}, we recall the expression of the Loewner energy in terms of the driving function. In \Cref{sec:loewner_arc_completion}, we provide short proofs of existence and uniqueness of solutions to $\Pc^2_{\IL}$ and discuss $\Pc^1_{\IL}$. In \Cref{sec:Mobius-Cauchy}, we study the relation between the Möbius energy with certain operators such as the Cauchy operator. In \Cref{sec:Grunksy}, we discuss the notion of fermionic Grunsky operator and relate it to the M\"obius energy. In \Cref{sec:log_det_I}, we express $\IL$ in terms of the fermionic Grunsky operator. In \Cref{sec:Mobius-Loewner}, we compare $\Mob$ and $\IL$. In \Cref{sec:wp_attainment}, we prove the existence of solution to $\Pc^1_{\Mob}$ and $\Pc^2_{\Mob}$. 
    In \Cref{sec:reg-P2mob} and \Cref{sec:reg-P1mob}, we prove the regularity of the solutions to $\Pc^2_{\Mob}$ and $\Pc^1_{\Mob}$ in  \Cref{th:main-Mob}, respectively. In \Cref{sec:uniqueness}, we prove the uniqueness statement of \Cref{th:main-Mob}. In \Cref{sec:balance}, we prove a balancing law for the singularities of $\Pc^1_{\Mob}$. In \Cref{sec:Round_pieces}, we prove the unique continuation property to the solution of $\Pc^2_{\IL}$ and $\Pc^2_{\Mob}$. In Section \ref{sec:distinct}, we show that the solutions to $\Pc^1_{\IL}$ and $\Pc^1_{\Mob}$ are different.

	\textbf{Notation:} Given two non-negative numbers $A,B$, we write $A\aleq B$ if $A\leqslant C B$ for a constant $C>0$ whose dependence is indicated when it is not universal. We write $A\approx B$ if $A\aleq B$ and $B\aleq A$. We will denote $\Chat=\C\cup\{\infty\}$. We distinguish between $\T = \R/(2\pi\Z)$ and $\s^1 = \dr\D$. We also use $\T_\ell = \R / \ell \Z$.
    Given $x\in \R$, $x\in\T$, $x \in \T_\ell$ or $x\in \s^1$ and $r>0$, we denote $B_r(x)$ the ball in $\R$, $\T$, $\T_{\ell}$ or $\s^1$ centered at $x$ with radius $r$.

	\textbf{Acknowledgement:} We thank Wolfgang Staubach for the comments on an earlier version of the manuscript. D.M. thanks Armin Schikorra for indicating the reference \cite{KW2025} and Fredrik Viklund for stimulating discussions leading to Proposition \ref{pr:distinct}. D.M. is funded by Swiss National Science Foundation, project SNF 200020\textunderscore 219429. Y.W. is supported by the Swiss State Secretariat for Education, Research and Innovation [MB25.00004] and the Simons Foundation International [SFI-MPS-PP-00012621-19]. 

\textbf{AI declaration:} All the mathematical ideas come from the authors. In the early stages of the project, we used AI tools (ChatGPT 5.6) to check the computations in Section~\ref{sec:differentiate-welding}, Section~\ref{sec:log_det_I}, and Section~\ref{sec:unique_mob} after writing them up. We also included AI's suggestions to simplify the computations in the appendix. 
We used AI (and Google, Kagi) for the literature search, and AI (and Grammarly) for grammar and spell-checking. All mathematical validation, final proof decisions, and final wording remain the sole responsibility of the 
authors.

\section{Preliminaries}\label{sec:Preliminaries}

	\subsection{Sobolev spaces and fractional Laplacian}\label{sec:Sobolev}
	
	Given an interval $I\subset \R$ or $I\subset \R/(\ell \Z)$ for some $\ell>0$, we denote for $0<s<1$ and $p\geqslant 1$,
	\begin{equation*}
		\W^{s,p}(I;\C)\coloneqq \left\{ u\in L^p(I) : [u]^p_{\W^{s,p}(I)}\coloneqq \int_{I\times I} \frac{|u(x)-u(y)|^p}{|x-y|^{1+sp}}\, \dd x\, \dd y<\infty \right\}. 
	\end{equation*}
	We denote $|D|^s\coloneqq (-\lap)^{s/2}$. If we denote $c_s \coloneqq 4^s\, \pi^{-1/2} \, \Gamma(\frac{1+s}{2})\, |\Gamma(-\frac{s}{2})|^{-1}>0$, then for any $f\colon \R\to \R$, it holds
	\begin{equation*}
		\forall x\in \R,\qquad |D|^s f(x) = c_s\, \pv \int_{\R} \frac{ f(y)-f(x) }{|y-x|^{1+s}}\, \dd y.
	\end{equation*}
    We denoted $\pv$ the principal value. 
	The Hilbert transform $\Hh$ is defined for any $\vp\in \Sr(\R)$ by	
	\begin{equation*}
		\forall \xi\in\R,\qquad  \Fr\left[\Hh(\vp)\right](\xi) = -\ii \sgn(\xi)\, \Fr[\vp](\xi),
	\end{equation*}
    where $\sgn(\xi)$ denotes the sign of $\xi$, defined to be $0$ at $\xi=0$. We have the explicit formula
	\begin{equation*}
		\forall x\in\R,\qquad  \Hh(\vp)(x) = \frac{1}{\pi}\, \pv \int_{\R} \frac{\vp(z)}{x-z}\, \dd z.
	\end{equation*}
	In particular, we have 
	\begin{equation}\label{eq:Hilbert-Laplace}
		\Hh[\vp'] = (-\lap)^{1/2}\vp.
	\end{equation}
	
	\subsection{Weil--Petersson curves and M\"obius energy}\label{sec:WP}
	
	In this section, we recall the notion of Weil--Petersson curves and record the estimates needed for computing their Möbius energy. There are many equivalent definitions of Weil--Petersson curves, the one we state below is from \cite[Definition 6]{B2025}.
	
	\begin{defi}\label{def:generalized_wp_curve}
		A Jordan curve $\Gamma\subset\Chat$ is a \emph{Weil--Petersson Jordan curve} if there is a an arclength parameterization $\gamma\colon I\to \Gamma$ with $I =  \R$ \textnormal{(}if $\infty \in \Gamma$\textnormal{)} or $I =  \T_\ell$ \textnormal{(}if $\infty \notin \Gamma$ and $\Gamma$ has total length $\ell$\textnormal{)}, which satisfies the following two conditions:
		\begin{enumerate}
			\item The curve is chord-arc: There exists $C_{\rm ca}>0$ such that
			\begin{equation*}
				\forall x,y\in I,\qquad \frac{1}{C_{\rm ca}} \leqslant \frac{|\gamma(x)-\gamma(y)|}{|x-y|} \leqslant 1.
			\end{equation*}
			
			\item The tangent vector $\gamma'$ lies in $\W^{1/2,2}(I)$:
			\begin{equation*}
				\int_{I\times I} \frac{|\gamma'(x)-\gamma'(y)|^2}{|x-y|^2}\, \dd x\, \dd y<\infty.
			\end{equation*}
		\end{enumerate}
	\end{defi}

	Blatt \cite{B2012} (see also Bishop \cite{B2025}) proved that $\Gamma$ is Weil--Petersson if and only if it has finite Möbius energy defined as follows: 
	\begin{equation}
		\Mob(\Gamma)  = \Mob(\gamma) \coloneqq \int_{I\times I} \left(\frac{1}{|\gamma(s)-\gamma(t)|^2} - \frac{1}{\dist_{I}(s,t)^2} \right)\, \dd s\, \dd t,
	\end{equation}
 where $\dist_I(s,t)$ is the Euclidean distance between $s$ and $t$ on $I$, which coincides with $|s-t|$ when $I = \mathbb R$. This formula coincides with \eqref{eq:intro_def_Mob} because we used the unit-speed parameterization.

 \begin{theorem}[Theorem 2.1 of \cite{FHW1994}]
 Let $\Gamma$ be a Weil--Petersson curve.
 If $\infty \notin \Gamma$ and $T$ is a Möbius transformation of $\Chat$.  Then 
 \begin{itemize}
     \item $\Mob (T(\Gamma)) = \Mob (\Gamma)$, if $\infty \notin T(\Gamma)$;
     \item $\Mob (T(\Gamma)) = \Mob (\Gamma) - 4$, if $\infty \in T(\Gamma)$.
 \end{itemize}
   In particular, if $\infty \notin \Gamma$, then $\Mob (\Gamma) \geqslant 4$, where the equality holds if and only if $\Gamma$ is a circle.
 \end{theorem}

	We next record the relevant estimates from \cite{B2012} related to the Möbius energy of $\Gamma$. For completeness, we provide short proofs of these estimates (this is merely the restriction $j=2$, $p=1$ of \cite{B2012} where the computations can be slightly simplified). We first prove in Lemma \ref{lem:distortion} that a bound on the Möbius energy controls the chord-arc constant. Then we prove in Lemma \ref{lem:tangent-control} that we control the tangent vector in $\W^{1/2,2}$ by the Möbius energy.
	
	\begin{lemma}\label{lem:distortion}
		For every $E_0<\infty$ there exists $c(E_0)>0$ with the following property. If $\gamma\colon I\to\C$ is an arc-length parameterization of a Weil--Petersson curve such that  $\Mob(\gamma)\leqslant E_0$, then we have 
		\begin{equation}\label{eq:chord-arc}
			\forall s,t\in I,\qquad |s-t|\geqslant |\gamma(s)-\gamma(t)|\geqslant c(E_0)\, |s-t|.
		\end{equation}
       
	\end{lemma}
	
	\begin{proof}
		Fix $s<t$. The upper estimate $|\gamma(t)-\gamma(s)|\leqslant |s-t|$ follows from the unit-speed parameterization.  We prove the lower estimate.
		
		Up to translating the time parameters, we can assume that $s=0$ and $t=|s-t|$.  If $|\gamma(0)-\gamma(t)|\geqslant t/20$, the lower bound holds if we choose $c(E_0) < 1/20$. Thus, we can assume that $|\gamma(0)-\gamma(t)|<t/20$.
		For every integer $j\geqslant 0$, we denote
		\begin{equation*}
		r_j\coloneqq 2^j |\gamma(0)-\gamma(t)|,
		\qquad
		A_j\coloneqq [r_j,2r_j],
		\qquad
		B_j\coloneqq [t-2r_j,t-r_j].
		\end{equation*}
		The interior of the intervals $A_j$ are pairwise disjoint, and so are the intervals $B_j$.
		We denote $J\coloneqq \{j\in\N : r_j\leqslant t/20\}$ and fix $j\in J$.
		
		Given $x\in A_j$ and $y\in B_j$, it holds
		\begin{align*}
			|\gamma(x)-\gamma(y)| &\leqslant |\gamma(x)-\gamma(0)| +|\gamma(0)-\gamma(t)| +|\gamma(t)-\gamma(y)|\\[2mm]
			&\leqslant x+|\gamma(0)-\gamma(t)|+(t-y)\\[2mm]
			&\leqslant 2r_j+r_j+2r_j = 5r_j.
		\end{align*}
		Moreover, we have using $r_j\leqslant t/20$,
		\begin{equation*}
			|x-y|\geqslant t-4r_j\geqslant \frac{4}{5} \, t.
		\end{equation*}
		Consequently, we obtain using $r_j\leqslant t/20$, for some unversal $c_0>0$,
		\begin{equation*}
			\frac{1}{|\gamma(x)-\gamma(y)|^2} -\frac{1}{|x-y|^2} \geqslant \frac{1}{25r_j^2}-\frac{25}{16t^2} \geqslant \frac{1}{25r_j^2} - \frac{1}{256r_j^2} \geqslant \frac{c_0}{r_j^2}.
		\end{equation*}
		Integrating over the rectangle $A_j\times B_j$, we obtain 
		\begin{equation*}
			\int_{A_j\times B_j} \left( \frac{1}{|\gamma(x)-\gamma(y)|^2} -\frac{1}{|x-y|^2} \right)\, \dd x\, \dd y \geqslant \frac{c_0}{r_j^2}\, |A_j|\, |B_j| = c_0.
		\end{equation*}
		We sum over all $j\in J$
		\begin{equation*}
			\Mob(\gamma) \geqslant c_0\, \# J \geqslant c_0\, \left\lfloor\log_2\left(\frac{t}{20\, |\gamma(0)-\gamma(t)|}\right) \right\rfloor.
		\end{equation*}
		We end up with the following estimate, which implies  \eqref{eq:chord-arc}
		\begin{equation*}
			|\gamma(0)-\gamma(t)|\geqslant C^{-1} \ee^{-C\, \Mob(\gamma)}\, t.\qedhere 
		\end{equation*}
	\end{proof}

    We now prove that in addition of the chord-arc constant, the Möbius energy also controls the $\W^{1/2,2}$-norm of the unit tangent vector.
    
	\begin{lemma}\label{lem:tangent-control}
		There exists $C>0$ verifying the following. Let $I\in\{\R,\R/(\ell\, \Z)\}$ for some $\ell>0$ and $\gamma\colon I \to\C$ be an arclength parametrized Weil--Petersson curve. Then it holds
		\begin{equation}\label{eq:Hhalf-control}
			[\gamma']_{\W^{1/2,2}(I)}^2 \coloneqq  \int_{I^2} \frac{|\gamma'(x)-\gamma'(y)|^2}{|x-y|^2}\dd x\dd y \leqslant C \Mob(\gamma).
		\end{equation}
	\end{lemma}
	
	\begin{proof}
		We first prove it for $I=\R$. For $w>0$, we apply the fundamental theorem of calculus and the equality $|\gamma'|=1$, to obtain 
		\begin{align*}
			w^2-|\gamma(s+w)-\gamma(s)|^2 &=\int_0^w\int_0^w \left(1-\scal{\gamma'(s+t_1)}{\gamma'(s+t_2)}_{\C} \right)\dd t_1\, \dd t_2\\[2mm]
			&=\frac{1}{2}\int_0^w\int_0^w |\gamma'(s+t_1)-\gamma'(s+t_2)|^2\, \dd t_1\, \dd t_2.
		\end{align*}
		Therefore, it holds, using \eqref{eq:chord-arc}
		\begin{align*}
			\frac{1}{ |\gamma(s+w)-\gamma(s)|^2 }-\frac{1}{w^2} &=\frac{w^2-|\gamma(s+w)-\gamma(s)|^2 }{w^2\, |\gamma(s+w)-\gamma(s)|^2 }\\[2mm]
			&\geqslant \frac{w^2-|\gamma(s+w)-\gamma(s)|^2 }{w^4 }\\[2mm]
			&\geqslant \frac{1}{2w^4} \int_0^w\int_0^w |\gamma'(s+t_1)-\gamma'(s+t_2)|^2\, \dd t_1\, \dd t_2.
		\end{align*}
		Integrating first over $s\in\R$ and $w>0$, with $\tau=t_1-t_2$ and a translation in $s$, we obtain
		\begin{align*}
			\Mob(\gamma) & \geqslant \int_{\R}\int_0^\infty \frac{1}{2w^4} \int_0^w\int_0^w |\gamma'(s+t_1)-\gamma'(s+t_2)|^2\, \dd t_1\, \dd t_2\, \dd w\, \dd s\\[2mm]
			& \geqslant \int_{\R}\int_0^\infty \frac{1}{2w^4} \int_0^w\int_0^w |\gamma'(x+t_1-t_2)-\gamma'(x)|^2\, \dd t_1\, \dd t_2\, \dd w\, \dd x\\[2mm]
			& \geqslant \int_{\R}\int_0^\infty \frac{1}{2w^4} \int_0^w\int_{-t_2}^{w-t_2} |\gamma'(x+\tau)-\gamma'(x)|^2\, \dd \tau\, \dd t_2\, \dd w\, \dd x\\[2mm]
			& \geqslant \int_{\R} \int_{\R} |\gamma'(x+\tau)-\gamma'(x)|^2  \int_0^\infty\int_0^w \frac{1}{2w^4}\, \mathbf{1}_{\{-t_2\leq\tau\leqslant w-t_2\}}\, \dd t_2\, \dd w\, \dd\tau \dd x\\[2mm]
			& \geqslant \int_{\R} \int_{\R} |\gamma'(x+\tau)-\gamma'(x)|^2  \int_{|\tau|}^\infty \frac{1}{2w^4} \left( \int_{\max(0,-\tau)}^{\min(w,w-\tau)} \dd t_2\right) \, \dd w\, \dd\tau \dd x\\[2mm]
			& \geqslant  \frac{1}{2}\, \int_{\R}\int_{\R} |\gamma'(x+\tau)-\gamma'(x)|^2 \left(\int_{|\tau|}^{\infty}\frac{w-|\tau|}{w^4}\dd w\right) \dd x\, \dd \tau.
		\end{align*}
		We compute the inner integral as follows
		\begin{equation*}
			\int_{|\tau|}^{\infty}\frac{w-|\tau|}{w^4}\dd w =\int_{|\tau|}^\infty \frac{ \dd w }{w^3} -|\tau|\, \int_{|\tau|}^\infty \frac{ \dd w }{w^4} =\frac{1}{2\, \tau^2}-\frac{1}{3\, \tau ^2} =\frac{1}{6\, \tau^2}.
		\end{equation*}
		We end up with 
		\begin{equation*}
		\Mob(\gamma) \geqslant c \int_{\R^2}\frac{|\gamma'(x)-\gamma'(y)|^2}{|x-y|^2}\, \dd x\, \dd y.
		\end{equation*}
		The cases $I=\R/(\ell\Z)$ are similar.
	\end{proof}

	\subsection{The Loewner energy}\label{sec:Loewner}

Given a Jordan curve $\Gamma$ in $\mathbb C$ (namely, $\infty \notin \Gamma$), let $\Omega$ and $\Omega^*$ denote respectively the bounded and unbounded connected components of $\mathbb C \setminus \Gamma$, the \emph{Loewner energy} (or called \emph{universal Liouville action}) of $\Gamma$ is given in \eqref{eq:def_IL} using two conformal maps $f \colon \mathbb D \to \Omega$ and $g\colon \mathbb D^* \to \Omega^*$ such that $g(\infty) = \infty$ as
\begin{equation}\label{eq:loewner_1}
		\IL(\Gamma) \coloneqq \frac{1}{\pi}\int_{\D} \left| \frac{f''}{f'} (z)\right|^2  |\dd z|^2+ \frac{1}{\pi}\int_{\D^*} \left|\frac{g''}{g'} (z)\right|^2 |\dd z|^2 + 4\, \log\frac{ |f'(0)| }{|g'(\infty)|}.
	\end{equation}
If $\infty \in \Gamma$, then we let $H$ and $H'$ denote the two connected components of $\mathbb C \setminus \Gamma$, $\mathbb H$ and $\mathbb H^*$ be the upper and lower half-planes,  the conformal maps $f \colon \mathbb H \to H$ and $g\colon \mathbb H^* \to H^*$ such that $f (\infty) = \infty$, $g(\infty) = \infty$. The Loewner energy is defined in this case as
\begin{equation}\label{eq:loewner_2}
		\IL(\Gamma) \coloneqq \frac{1}{\pi}\int_{\mathbb H} \left| \frac{f''}{f'} (z)\right|^2  |\dd z|^2+ \frac{1}{\pi}\int_{\mathbb H^*} \left|\frac{g''}{g'} (z)\right|^2 |\dd z|^2.
	\end{equation}
Although it is not obvious from these expressions, the definition of the Loewner energy does not depend on the choice of $f$ and $g$, and is also invariant if we replace $\Gamma$ by $A \circ \Gamma$, where $A \in \operatorname{PSL}(2,\mathbb C)$ is a M\"obius transformation of $\Chat$. In particular, if $A$ sends a point on a bounded curve $\Gamma$ to $\infty$, the Loewner energy remains unchanged although its expression changes (see Theorem~\ref{thm:equiv} below).
It was proved in \cite{TT2006,B2025} that $\IL (\Gamma) < \infty$ if and only if $\Gamma$ is a Weil--Petersson Jordan curve. 


There is another equivalent definition of $\IL$ using the loop Loewner driving function, which explains the origin of the name ``Loewner energy'' \cite{FS2017,RW2021,W2019}. 
It was proved to be the same as the expression \eqref{eq:def_IL} in \cite{W2019-2}. 
As the Loewner driving function will be instrumental to our discussion about the regularity of Loewner energy optimizers, we give a detailed definition here. Given a curve $\Gamma \subset \Chat$, we first choose
\begin{itemize}
    \item a point on $\Gamma$ that we call the \emph{root} and denote as $\Gamma(-\infty) = \Gamma (+\infty)$;
    \item an orientation of $\Gamma$;
    \item another point on $\Gamma$ that we denote as $\Gamma (0)$;
    \item a conformal map $h_0$ sending $\Chat \setminus \Gamma [-\infty, 0] \to \mathbb H$ (the upper halfplane) such that $h_0(\Gamma(-\infty)) = \infty$ and $h_0 (\Gamma (0)) = 0$, $\Gamma [-\infty, 0]$ is the subarc of $\Gamma$ from the root to $\Gamma(0)$ along the orientation of $\Gamma$. 
\end{itemize}
With these choices fixed, then for every subarc $A$ of $\Gamma$ starting from the root, there exists a unique conformal map $h_A \colon \Chat \setminus A \to \mathbb H$, sending the root to $\infty$ and the other end of $A$ to $0$, such that the expansion at $\infty$ of 
\begin{equation}
    \label{eq:hydro}
h_A \circ h_0^{-1} (z) = z - W_A + \frac{C_A}{z} + o \left(\frac{1}{z}\right).
\end{equation}
The constant $C_A \in \mathbb R$ is the \emph{relative halfplane capacity} of $A$, which is positive if $A$ includes $\Gamma [-\infty, 0]$, and $C_A \le 0$ if $A$ is included in $\Gamma [-\infty, 0]$, and $A \mapsto C_A$ is increasing for the inclusion.  We now parametrize $\Gamma$ \emph{by capacity}, in the way such that 
$$C_{\Gamma [-\infty, t]} = 2t.$$
This capacity parameterization depends on all the choices in the bullet points that we made before, but is otherwise determined by $\Gamma$, and compatible with the notation of $\Gamma (0)$ and $\Gamma (\pm \infty)$ above. 
Once we know the capacity parameterization, we write $h_t \coloneqq h_{\Gamma [-\infty, t]}$, and $W(t) \coloneqq W_{\Gamma [-\infty, t]}$.  In particular, $W(0) = 0$.
\begin{defi}
    The function $W \colon \mathbb R \to \mathbb R$ is called the \emph{loop Loewner driving function} of $\Gamma$.
\end{defi}

 Although $W$ is not canonically associated with $\Gamma$ since it depends on the choices in the bullet points, the Dirichlet energy of $W$ is, as shown by the next theorem. 

\begin{theorem}[\cite{W2019-2}] \label{thm:equiv}
   The Dirichlet energy of a driving function $W$ of $\Gamma$ coincides with $\IL (\Gamma)$ defined in \eqref{eq:loewner_1} and $\eqref{eq:loewner_2}$: If $W$ is absolutely continuous, then
\begin{equation*}
    \IL(\Gamma) = \frac12 \int_{-\infty}^{\infty} \left(\frac{\dd W(t)}{\dd t}\right)^2 \, \dd t
\end{equation*}
and $\IL(\Gamma) = \infty$ otherwise. In particular, the Dirichlet energy of $W$ is finite if and only if $\Gamma$ is Weil--Petersson.
\end{theorem}

The driving function description of the Loewner energy also allows us to describe the partial energy of a subarc. 

Let $A$ be an unparametrized arc, that is, the image of a continuous embedding of $[0,1]$ in $\Chat$. We define its driving function by choosing one end as root $A(-\infty)$, the other end as $A(0)$, and a conformal map $h_0 \colon \Chat \setminus A \to \mathbb H$ as above. It then provides in a similar way a capacity parametrization of $A \colon [-\infty, 0] \to \Chat$ of the arc, and a driving function $W \colon (-\infty, 0] \to \mathbb R$.
\begin{defi}
   The \emph{arc Loewner energy} of unparametrized arc $A$ is  
   \begin{equation*}
    {\rm I^A} (A) = \frac12 \int_{-\infty}^{0} \left(\frac{\dd W(t)}{\dd t}\right)^2 \, \dd t.
\end{equation*}
  We say that an arc  $ A \subset \Chat$ is a \emph{Weil--Petersson} arc if it is a subarc of a Weil--Petersson Jordan curve, which is equivalent to ${\rm I^A} (A) < \infty$.
\end{defi}

In other words, the arc energy of $A$ is naturally viewed as a part of the loop energy of $\Gamma$ if $A = \Gamma [-\infty, 0]$, and the driving function of $A$ is the restriction of the loop driving function $W$ of $\Gamma$ to $(-\infty, 0]$.

One may wonder if $\frac12 \int_{0}^{\infty} \left(\frac{\dd W(t)}{\dd t}\right)^2 \, \dd t$ also means something. The first thing one notices is that the restriction of $W$ on $[0,\infty)$ depends not only on $\Gamma [0,\infty]$, but also on $\Gamma [-\infty, 0]$. So this integral should describe $\Gamma [0,\infty]$ with respect to $\Gamma [-\infty, 0]$. 

The correct notion is the \emph{chordal Loewner energy} in a simply connected domain. 
Consider a connected and simply connected open set $D$ with distinct prime ends $a, b$. Let $C$ be a \emph{chord}, namely a simple, unparametrized oriented curve in $D$ going from $a$ to $b$.  
We fix a conformal map $h_0$ from $D \to \mathbb H$ sending $a$ to $0$ and $b$ to $\infty$. For each initial arc $A$ of $C$ (starting from $a$), there is a unique map $h_A$ satisfying \eqref{eq:hydro}. Now, the capacities are all nonnegative, and we parameterize $C$ by capacity with $C(0) = a$, and obtain a driving function $W \colon [0,\infty) \to \mathbb R$ the same way as in \eqref{eq:hydro}. 

\begin{defi}
    The chordal Loewner energy of a chord $C$ in $(D;a,b)$ is 
    $$
    {\rm I}_{D;a,b} (C) \coloneqq \frac12 \int_{0}^{\infty} \left(\frac{\dd W(t)}{\dd t}\right)^2 \, \dd t$$
    and it does not depend on the choice of the conformal map $h_0$.
\end{defi}
In fact, we can replace $h_0$ by $\lambda h_0$  with $\lambda >0$ (and this is the only ambiguity in the definition of the chordal driving function), then the driving function is changed to $t \mapsto \lambda W(\lambda^{-2} t)$, which has the same Dirichlet energy as $W$.

 It follows immediately from the definition of the chordal driving function:
\begin{lemma}\label{lem:minI}
    If $\varphi \colon D \to D'$ is a conformal map, then for any chord $C$, 
    $${\rm I}_{D;a,b} (C) = {\rm I}_{D'; \varphi(a),\varphi(b)} (\varphi(C)).$$
   In particular, ${\rm I}_{D;a,b} (C) = 0$ if and only if $W \equiv 0$, which is equivalent to $h_0 (C) = \ii \mathbb R_+$. In other words, $C$ is the hyperbolic geodesic in $D$ connecting $a$ to $b$.
    \end{lemma}

    By comparing to the loop driving function, we obtain:
     \begin{lemma}\label{lem:add}
 If $\Gamma$ is a Jordan curve and $W \colon \mathbb R \to \mathbb R$ is its driving function, then
    $${\rm I}_{\Chat \setminus \Gamma [-\infty, 0];\Gamma(0),\Gamma(\infty)} (\Gamma [0,\infty]) =  \frac12 \int_{0}^{\infty} \left(\frac{\dd W(t)}{\dd t}\right)^2 \, \dd t.$$
    In other words, if $\Gamma$ is a Jordan curve in $\Chat$ and decomposed into two arcs $\Gamma_0$ and $\Gamma_1$, such that $\Gamma_0 \cap \Gamma_1 = \{a,b\}$, then 
    $$\IL (\Gamma) = {\rm I^A} (\Gamma_0) + {\rm I}_{\Chat \setminus \Gamma_0; a,b } (\Gamma_1).$$
\end{lemma}

\begin{remark}
    Finally, we emphasize that all versions of the Loewner energy are associated with unparametrized curves, although the definition of the driving function starts from determining the capacity parametrization by fixing the root, the conformal map $h_0$, etc.; the Dirichlet energy of the driving function does not depend on any of these choices.
\end{remark}



	\subsection{Optimizations of the Loewner energy}
	\label{sec:loewner_arc_completion}
	
	In this section, we recall results on $\Pc^1_{\IL}$ and $\Pc^2_{\IL}$ and prove the regularity of the solution to the latter. Let us start with the latter. 

    \begin{theorem}\label{thm:loewner_arc_completion}
		Let $\Gamma_0$ be a $\Cd^{1,\alpha}$ arc \textnormal{(}in arclength parameterization\textnormal{)} with endpoints $a$ and $b$, where $\alpha > 1/2$. The unique chord $\Gamma_1$ in $\Chat \setminus \Gamma_0$ connecting $a$ and $b$ that minimizes the Loewner energy of $\Gamma = \Gamma_0 \cup \Gamma_1$ is the hyperbolic geodesic in $\Chat \setminus \Gamma_0$. Moreover, $\Gamma$ is $\Cd^{1,1/2}$ regular.
	\end{theorem}

\begin{proof}
By Lemma~\ref{lem:minI} and Lemma~\ref{lem:add}, the curve $\Gamma_1$ minimizing the energy of $\Gamma_0 \cup \Gamma_1$ is given by the hyperbolic geodesic in $\Chat \setminus \Gamma_0$ connecting $a$ and $b$. 

In particular, the interior of $\Gamma_1$ is analytic, hence the regularity of $\Gamma$ is determined by that at the junctions $a$ and $b$. Let us consider the regularity at $a$ and root the curve $\Gamma$ at $b$ (the role of $a$ and $b$ are exactly symmetric), so that the driving function of $\Gamma_0$ is a function $W \colon (-\infty, 0] \to \mathbb R$, where $W(0)$ corresponds to the point $a$. 

    It was shown in Theorem 1.3 of \cite{RW2021} that the H\"older regularity of $W$ drops by $1/2$ compared to that of $\Gamma_0$, namely, it is $\Cd^{1/2+ \alpha}$ up to time $0$. The driving function of $\Gamma$ is that of $W$ completed to a function $\widetilde W \colon \mathbb R \to \mathbb R$ such that $\widetilde W$ is constant on $[0,\infty)$. In particular, $\widetilde W$ is Lipschitz in a neighborhood of $0$. Theorem 4.8 of \cite{W2014} then shows that $\Gamma$ is $\Cd^{1,1/2}$ in the neighborhood of $a$ which completes the proof.
\end{proof}

\begin{remark}
    We also know that the solution $\Gamma$ is not $\Cd^{1,\alpha'}$ for any $\alpha' > 1/2$, as it would imply that its driving function is $\Cd^{1,\alpha'-1/2}$, which is not true if $W'(0-) \neq 0$.
\end{remark}

Now we turn to $\Pc^1_{\IL}$. Given $n\geqslant 2$ distinct points $Z\coloneqq (z_1,\ldots,z_n)\in \C^n$ and a relative isotopy class $\tau$ of curves passing through $Z$ in its order, we seek to understand $\Pc^1_{\IL}(Z,\tau)$ defined in Equation \eqref{eq:Pc1-IL}.
According to Lemma~\ref{lem:minI} and Lemma~\ref{lem:add}, if such a minimizer $\gamma_{Z,\tau}\colon \T \to \Chat$ exists, then it has the \emph{piecewise geodesic property}: If $a_i\in\T$ is such that $\gamma_{Z,\tau}(a_i)=z_i$, then $\gamma_{Z,\tau}|_{(a_i,a_{i+1})}$ parametrizes a hyperbolic geodesic in the simply connected domain $\Chat \setminus \gamma_{Z,\tau}(\T\setminus (a_i,a_{i+1}))$. The existence of the minimizer follows from a simple compactness result, but the uniqueness is much harder and was proved in \cite{BJRW2025} (see also \cite{W2025} which reinterpreted the welding of the minimizer as the boundary of a pleated plane in Anti-de Sitter $3$-space).

\begin{theorem}[\cite{BJRW2025}]
		For any distinct points $Z= (z_1,\ldots,z_n)\in\C^n$ and relative isotopy class $\tau$,  there exists a unique minimizer $\gamma_{Z,\tau}$ of the Loewner energy. 
	\end{theorem}
In terms of regularity, the piecewise geodesic property implies that $\gamma_{Z,\tau}$ is analytic outside of $Z$. The regularity at each point in $Z$ is reduced to the study of the regularity of a finite energy \emph{geodesic pair}, namely, a chord in a simply connected domain $D$ of the form $\Gamma_0 \cup \Gamma_1$, where $\Gamma_0$ and $\Gamma_1$ meet at a point $a$ and such that $\Gamma_0$ is a hyperbolic geodesic in $D \setminus \Gamma_1$ and $\Gamma_1$ is a hyperbolic geodesic in $D \setminus \Gamma_0$. 
Lemma 2.19 in \cite{BJRW2025}  (based on an explicit computation in Eq.(1) of \cite{MRW2026}) shows that a finite energy geodesic pair $\gamma \colon (-L,L') \to \Gamma_0 \cup \Gamma_1$  parametrized by arclength is continuously differentiable, and its derivative is log-Lipschitz: 
$$\text{if }\gamma (0) = a, \qquad |\gamma'(t) - \gamma'(0)| \leqslant C |t \log |t||$$
for some constant $C$.

\section{Energies and Fermionic Grunsky operator}

	\subsection{Möbius energy and Cauchy operator}\label{sec:Mobius-Cauchy}

	In this section, we express the Möbius energy of a Weil--Petersson curve as the Hilbert--Schmidt norm of an operator in Proposition \ref{pr:Mob-Cauchy}. For the rest of the section, we will use the following notations. Let $\Gamma=\gamma(\T)\subset \C$ be a Weil--Petersson curve and $\tau(x)\coloneqq \frac{\gamma'}{|\gamma'|}\circ\gamma^{-1}\colon \Gamma\to\s^1$ the tangent vector. 
	We define the \emph{Cauchy operator} for $u\in \Cd^{\infty}(\Gamma;\C)$ as
	\begin{equation}\label{eq:def-Cr}
		\Cr_{\Gamma}u(x) \coloneqq \frac{u(x)}{2} + \frac{1}{2\ii \pi} \pv \int_{\Gamma} \frac{u(y)}{y-x}\, \tau(y)\, |\dd y|.
	\end{equation}
	We will also use the following notation:
	\begin{equation}\label{eq:def-Cc}
		\Cc_{\Gamma} u(x) \coloneqq \frac{1}{2\ii \pi} \pv \int_{\Gamma} \frac{u(y)}{y-x}\, \tau(y)\, |\dd y|.
	\end{equation}
	The main result of this section is the following.
    We denote 
    $$
    \Ll^2(\Gamma) \coloneqq \left\{u \colon \Gamma \to \C \mid \int_\Gamma| u(z)|^2 \, |\dd z| < \infty\right\}.
    $$
	\begin{proposition}\label{pr:Mob-Cauchy}
		Let $\Gamma\subset \C$ be a closed Weil--Petersson curve. Then it holds
		\begin{equation*}
			\Mob(\Gamma)-4 = 2\pi^2\, \|\Cr_{\Gamma}-\Cr_{\Gamma}^*\|^2_{\rm HS(L^2(\Gamma))}=  2\pi^2\, \|\Cc_{\Gamma}-\Cc_{\Gamma}^*\|_{\rm HS(L^2(\Gamma))}^2 .
		\end{equation*}
	\end{proposition}
	
    We first express the Möbius energy using the conformal angle of $\Gamma$, defined as follows; see, for instance, \cite[Definition 3.2]{ohara2003}. 
	\begin{defi}\label{def:conformal-angle}
		Let $x,y\in\Gamma$ be distinct. Let $S_{x,y}$ be the circle tangent to $\Gamma$ at $x$ and passing through $y$ oriented as the curve $\Gamma$ at $x$. The \emph{conformal angle} $\theta_{\Gamma}(x,y)\in[0,\pi]$ is the angle between $S_{x,y}$ and $S_{y,x}$. 
	\end{defi}

    Note that $\theta_{\Gamma} \equiv 0$ if and only if $\Gamma$ is a circle.

\begin{figure}[ht]
    \begin{center}    \includegraphics[width=0.45\textwidth]{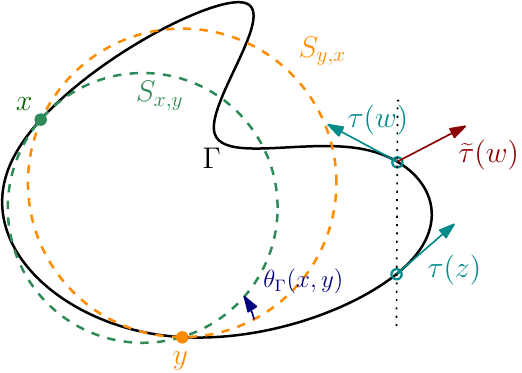}
    \end{center}
    \caption{Illustration of the conformal angle $\theta_\Gamma (x,y)$ in Definition~\ref{def:conformal-angle}, and $\tau (z)$ and $\widetilde \tau(w)$ in the proof of Proposition~\ref{pr:Mob-Cauchy}.}
    \label{fig::circles}
\end{figure}
	
	\begin{proof}[Proof of Proposition \ref{pr:Mob-Cauchy}]
	On the one hand, by \cite[Equations (3.5) and (3.6)]{ohara2003}, the Möbius energy has the following expression
	\begin{equation}\label{eq:cos-law}
		\Mob(\Gamma)-4 = \int_{\Gamma\times \Gamma} \frac{1-\cos\theta_{\Gamma}(x,y)}{|x-y|^2}\, |\dd x|\, |\dd y|.
	\end{equation}
	On the other hand, the kernel of $\Cr_{\Gamma} - \Cr_{\Gamma}^*=\Cc_{\Gamma}-\Cc_{\Gamma}^*$, defined in \eqref{eq:def-Cr} and \eqref{eq:def-Cc}, is given by Equation (7.1) in \cite{KS1978} as follows
	\begin{equation}\label{eq:kernel-a}
		\forall w,z\in\Gamma,\qquad a_{\Gamma}(w,z) \coloneqq  \frac{1}{2\ii \pi}\, \frac{\widetilde{\tau}(w) - \tau(z)}{z-w},
	\end{equation}
	where $\widetilde{\tau}(w)$ is the Euclidean reflection of $\tau(w)$ by the chord determined by $z$ and $w$. Therefore, it holds
	\begin{equation*}
		|a_{\Gamma}(w,z)|^2 = \frac{1}{2\pi^2}\, \frac{1 - \scal{\widetilde{\tau}(w)}{\tau(z)}_{\C} }{|z-w|^2}
	\end{equation*}
	We now claim that the conformal angle $\theta_{\Gamma}(z,w)$ verifies 
	\begin{equation}\label{eq:costheta}
		\cos\theta_{\Gamma}(z,w) = \scal{\tau(w)}{\widetilde{\tau}(z)}_{\C}.
	\end{equation}
	Indeed, we consider two points $z$ and $w$ on a circle, say $z=\ee^{\ii s},w=\ee^{\ii t}\in\s^1$. The tangent vectors to $\s^1$ at $z$ and $w$ are respectively $\ii z=\ee^{\ii \left(s+\pi/2 \right)}$ and $\ii w = \ee^{\ii\left(t+\pi/2\right)}$. Therefore, we have 
		\begin{equation*}
			z-w = \ee^{\ii\frac{s+t}{2}} \left[ \ee^{\ii \frac{s-t}{2}} - \ee^{\ii \frac{t-s}{2}} \right] = 2\ii \ee^{\ii \frac{s+t}{2}}\, \sin\left(\frac{s-t}{2}\right).
		\end{equation*}
		Hence, there holds
		\begin{equation}\label{eq:arguments}
			2\, \arg(z-w) = s+t+\pi = \arg(\ii z) + \arg(\ii w).
		\end{equation}
		If $S_{x,y}\subset \C$ is the circle tangent to $\Gamma$ at $x$ and passing through $y$, then we deduce from \eqref{eq:arguments} that the tangent vector $T_{S_{x,y}}(y)$ to $S_{x,y}$ at $y$ verifies 
		\begin{equation*}
			2\, \arg(y-x) = \arg T_{S_{x,y}}(y) +  \arg\tau(x).
		\end{equation*}
		Therefore, we obtain 
		\begin{equation}\label{eq:theta_tau}
			\theta_{\Gamma}(y,x) = \angle (S_{y,x},S_{x,y}) = \arg T_{S_{y,x}}(y) - \arg T_{S_{x,y}}(y) = \arg \tau(y) - 2\, \arg(y-x) +  \arg\tau(x).
		\end{equation}
		Moreover, the reflection of $\tau(y)$ with respect to the segment $y-x$ change the argument to $ 2\, \arg(y-x)-\arg \tau(y)$. Therefore, we obtain 
		\begin{equation*}
			\scal{\tau(x)}{\widetilde{\tau}(y)}_{\C} = \cos\Big( \arg\tau(x) - \big( 2\, \arg(y-x)-\arg \tau(y) \big) \Big) = \cos\theta_{\Gamma}(y,x).
		\end{equation*}
	This proves \eqref{eq:costheta}.
	
	As a consequence of \eqref{eq:cos-law} and \eqref{eq:costheta} (see also \cite[Equation (10)]{BB2007}), we have
	\begin{equation*}
		\|\Cr_{\Gamma}-\Cr_{\Gamma}^*\|_{\rm HS(L^2(\Gamma))}^2 = \int_{\Gamma\times \Gamma} |a_{\Gamma}(w,z)|^2 \ |\dd w|\, |\dd z| = \frac{1}{2\pi^2}\, \left(\Mob(\Gamma)-4\right).
	\end{equation*}
	This concludes the proof of Proposition \ref{pr:Mob-Cauchy}.
	\end{proof}
	
	Thanks to Proposition \ref{pr:Mob-Cauchy}, we now derive a relation between the Möbius energy and the $\W^{1/2,2}$-seminorm of the tangent vector.
	
	\begin{proposition}
		Let $\Gamma\subset \C$ be a Weil--Petersson curve. We denote by $\tau$ its tangent vector and by $n\coloneqq \ii\tau$ its inward-pointing normal vector. We have
		\begin{equation}\label{eq:Mobius-tangent} 
			\Mob(\Gamma)-4  = \frac{1}{2}\, [\tau]_{\W^{1/2,2}(\Gamma)}^2 + 2 \int_{\Gamma\times \Gamma} \frac{\scal{x-y}{n(x)}_{\C} \, \scal{x-y}{n(y)}_{\C}}{|x-y|^4}\, |\dd x|\,|\dd y|.
		\end{equation}
	\end{proposition}
    \begin{remark}
    \begin{itemize}
        \item This decomposition coincides with the one obtained in \cite[Theorem 2.1]{IN2014}.
        \item The second term of the right-hand side in \eqref{eq:Mobius-tangent} is related to the Hilbert--Schmidt norm of the standard \textnormal{(}or ``bosonic''\textnormal{)} Grunsky operator, see Lemma 6.2 in \cite{M2026}.
    \end{itemize}
    \end{remark}
    
	\begin{proof} 
    For $a, b,x,y\in \R$, there holds
    \begin{equation*}
        \left\langle \begin{pmatrix} a\\ b \end{pmatrix}, \begin{pmatrix} x\\ y \end{pmatrix} \right\rangle_{\R^2} = a x + b y = \Re\left[ (a+\ii b)\, \overline{(x+\ii y)} \right].
    \end{equation*}
	For $f\in \Cd^{\infty}(\Gamma;\C)$, we have 
	\begin{equation*}
		\begin{aligned} 
		2\ii\pi\, \Cc_{\Gamma} f(x) & = \pv \int_{\Gamma} \frac{\tau(y)}{y-x}\, f(y)\, |\dd y|\\[2mm]
        & = \pv \int_{\Gamma} \frac{\overline{ (y-x) }\, \tau(y)}{|x-y|^2}\, f(y)\, |\dd y| \\[2mm]
		& = \pv \int_{\Gamma} \frac{\Re\left[ \overline{ (y-x) }\, \tau(y)\right] + \ii \Im\left[ \overline{ (y-x) }\, \tau(y)\right]}{|x-y|^2}\, f(y)\, |\dd y|\\[2mm]
		& = \pv \int_{\Gamma} \frac{ \scal{y-x}{\tau(y)}_{\C} - \ii \Re\left[ \overline{ (y-x) }\, \ii\tau(y)\right]}{|x-y|^2}\, f(y)\, |\dd y|\\[2mm]
		& = \pv \int_{\Gamma} \frac{ \scal{y-x}{\tau(y)}_{\C} - \ii \scal{y-x}{n(y)}_{\C}}{|x-y|^2}\, f(y)\, |\dd y|.
		\end{aligned} 
	\end{equation*}
	We denote
	\begin{equation} \label{eq:TN}
		\begin{cases}
			\displaystyle T_{\Gamma} f(x) \coloneqq \pv \int_{\Gamma} \frac{ \scal{y-x}{\tau(y)}_{\C}}{|x-y|^2}\, f(y)\, |\dd y|, \\[5mm]
			\displaystyle N_{\Gamma}f(x) \coloneqq \pv \int_{\Gamma} \frac{ \scal{x-y}{n(y)}_{\C}}{|x-y|^2}\, f(y)\, |\dd y|.
		\end{cases} 
	\end{equation}
	We have
	\begin{equation*}
		\begin{aligned} 
			\Cc_{\Gamma} - \Cc_{\Gamma}^* & = \frac{1}{2\ii \pi}\left(T_{\Gamma} + \ii N_{\Gamma}\right) + \frac{1}{2\ii \pi}\left(T_{\Gamma}^* - \ii N_{\Gamma}^*\right) \\[2mm]
			& = \frac{1}{2\ii \pi}\left( T_{\Gamma}+T_{\Gamma}^*\right) + \frac{1}{2\pi}\left(N_{\Gamma} - N_{\Gamma}^*\right).
		\end{aligned} 
	\end{equation*}
	Since the kernels of $T_{\Gamma}$ and $N_{\Gamma}$ are real-valued, we obtain by Proposition \ref{pr:Mob-Cauchy}
	\begin{equation}\label{eq:Mobius-NP}
		\frac{1}{2\pi^2}\left(\Mob(\Gamma)-4\right) = \|\Cc_{\Gamma} - \Cc_{\Gamma}^*\|_{\rm HS}^2 = \frac{1}{4\pi^2}\left(  \|T_{\Gamma}+T_{\Gamma}^*\|_{\rm HS}^2 +\|N_{\Gamma}-N_{\Gamma}^*\|_{\rm HS}^2 \right).
	\end{equation}	
	Each term is given by
	\begin{equation*}
		\begin{aligned}
			\|N_{\Gamma}-N_{\Gamma}^*\|_{\rm HS}^2 & = \int_{\Gamma\times \Gamma} \frac{\left|\scal{x-y}{n(x)+n(y)}_{\C}\right|^2}{|x-y|^4}\, |\dd x|\, |\dd y|, \\[2mm]
			\|T_{\Gamma}+T_{\Gamma}^*\|_{\rm HS}^2 & = \int_{\Gamma\times \Gamma} \frac{\left|\scal{x-y}{\tau(x)-\tau(y)}_{\C}\right|^2}{|x-y|^4}\, |\dd x|\, |\dd y|.
		\end{aligned}
	\end{equation*}
	Moreover, it holds
	\begin{equation*}
		\begin{aligned} 
			\left|\scal{x-y}{\tau(x)-\tau(y)}_{\C}\right|^2 & = \left| \Re\left[ \overline{(x-y)}\, (\tau(x)-\tau(y)) \right] \right|^2 \\[2mm]
			& = \left| \overline{(x-y)}\, (\tau(x)-\tau(y)) \right|^2 - \left| \Im\left[ \overline{(x-y)}\, (\tau(x)-\tau(y)) \right] \right|^2 \\[2mm]
			& = |x-y|^2\, |\tau(x)-\tau(y)|^2 - \left| \Re\left[ \overline{(x-y)}\, \ii (\tau(x)-\tau(y)) \right] \right|^2 \\[2mm]
			& = |x-y|^2\, |\tau(x)-\tau(y)|^2 - \left| \scal{x-y}{n(x)-n(y)}_{\C} \right|^2 \\[2mm]
			& = |x-y|^2\, |\tau(x)-\tau(y)|^2 - \left| \scal{x-y}{n(x)}_{\C} \right|^2 - \left| \scal{x-y}{n(y)}_{\C} \right|^2 \\[2mm]
            & \quad + 2\, \scal{x-y}{n(x)}_{\C} \, \scal{x-y}{n(y)}_{\C}.
		\end{aligned} 
	\end{equation*}
	Hence, we obtain 
	\begin{equation*}
		\begin{aligned}
			& \|T_{\Gamma}+T_{\Gamma}^*\|_{\rm HS}^2 +\|N_{\Gamma}-N_{\Gamma}^*\|_{\rm HS}^2 \\[2mm]
            & = \int_{\Gamma\times \Gamma} \frac{|\tau(x)-\tau(y)|^2}{|x-y|^2}\, |\dd x|\, |\dd y|+ 4 \int_{\Gamma\times \Gamma} \frac{\scal{x-y}{n(x)}_{\C} \, \scal{x-y}{n(y)}_{\C}}{|x-y|^4}\, |\dd x|\,|\dd y|.
		\end{aligned}
	\end{equation*}
	We end up with 
	\begin{equation*}
		\begin{aligned} 
			\Mob(\Gamma)-4 & = 2\pi^2 \|\Cc_{\Gamma} - \Cc_{\Gamma}^*\|_{\rm HS}^2 \\[2mm]
            &  = \frac{1}{2}\left(  \|T_{\Gamma}+T_{\Gamma}^*\|_{\rm HS}^2 +\|N_{\Gamma}-N_{\Gamma}^*\|_{\rm HS}^2 \right) \\[2mm]
			& = \frac{1}{2}\, [\tau]_{\W^{1/2,2}(\Gamma)}^2 + 2\int_{\Gamma\times \Gamma} \frac{\scal{x-y}{n(x)}_{\C} \, \scal{x-y}{n(y)}_{\C}}{|x-y|^4}\, |\dd x|\,|\dd y|.
		\end{aligned} 
	\end{equation*}
	\end{proof}

	\subsection{Fermionic Grunsky operator}\label{sec:Grunksy}
	
	In this section, we discuss the notion of half-order differential forms and the fermionic Grunsky operator introduced in \cite{KSS2025}.

	\subsubsection{Half-order differential forms}

    We follow closely the presentation of \cite{KSS2025}. 
Given two domains (i.e., non-empty, connected, open subset) $\Omega_1,\Omega_2\subset \Chat$ and assume there is a biholomorphism $\sigma\colon \Omega_1\to \Omega_2$, then there are two choices of 
  meromorphic functions $\sqrt{\sigma'} \colon\Omega_1\to \C$ such that $
    (\sqrt{\sigma'})^2 = \sigma'$.  The function $\sqrt{\sigma'}$ has at most one pole of order $1$ in case there is a
    point $z$ where $\sigma (z) = \infty$.
    
    Given a choice of $\sqrt{\sigma'}$, we use $\widehat{\sigma}$ for the short-hand notation for the pair $\left(\sigma,\sqrt{\sigma'} \right)$. We denote the set of all such choices by 
	\begin{equation*}
		\Confhat(\Omega_1,\Omega_2) \coloneqq \left\{ \widehat \sigma = \left( \sigma,\sqrt{\sigma'} \right) \mid \sigma \colon \Omega_1\to \Omega_2 \text{ biholomorphic and }\sqrt{\sigma'}^2 = \sigma' \right\}. 
	\end{equation*}
    \begin{defi}
        Let $\Omega_1$ be a domain of  $\Chat$.
      If $\infty \notin \Omega_1$,
        the space of \emph{holomorphic half-order differential forms} on $\Omega_1$ is the space of holomorphic functions $h\colon \Omega_1 \to \mathbb C$, denoted as $h(z)\, \dd z^{1/2}$. 
          If $\infty \in \Omega_1$, then we require that $h(z)=O(|z|^{-1})$ as $|z|\to\infty$.
        
        If $\Omega_2$ is another domain and there exists $\widehat \sigma \in \Confhat(\Omega_2,\Omega_1)$,  the \emph{pullback} of $h(z)\, \dd z^{1/2}$ by $\widehat{\sigma}\in\Confhat(\Omega_2,\Omega_1)$ is a holomorphic half-order differential form on $\Omega_2$:
	\begin{equation*}
		\widehat{\sigma}^*\left(h(z)\, \dd z^{1/2}\right) \coloneqq h\circ \sigma(z)\, \sqrt{\sigma'(z)}\, \dd z^{1/2}.
	\end{equation*}
    
    	Similarly, an \emph{anti-holomorphic half-order differential form} on $\Omega_1$ is of the form $\overline h(z)\, \dd \bar{z}^{1/2}$, where $\overline h$ is an anti-holomorphic function $\Omega_1\to \C$ \textnormal{(}and $h(z)=O(|z|^{-1})$ if $\infty \in \Omega_1$\textnormal{)} and transforms like
	\begin{equation*}
		\widehat{\sigma}^*\left(\overline h(z)\, \dd \bar{z}^{1/2}\right) \coloneqq \overline h\circ \sigma (z)\, \overline{ \sqrt{\sigma'(z)} }\, \dd \bar{z}^{1/2}, \qquad \text{for } \widehat \sigma \in \Confhat(\Omega_2,\Omega_1).
	\end{equation*}
    \end{defi}

\begin{remark}
      The decay condition on $h$ when $\infty \in \Omega_1$ ensures that the pullback of $h(z) \, \dd z^{1/2}$ to a domain $\Omega_2$ not containing $\infty$ extends to a holomorphic half-order differential form.
\end{remark}

	For the rest of the section, we consider $\Gamma\subset \C$ a Weil--Petersson curve. We denote $\Omega$ the component of $\Chat \setminus \Gamma$ that does not contain $\infty$ and $\Omega^*$ the component that contains $\infty$. We consider two biholomorphisms $f\colon \D\to \Omega$ and $g\colon \Chat \setminus \overline{\D} \to \Omega^*$ that fixes $\infty$.

	The space of $\Ll^2$ holomorphic half-differential forms on $\Omega$ is defined as
	\begin{equation*}
		\Ar^{1/2}(\Omega) \coloneqq \left\{ h(z)\, \dd z^{1/2} : h \text{ holomorphic in }\Omega,\ \lim_{r\to 1^-}\int_{f(\s^1_r)} |h(z)|^2\, |\dd z|<\infty \right\}.
	\end{equation*}
    
	The space of $\Ll^2$ anti-holomorphic half-differential forms is 
	\begin{equation*}
		\overline{  \Ar^{1/2}(\Omega) } \coloneqq \left\{ \bar{h}\, \dd \bar{z}^{1/2} : h\, \dd z^{1/2}\in \Ar^{1/2}(\Omega) \right\}.
	\end{equation*}
	We consider the following scalar products on $ \Ar^{1/2}(\Omega)$ and  $\overline{  \Ar^{1/2}(\Omega)}$
	\begin{equation*}
		\scal{h_1\, \dd z^{1/2}}{h_2\, \dd z^{1/2}}_{\Omega} \coloneqq \lim_{r\to 1^-} \int_{f(\s^1_r)} h_1(z)\, \overline{ h_2(z) }\, |\dd z|,
	\end{equation*}
	\begin{equation*}
		\scal{\overline{h_1}\, \dd \bar{z}^{1/2}}{ \overline{h_2}\, \dd \bar{z}^{1/2}}_{\Omega} \coloneqq \lim_{r\to 1^-} \int_{f(\s^1_r)} \overline{ h_1(z) }\, h_2(z)\, |\dd z|,
	\end{equation*}
and $\Ar^{1/2}(\Omega)$ and  $\overline{  \Ar^{1/2}(\Omega)}$ are orthogonal:
\begin{equation*}
		\scal{h_1\, \dd \bar{z}^{1/2}}{ \overline{h_2}\, \dd \bar{z}^{1/2}}_{\Omega} \coloneqq \lim_{r\to 1^-} \int_{f(\s^1_r)} h_1(z) \, h_2(z)\, \dd z = 0.
	\end{equation*}
Moreover, these scalar products are invariant under pullback by $\widehat \sigma$.

	The direct sum of holomorphic half-differential forms and anti-holomorphic half-differential forms will be denoted
	\begin{equation*}
		\Ar^{1/2}_h(\Omega) \coloneqq \Ar^{1/2}(\Omega) \operp  \overline{\Ar^{1/2}(\Omega)}.
	\end{equation*}
	We denote the set of traces of (anti-)holomorphic half-differentials by
	\begin{equation*}
		W_{\Omega} \coloneqq \left\{ h|_{\Gamma} \in \Ll^2(\Gamma) : h\, \dd z^{1/2}\in \Ar^{1/2}(\Omega) \right\},
	\end{equation*}
	and
	\begin{equation*}
		\overline{W}_{\Omega} \coloneqq \left\{ \overline{H}|_{\Gamma} \in \Ll^2(\Gamma) : \overline{H}\, \dd \bar{z}^{1/2}\in \overline{\Ar^{1/2}(\Omega)} \right\}.
	\end{equation*}
	By \cite[Theorem 2.16]{KSS2025}, if $\tau \colon \Gamma \to \s^1 \subset \C$ denotes the tangent vector along $\Gamma$, then every function in $\Ll^2(\Gamma)$ can be written in the form of $h|_{\Gamma} + \overline{H|_{\Gamma}\, \tau}$. In other words, we have the decomposition
	\begin{equation}\label{eq:decompo_L2}
		\Ll^2(\Gamma) = W_{\Omega} \oplus (\overline{W}_{\Omega}\, \bar{\tau}). 
	\end{equation} 
	Moreover, the following map is an isometric isomorphism
	\begin{equation}\label{eq:def-b}
		\bb_{\Omega} \colon \left( h\, \dd z^{1/2} + \overline{H}\, \dd \bar{z}^{1/2} \right) \in \Ar^{1/2}_h(\Omega) \mapsto \left( h|_{\Gamma} + \overline{H}|_{\Gamma}\, \overline{\tau} \right) \in \Ll^2(\Gamma).
	\end{equation}
	The isomorphism $\bb_{\Omega^*}$ is defined similarly with $\Omega$ replaced by $\Omega^*$. 
	
	We will also use the following two operators: for $\alpha\in \Ar_h^{1/2}(\Omega)$,
	\begin{equation*}
		\forall z\in\Omega,\qquad \Jj^{1/2}\alpha(z) = \left( \frac{1}{2\ii \pi}\, \lim_{r\to 1^-} \int_{f(\s^1_r)} \alpha(w)\, \frac{\dd w^{1/2}}{w-z} \right)\, \dd z^{1/2},
	\end{equation*}
	and 
	\begin{equation*}
		\forall z\in\Omega^*,\qquad \Jj_*^{1/2}\alpha(z) = \left( \frac{1}{2\ii \pi}\, \lim_{r\to 1^-} \int_{f(\s^1_r)} \alpha(w)\, \frac{\dd w^{1/2}}{w-z} \right)\, \dd z^{1/2}.
	\end{equation*}
	
	\subsubsection{The Cauchy operator}
	Let $\alpha = h\, \dd z^{1/2} + \overline{H}\, \dd \bar{z}^{1/2}\in \Ar_h^{1/2}(\Omega)$. By Sokhotski--Plemelj theorem (see Remark~\ref{rem:SP} for its validity), the operator $\Cr_{\Gamma}$ defined in \eqref{eq:def-Cr} satisfies
	\begin{equation}\label{eq:jump}
		\Cr_{\Gamma}(\bb_{\Omega}\alpha) = \bb_{\Omega} \Jj^{1/2} \alpha
	\end{equation}
	Therefore, we have 
	\begin{equation*}
		\Jj^{1/2} = \bb_{\Omega}^{-1}\Cr_{\Gamma}\bb_{\Omega}\colon \Ar_h^{1/2}(\Omega)\to \Ar_h^{1/2}(\Omega).
	\end{equation*}
	Since $\bb_{\Omega}$ is an isometry, we have $\bb_{\Omega}^{-1} = \bb_{\Omega}^*$ and 
	\begin{equation}\label{eq:jump2}
		\Jj^{1/2} - \left(\Jj^{1/2}\right)^* = \bb_{\Omega}^{-1}\left(\Cc_{\Gamma} - \Cc_{\Gamma}^*\right) \bb_{\Omega}.
	\end{equation}
	Hence, we have the following equality of Hilbert--Schmidt norms
	\begin{equation}\label{eq:norm-Cauchy-J}
		\|\Cc_{\Gamma}-\Cc_{\Gamma}^*\|_{\rm HS\left(L^2(\Gamma)\right)} = \left\| \Jj^{1/2} - \left(\Jj^{1/2}\right)^* \right\|_{\rm HS\left(\Ar_h^{1/2}(\Omega)\right)}.
	\end{equation}
	By Cauchy's integral formula (see also Equation (33) and Proposition 3.3 in \cite{KSS2025}), the map $\Jj^{1/2}$ restricted to $ \Ar^{1/2}(\Omega)$ is the identity but we also have that $\Jj^{1/2}\colon \overline{\Ar^{1/2}(\Omega)}\to \Ar^{1/2}(\Omega)$. Hence, we can decompose $\Jj^{1/2}$ by blocks as
	\begin{equation}\label{eq:decompo-J-blocks}
		\Jj^{1/2} = \begin{pmatrix}
			I & \mathrm{j}^{1/2} \\
			0 & 0 
		\end{pmatrix}, \qquad \text{ with }\quad \mathrm{j}^{1/2}\coloneqq \left( \Jj^{1/2}\colon \overline{\Ar^{1/2}(\Omega)}\to \Ar^{1/2}(\Omega) \right).
	\end{equation}
	Hence, we have by blocks
	\begin{equation*}
		\Jj^{1/2} - \left(\Jj^{1/2}\right)^* = \begin{pmatrix}
			0 & \mathrm{j}^{1/2} \\
			- \left( \mathrm{j}^{1/2} \right)^* & 0 
		\end{pmatrix}.
	\end{equation*}
\begin{defi}[Proposition 4.5 of \cite{KSS2025}]
   The \emph{fermionic Grunsky operator} $\mathrm{Gr}^{1/2}_f\colon  \overline{\Ar^{1/2}(\D) } \to \Ar^{1/2}(\D)$ is given by
	\begin{equation}\label{eq:def_Grunsky}
		\begin{aligned} 
			\mathrm{Gr}^{1/2}_f = -\left. \widehat{f}^*\, \Jj^{1/2}\, (\widehat{f}^{-1})^*\right|_{\overline{\Ar^{1/2}(\D) }}  = -\left. \widehat{f}^*\, \mathrm{j}^{1/2}\, (\widehat{f}^{-1})^*\right|_{\overline{\Ar^{1/2}(\D) }} 
		\end{aligned} 
	\end{equation}
    where we have chosen a branch of $\sqrt{f'}$ to define $\widehat f \in \Confhat(\D,\Omega)$. The definition of $\mathrm{Gr}^{1/2}_f$ does not  depend on this choice though.  
\end{defi}
   
	We obtain 
	\begin{equation}\label{eq:J-G}
		\Jj^{1/2} - \left(\Jj^{1/2}\right)^* = - (\widehat{f}^*)^{-1}\,\begin{pmatrix}
			0 & \mathrm{Gr}^{1/2}_f \\
			- \left( \mathrm{Gr}^{1/2}_f \right)^* & 0 
		\end{pmatrix}\, \widehat{f}^*.
	\end{equation}
	By \eqref{eq:norm-Cauchy-J}, we deduce that 
	\begin{equation*}
		\|\Cc_{\Gamma}-\Cc_{\Gamma}^*\|_{\rm HS\left(L^2(\Gamma)\right)}^2  =2\, \left\| \mathrm{Gr}^{1/2}_f \right\|^2_{\rm HS\left(\overline{\Ar^{1/2}(\D) } \right) }.
	\end{equation*}
	Thanks to Proposition \ref{pr:Mob-Cauchy}, we have proved the following result.
	\begin{proposition}\label{pr:Mob-Grunsky}
		Let $\Gamma\subset \C$ be a Weil--Petersson curve bounding an open set $\Omega$. Let $f\colon\D\to \Omega$ be a uniformizing map. Then it holds
		\begin{equation*}
			\Mob(\Gamma)-4 = 4\pi^2\, \left\| \mathrm{Gr}^{1/2}_f \right\|^2_{\rm HS\left(\overline{\Ar^{1/2}(\D) } \right) }.
		\end{equation*}
	\end{proposition}
	
    \begin{remark}\label{rem:SP}
		This result answers Problem 5.6 in \cite{KSS2025}. If $\Gamma$ is a Weil--Petersson curve, then Proposition \ref{pr:Mob-Grunsky} shows that the operator $\mathrm{Gr}^{1/2}_f$ is Hilbert--Schmidt, which was already proved by \cite{KSS2025}. 
        Conversely, let $\Gamma$ be a rectifiable Jordan curve that satisfies the following Carleson condition 
		\begin{equation}\label{eq:Carleson}
			\sup_{p\in\Gamma,\, r>0} \frac{ |\Gamma\cap \D(p,r)| }{ r} <\infty.
		\end{equation}
    By Theorem~2.3 in \cite{L2018}, the formula \eqref{eq:jump} holds. Since $\Gamma$ is rectifiable, the decompositions \eqref{eq:decompo_L2} and \eqref{eq:decompo-J-blocks} also hold, see Theorem~2.16 and Proposition~3.3 in \cite{KSS2025}. Therefore, we can define the operator $\mathrm{Gr}^{1/2}_f$ by the formula \eqref{eq:J-G}.

		In this case, if $\mathrm{Gr}^{1/2}_f$ is Hilbert--Schmidt, 
        then 
        \eqref{eq:norm-Cauchy-J}  implies that the operator $\Cc_{\Gamma}-\Cc_{\Gamma}^*$ is Hilbert--Schmidt.  
        Following the proof of Proposition \ref{pr:Mob-Cauchy}, we obtain that the quantity $(1-\cos\theta_{\Gamma}(x,y))/|x-y|^2$ is integrable on $\Gamma\times \Gamma$. Hence, we obtain  $\Mob(\Gamma)<\infty$. Therefore $\Gamma$ is a Weil--Petersson curve by Lemmas \ref{lem:distortion} and \ref{lem:tangent-control}.
	\end{remark}
	
	\subsubsection{Conformal welding and Grunsky operator} 
	In this section, we define an action of the conformal welding on $\Ll^2(\s^1)$ and relate it to the fermionic Grunsky operator.
	
	We denote $\vp\coloneqq g^{-1}\circ f$, where $f$ and $g$ are defined at the beginning of Section \ref{sec:Grunksy}. We consider a lifting $\vp(\ee^{\ii \theta}) = \ee^{\ii h(\theta)}$ with 
	\begin{equation*}
		\begin{cases} 
			h'>0,\\
			h(\theta +2\pi) = h(\theta)+2\pi.
		\end{cases}
	\end{equation*}
	Then it holds
	\begin{equation*}
		\vp'(\ee^{\ii \theta}) = -\ii \ee^{-\ii \theta}\, \frac{\dd \vp(\ee^{\ii \theta})}{\dd \theta} = h'(\theta)\, \ee^{\ii (h(\theta)-\theta)}.
	\end{equation*}
	Since $h(\theta+2\pi)-(\theta+2\pi)  =h(\theta)-\theta$, we can define
	\begin{equation*}
		\sqrt{\vp'(\ee^{\ii \theta})} \coloneqq \sqrt{h'(\theta)}\, \ee^{\ii \frac{h(\theta)-\theta}{2}}.
	\end{equation*}
	We define the operator $W_{\vp}\colon \Ll^2(\s^1)\to \Ll^2(\s^1)$ by the formula
	\begin{equation*}
		\forall u\in \Ll^2(\s^1),\qquad W_{\vp} u \coloneqq \sqrt{\vp'}\, u\circ\vp.
	\end{equation*}
	\begin{claim}
		The operator $W_{\vp}$ is an isometry.
	\end{claim}
	\begin{proof}
		By change of variable $w=\vp(z)$, we have 
		\begin{equation*}
			\int_{\s^1} \left| \sqrt{\vp'(z)}\, u\circ\vp(z)\right|^2\, |\dd z| = \int_{\s^1} \left|  u\circ\vp(z)\right|^2\, |\vp'(z)|\, |\dd z| = \int_{\s^1} |u(w)|^2\, |\dd w|. \qedhere
		\end{equation*}
	\end{proof}
	The relationship between the operator $W_{\vp}$ and the fermionic Grunsky operator is obtained through the Faber operator defined as follows.
	
	We decompose $\Ll^2(\s^1) = H_+\operp H_-$ as in \eqref{eq:decompo_L2}, so that we have 
	\begin{equation*}
		\begin{cases} 
			H_+ \coloneqq \Span_{\,\Ll^2(\s^1)}\left\{z\mapsto z^n :n\geqslant 0 \right\}, \\[2mm]
			H_- \coloneqq \Span_{\,\Ll^2(\s^1)}\left\{z\mapsto z^n : n\leqslant -1 \right\}
		\end{cases} 
	\end{equation*}
	The \textit{Faber operator} $\mathrm{I}^{1/2}_{\widehat{f}}\colon \overline{\Ar^{1/2}(\D)} \to \Ar^{1/2}(\Omega)$ is defined, using Definition 3.6, Equation (33) and the text above Remark 3.9 in \cite{KSS2025}, by the following formula for any $w\in\Omega^*$ and $\alpha = \overline{H}\, \dd \bar{z}^{1/2}\in \overline{\Ar^{1/2}(\D)}$
	\begin{equation*}
		\mathrm{I}^{1/2}_{\widehat{f}}\alpha(w) \coloneqq -\Jj^{1/2}_*(\widehat{f}^{-1})^*\alpha(w) =  \left(-\frac{1}{2\ii \pi} \int_{\s^1} \overline{H(\eta)}\, \frac{\sqrt{f'(\eta)}}{f(\eta)-w}\, |\dd \eta|\right) \dd w^{1/2}.
	\end{equation*}
	Pulling back by $\widehat{g}$, we obtain the following operator $U$ and $U(\alpha) \in \Ar^{1/2}(\D^*)$ is given by 
	\begin{equation}\label{eq:defU}
		U(\alpha)(\zeta)\coloneqq \widehat{g}^*\mathrm{I}^{1/2}_{\widehat{f}}\alpha(\zeta) = \left( -\frac{1}{2\ii \pi} \int_{\s^1} \overline{H(\eta)}\, \frac{\sqrt{f'(\eta)}\, \sqrt{g'(\zeta)}}{f(\eta)-g(\zeta)}\, |\dd \eta|\right) \dd \zeta^{1/2}.
	\end{equation}
	Applying the operator $W_{\vp}\circ\bb_{\D^*}$, we have\footnote{Given $\beta \in \Ar^{1/2}(\Omega^*)$, there holds $\bb_{\D^*}(\widehat{g}^*\beta) = \sqrt{g'}\, (\bb_{\Omega^*}\beta)\circ g$. Similarly on $\D$ with $g$ replaced by $f$.} 
    \begin{equation*}
        \begin{aligned} 
            W_{\vp}(\bb_{\D^*}U(\alpha)) & =  W_{\vp} \left( \sqrt{g'}\, \left(\bb_{\Omega^*} {\rm I}_{\widehat{f}}^{1/2}\alpha\right)\circ g \right) \\[2mm]
            & =   \sqrt{f'}\, \left(\bb_{\Omega^*} {\rm I}_{\widehat{f}}^{1/2}\alpha\right)\circ f\\[2mm]
            & =\bb_{\D}\widehat{f}^*\left( \bb_{\Omega}^{-1} \bb_{\Omega^*} {\rm I}_{\widehat{f}}^{1/2} \alpha \right) \\[2mm]
            & = - \bb_{\D}\widehat{f}^*\left( \bb_{\Omega}^{-1} \bb_{\Omega^*} \Jj_{*}^{1/2} (\widehat{f}^{-1})^*\alpha \right).
        \end{aligned} 
    \end{equation*}
    Using Equation (35) p.\ 5338 in \cite{KSS2025}, we have 
    \begin{equation*}
        W_{\vp}(\bb_{\D^*}U(\alpha)) =\bb_{\D} \widehat{f}^* \left( (\widehat{f}^{-1})^*\alpha-\Jj^{1/2} (\widehat{f}^{-1})^*\alpha  \right).
    \end{equation*}
    By \eqref{eq:def_Grunsky}, we obtain
	\begin{equation}\label{eq:W-graphG}
			W_{\vp}(\bb_{\D^*}U(\alpha))  = \bb_{\D}\left( \mathrm{Gr}_{f}^{1/2}\alpha + \alpha \right).
	\end{equation}
	\begin{claim}
		Given $\alpha\in \overline{\Ar^{1/2}(\D)}$, we have $\bb_{\D^*}U(\alpha) \in H_-$.
	\end{claim}
	\begin{proof}
        Recall that $U(\alpha)$ is given by the image 
		\begin{equation*}
			\alpha\in \overline{\Ar^{1/2}(\D)} \xrightarrow[]{(\widehat{f}^{-1})^*}\overline{\Ar^{1/2}(\Omega)} \xrightarrow[]{-\Jj^{1/2}_*} \Ar^{1/2}(\Omega^*) \xrightarrow[]{\widehat{g}^*} \Ar^{1/2}(\D^*).
		\end{equation*}
		Since half-order differentials on $\D^*$ are extendible at $\infty$, the map $t^{-1}U(\alpha)(\frac{1}{t})$ is holomorphic at $0$ and we obtain
		\begin{equation*}
			U(\alpha)(\zeta) = \sum_{n\in\N} u_n\, \zeta^{-n-1}\, \dd \zeta^{1/2}. 
		\end{equation*}
		Therefore, $\bb_{\D^*}U(\alpha)\in H_-$.
	\end{proof}
	We now express the fermionic Grunsky operator in terms of the operator $W_{\vp}$. To this end, we decompose $W_{\vp}$ by blocks corresponding to the decomposition $\Ll^2(\s^1) = H_+\operp H_-$ as follows
	\begin{equation*}
		W_{\vp} = \begin{pmatrix}
			A & B \\
			C & D 
		\end{pmatrix}.
	\end{equation*}
	Using \eqref{eq:W-graphG} and $\mathrm{Gr}^{1/2}_f\colon  \overline{\Ar^{1/2}(\D) } \to \Ar^{1/2}(\D)$, we obtain 
	\begin{equation}\label{eq:formulaD}
		D\circ \bb_{\D^*}U = \bb_{\D}.
	\end{equation}
	By \cite[Corollary 3.7]{KSS2025} the Faber operator $\mathrm{I}^{1/2}_{\widehat{f}}\colon \overline{\Ar^{1/2}(\D)} \to \Ar^{1/2}(\Omega^*)$ is an isomorphism and thus by \eqref{eq:defU}, the operator $U\colon \overline{\Ar^{1/2}(\D)}\to \Ar^{1/2}(\D^*)$ is also an isomorphism. Hence, we obtain that $D = \left.\bb_{\D}\circ U^{-1}\circ \bb_{\D^*}^{-1}\right|_{H_-}$ is invertible. Restricting \eqref{eq:W-graphG} to $\bb_{\D}^{-1}(H_-)$, we deduce that the image of $W_{\vp}|_{H_-}$ coincides with the graph of $\bb_{\D}\mathrm{Gr}_{f}^{1/2}$:
	\begin{equation*}
		W_{\vp}(H_-) = \mathrm{Graph}\left(\bb_{\D}\mathrm{Gr}_{f}^{1/2}\right).
	\end{equation*}
	We also get from \eqref{eq:W-graphG} and \eqref{eq:formulaD}
	\begin{equation*}
		\begin{aligned} 
			W_{\vp}\circ D^{-1} & = W_{\vp} \circ \bb_{\D^*}\circ U\circ \left( \bb_{\D}|_{H_-}\right)^{-1}\\[2mm]
			& = \bb_{\D}\left( \mathrm{Gr}_{f}^{1/2} + \id \right)\circ \left( \bb_{\D}|_{H_-}\right)^{-1}.
		\end{aligned} 
	\end{equation*}
	Restricting to the block sending $H_-\to H_+$, we end up with the following operator 
	\begin{equation}\label{def:G}
		G \coloneqq \bb_{\D} \mathrm{Gr}_{f}^{1/2}  \left( \bb_{\D}|_{H_-}\right)^{-1} = B\, D^{-1} \colon H_-\to H_+.
	\end{equation}
	On the other hand, since $W_{\vp}$ is an isometry, we have $W_{\vp}^*W_{\vp}=  \id_{\Ll^2(\s^1)}$ and thus,
	\begin{equation*}
		\id = \begin{pmatrix}
			A^* & C^*\\
			B^* & D^*
		\end{pmatrix} \begin{pmatrix}
			A & B\\
			C & D
		\end{pmatrix} = \begin{pmatrix}
			A^*A + C^*C & A^*B + C^*D \\
			B^*A + D^*C & B^*B + D^*D
		\end{pmatrix}.
	\end{equation*}
	We have $B^*A+D^*C=0$, which implies that $(D^*)^{-1}B^* = -CA^{-1}$. Hence, it holds
	\begin{equation}\label{eq:CAG}
		CA^{-1} = -G^*.
	\end{equation}
	Moreover, the equation $B^*B +D^*D=\id_{H_-}$ implies that 
	\begin{equation*}
		\id_{H_-} = D^*G^*G D + D^*D = D^*(G^*G + \id_{H_-})D.
	\end{equation*}
	We obtain 
	\begin{equation*}
		G^*G + \id_{H_-} = (DD^*)^{-1}.
	\end{equation*}
	Thus, we have 
	\begin{equation}\label{eq:logdetDD}
		\log\det( DD^*) = \log\det (\id_{H_-}+G^*G)^{-1}.
	\end{equation}
	We also deduce that $AA^*+BB^*=\id_{H_+}$ from the decomposition by blocks of the identity $W_{\vp} W_{\vp}^* = \id_{\Ll^2(\s^1)}$. Hence, we obtain 
	\begin{equation*}
			\log\det(AA^*) = \log\det(\id_{H_+}-BB^*).
	\end{equation*}
	Since $\id_{H_+}-BB^*$ and $\id_{H_-}-B^*B$ have same eigenvalues\footnote{If $BB^*x=\lambda x$ and $\lambda x\neq 0$, then $B^*x\neq 0$. Applying $B^*$, we get that $\lambda$ is also an eigenvalue for $B^*B$. Hence, $\id_{H_+}-BB^*$ and $\id_{H_-}-B^*B$ have same non-zero eigenvalues.}, we obtain 
	\begin{equation*}
		\log\det(AA^*) = \log\det(\id_{H_-}-B^*B).
	\end{equation*}
	Using the relation $D^*D=\id_{H_-} - B^*B $ together with \eqref{eq:logdetDD}, we end up with the following identities. 
    \begin{lemma}\label{lem:logdetA} 
    We have
    $$\log\det(AA^*) = \log\det(D^*D) = \log\det (\id_{H_-}+G^*G)^{-1} = - \log\det (\id_{H_-}+G^*G).
	$$    		
    \end{lemma}

	\subsubsection{Variations of the welding operator.}\label{sec:differentiate-welding}
	In this section, we perform preliminary computations to express the Loewner energy in terms of the fermionic Grunsky operator. 
	
	Let $(\Gamma_t)_{t\in(-1,1)}$ be a smooth family of Weil--Petersson curves with $\Gamma_0\coloneqq \Gamma$. For each parameter $-1<t<1$, we denote $\vp_t$ the conformal welding of $\Gamma_t$ with a normalization so that $\vp_0=\vp$ is the conformal welding of $\Gamma$. We define $u_t\coloneqq \vp^{-1}\circ\vp_t$ and we denote $\IL(\vp_t)\coloneqq \IL(\Gamma_t)$.
	
	Then we have $\vp_t=\vp \circ u_t$ and $W_{\vp_t} = W_{u_t}\circ W_{\vp}$.
	Denoting $\dot{u} \coloneqq \dr_t u_t|_{t=0}$, we have for any function $F\in\Cd^{\infty}(\s^1)$,
	\begin{equation}\label{eq:der_Wu}
		\left. \frac{\dd}{\dd t}\right|_{t=0} W_{u_t}F = F'\, \dot{u} + \frac{\dot{u}'}{2}\, F.
	\end{equation}
	Moreover, we have $|u_t|^2=1$ for any $-1<t<1$. Differentiating at $t=0$ where $u_0(z)=z$, we obtain $\Re\left( \bar{z}\,\dot{u}\right)=0$. Thus, we have $\bar{z}\, \dot{u} = -z\, \overline{\dot{u}}$, and hence $\dot{u} = -z^2\, \overline{\dot{u}}$. Considering the Fourier coefficients $(c_n)_{n\in\Z}$ of $\dot{u}$, we obtain the following relation
	\begin{equation*}
			\dot{u}(z) = \sum_{n\in\Z} c_n z^n  = -\sum_{n\in\Z} \overline{c_n}\, z^{2-n}  = -\sum_{m\in\Z} \overline{c_{2-m}}\, z^m .
	\end{equation*}
	We deduce that for all $n\in\N$, there holds 
	\begin{equation}\label{eq:relation_c}
		c_n= -\overline{c_{2-n}}.
	\end{equation}  
	\begin{claim}\label{cl:Mobius} 
	Up to Möbius transformation, we can assume $c_0=0$ and $c_1=0$. 
	\end{claim} 
	\begin{proof}
		We have $c_1=-\overline{c_1}$. Hence, there exists $\beta\in\R$ such that $c_1=\ii \beta$, and the first three modes of $\dot{u}$ are given by $c_0+c_1 z +c_2 z^2 = c_2 +\ii \beta z - \overline{c_0}\, z^2$. Such maps coincides with the derivatives of Möbius transformations. Indeed, for $-1<t<1$, we consider the following family of Möbius transformations with $m_0=\id$, $a_t\in\C$ and $\theta_t\in\R$
		\begin{equation*}
			 m_t(z) \coloneqq \ee^{\ii \theta_t}\frac{z- a_t}{1-\overline{a_t}\, z}.
		\end{equation*}
		Then it holds
		\begin{equation*}
			\dot{m}(z) =  - \dot{a} + \ii \dot{\theta}\, z+ \overline{\dot{a}}\, z^2.
		\end{equation*}
		Therefore, there exists a family of Möbius transformation with 
		\begin{equation*}
			\dot{m}(z) =  -\left( c_0+c_1 z +c_2 z^2\right).
		\end{equation*}
		We then have 
		\begin{equation*}
			\left.\frac{\dd}{\dd t}\right|_{t=0} u_t\circ m_t(z) = \dot{u}(z) + \dot{m}(z),
		\end{equation*}
		for which the components in the directions $1$, $z$ and $z^2$ are $0$.
	\end{proof}
	
	We consider $v(z) \coloneqq \sum_{n\geqslant 2} c_n z^n$, which is the part of $u$ that extends holomorphically to $\D$. We obtain 
	\begin{equation}\label{eq:decompo-u}
		\begin{aligned} 
			v(z)-z^2\overline{v(z)} & = \sum_{n\geqslant 2} \left( c_n\, z^n - \overline{c_n}\, z^{2-n} \right) = \sum_{n\geqslant 2} \left( c_n\, z^n + c_{2-n}\, z^{2-n} \right) = \dot{u}(z).
		\end{aligned} 
	\end{equation}
	We differentiate in $z$ 
	\begin{equation*}
		\begin{aligned}
			\dot{u}'(z)= \sum_{n\geqslant 2} \left( n\, c_n\, z^{n-1}-(2-n)\, \overline{c_n}\, z^{1-n} \right).
		\end{aligned}
	\end{equation*}
	We now come back to \eqref{eq:der_Wu}. 
	We consider the direction of $F(z)=z^{-k}$ for $k\geqslant 0$ in \eqref{eq:der_Wu}, we obtain 
	\begin{equation*}
		\begin{aligned} 
			& (z^{-k})'\, \dot{u}(z) + \frac{\dot{u}'(z)}{2}\, z^{-k} \\[2mm]
			& = -\sum_{n\geqslant 2} k\left( c_n\, z^{n-k-1}-\overline{c_n}\, z^{-k-n+1} \right) + \frac{1}{2}\sum_{n\geqslant 2} \left( n\, c_n\, z^{n-k-1}-(2-n)\overline{c_n}\, z^{1-n-k} \right) \\[2mm]
			& =\sum_{n\geqslant 2} \left[ \left(-k+\frac{n}{2}\right)\, c_n\, z^{n-k-1}  - \left(-k+ \frac{2-n}{2}\right)\, \overline{c_n}\, z^{1-n-k} \right]  \\[2mm]
			& = \sum_{m\geqslant -k+1} \frac{m-k+1}{2}\, c_{m+k+1}\, z^m - \sum_{p\leqslant -k-1} \frac{1+p-k}{2}\, \overline{c_{1-k-p}}\, z^p.
		\end{aligned} 
	\end{equation*} 
	Since $\overline{c_{1-k-p}} = -c_{p+k+1}$ by \eqref{eq:relation_c} and $c_1=0$ by Claim \ref{cl:Mobius}, we obtain 
	\begin{equation}\label{eq:derB1} 
			(z^{-k})'\, \dot{u}(z) + \frac{\dot{u}'(z)}{2}\, z^{-k} = \sum_{m\in\Z} \frac{m-k+1}{2}\, c_{m+k+1}\, z^m .
	\end{equation}
	We rewrite the relation $W_{\vp_t} = W_{u_t}\circ W_{\vp}$ by blocks as
	\begin{equation*}
		\begin{pmatrix}
			A_{\vp_t} & B_{\vp_t}\\
			C_{\vp_t} & D_{\vp_t}
		\end{pmatrix} = \begin{pmatrix}
			A_{u_t} & B_{u_t}\\
			C_{u_t} & D_{u_t}
		\end{pmatrix}\, \begin{pmatrix}
			A_{\vp} & B_{\vp}\\
			C_{\vp} & D_{\vp}
		\end{pmatrix}.
	\end{equation*}
	Therefore, we have 
	\begin{equation}\label{eq:derA}
		\left. \frac{\dd}{\dd t}\right|_{t=0} A_{\vp_t} = \left(\left. \frac{\dd}{\dd t}\right|_{t=0} A_{u_t} \right)\, A_{\vp} + \left( \left. \frac{\dd}{\dd t}\right|_{t=0} B_{u_t}\right)\, C_{\vp}.
	\end{equation}
	The second derivative of the right-hand side is given by \eqref{eq:derB1}
	\begin{equation}\label{eq:der-B}
			\left. \frac{\dd}{\dd t}\right|_{t=0} B_{u_t} = \left(\frac{1+p-k}{2}\, c_{p+k+1} \right)_{p,k\geqslant 1}.
	\end{equation}

	\subsection{Log-determinant and Loewner energy.}\label{sec:log_det_I}
	Throughout this subsection, we use the notations of Section \ref{sec:Grunksy}.
	The goal of this subsection is to prove the following result.
	\begin{theorem}\label{th:Loewner-Grunsky}
		For any Weil--Petersson curve $\Gamma\subset \C$ with conformal welding $\vp$ and conformal map $f\colon \D\to \Omega$, it holds
		\begin{equation*}
			\IL(\Gamma) = -24\log\det(A_{\vp}A_{\vp}^*) = 24\log\det\left(I+\left[\mathrm{Gr}^{1/2}_{f}\right]^*\, \mathrm{Gr}^{1/2}_{f}\right).
		\end{equation*}
	\end{theorem}
    The second identity follows from Lemma~\ref{lem:logdetA} and \eqref{def:G}. For the first identity, 
	we differentiate along a path $(\vp_t)_{t\in(-1,1)}$ as in Section \ref{sec:differentiate-welding} and obtain 
	\begin{equation}\label{eq:derivative-logdet}
		\begin{aligned} 
			\left. \frac{\dd }{\dd t}\right|_{t=0} \log \det(A_{\vp_t} A^*_{\vp_t}) & = \tr\left[ (A_{\vp} A^*_{\vp})^{-1} \left. \frac{\dd }{\dd t}\right|_{t=0}(A_{\vp_t} A^*_{\vp_t})  \right] \\[2mm]
			& = \tr\left[ (A_{\vp} A^*_{\vp})^{-1} \left( \left( \left. \frac{\dd }{\dd t}\right|_{t=0}A_{\vp_t}\right)\ A^*_{\vp} + A_{\vp}\ \left(\left. \frac{\dd }{\dd t}\right|_{t=0} A^*_{\vp_t}\right)  \right)  \right] \\[2mm]
			& = \tr\left[ A_{\vp}^{-1}  \left( \left. \frac{\dd }{\dd t}\right|_{t=0}A_{\vp_t} \right) \right] + \tr\left[ (A_{\vp}^*)^{-1}  \left( \left. \frac{\dd }{\dd t}\right|_{t=0}A_{\vp_t}^* \right) \right] \\[2mm]
			& = 2\, \Re\, \tr\left[ A_{\vp}^{-1}  \left( \left. \frac{\dd }{\dd t}\right|_{t=0}A_{\vp_t} \right) \right].
		\end{aligned} 
	\end{equation} 
	Using \eqref{eq:derA}, we obtain
	\begin{equation*} 
			\left. \frac{\dd }{\dd t}\right|_{t=0} \log \det(A_{\vp_t} A^*_{\vp_t}) = 2\, \Re\, \tr\left[ \left. \frac{\dd }{\dd t}\right|_{t=0}A_{u_t}  \right] + 2\, \Re\, \tr\left[ C_{\vp}\, A_{\vp}^{-1} \left( \left. \frac{\dd}{\dd t}\right|_{t=0} B_{u_t}\right) \right].
	\end{equation*}
	Since $W_{u_t}$ is unitary with $W_{u_0}=\id_{\Ll^2(\s^1)}$, the first term vanishes. We obtain, using \eqref{eq:CAG} for the second term on the right-hand side 
	\begin{equation}\label{eq:der1-logdet}
		\left. \frac{\dd }{\dd t}\right|_{t=0} \log \det(A_{\vp_t} A^*_{\vp_t}) = -2\, \Re\, \tr\left[ G^* \left( \left. \frac{\dd}{\dd t}\right|_{t=0} B_{u_t}\right) \right].
	\end{equation}
	We decompose $G^*$ in the basis $(z^m)_{m\geqslant 0}$ as follows
	\begin{equation}\label{eq:decompo-Gstar}
		G^*(z^m) = \sum_{k\geqslant 1} r_{k,m}\, z^{-k}.
	\end{equation}
	Therefore, we have by \eqref{eq:der-B}
	\begin{equation}\label{eq:der2-logdet}
		-2\, \tr\left[ G^* \left( \left. \frac{\dd}{\dd t}\right|_{t=0} B_{u_t}\right) \right] = -2\sum_{k\geqslant 1, m\geqslant 0} \frac{1+m-k}{2}\, r_{k,m}\, c_{m+k+1}.
	\end{equation}
	The kernel of $G$ is given by \cite[Corollary 4.6]{KSS2025} as
	\begin{equation}\label{eq:def-Grunsky}
		\mathrm{Gr}^{1/2}_{f}\left(\overline{H}\, \dd \bar{z}^{1/2}\right)(z) = \frac{1}{2\ii \pi}\int_{\zeta\in\s^1} K(z,\zeta)\, \overline{H}(\zeta)\, |\dd \zeta|\, \dd z^{1/2},
	\end{equation}
	with 
	\begin{equation}\label{eq:kernelK}
		K(z,\zeta)\coloneqq \frac{1}{\zeta-z} - \frac{\sqrt{f'(\zeta)}\, \sqrt{f'(z)}}{f(\zeta)-f(z)} = \sum_{m,n\geqslant 0} a_{m,n}\, z^n\, \zeta^m.
	\end{equation}
	If $u(\zeta) = \bb_{\D}\alpha(\zeta) \in H_-$ with $\alpha = \overline{H}(z)\, \dd \bar{z}^{1/2}\in \overline{\Ar^{1/2}(\D)}$, then $u(\zeta) = \overline{H(\zeta)}\, \overline{\tau_{\s^1}(\zeta)}$, with $\tau_{\s^1}(\zeta) = \ii\zeta$, by \eqref{eq:def-b}. 
	We obtain the following expression of $G$ for $u\in H_-$,
	\begin{equation*}
		\begin{aligned} 
			Gu(z) & = \bb_{\D}\, \mathrm{Gr}_{f}^{1/2}\left(\overline{H}(z)\, \dd \bar{z}^{1/2}\right)\\[2mm]
			& = \bb_{\D} \left( \frac{1}{2\ii \pi}\int_{\zeta\in\s^1} K(z,\zeta)\, \overline{H}(\zeta)\, |\dd \zeta|\, \dd z^{1/2}\right) \\[2mm]
			& = \frac{1}{2\ii \pi}\int_{\s^1} K(z,\zeta)\, \overline{H}(\zeta)\, |\dd \zeta| \\[2mm]
			& = \frac{1}{2\ii \pi}\int_{\s^1} K(z,\zeta)\, \frac{u(\zeta)}{-\ii\bar{\zeta}}\, |\dd \zeta| \\[2mm]
			& = \frac{1}{2\ii \pi}\int_{\s^1} K(z,\zeta)\, u(\zeta)\, \ii \zeta\, |\dd \zeta| \\[2mm]
			& = \frac{1}{2\ii \pi}\int_{\s^1} K(z,\zeta)\, u(\zeta)\, \dd \zeta.
		\end{aligned} 
	\end{equation*}
	\begin{claim}
		It holds $r_{k,m} = \overline{a_{k-1,m}}$.
	\end{claim}
	\begin{proof}
		In the basis $(z^{-k})_{k\geqslant 1}$, the operator $G\colon H_-\to H_+$ is given by 
		\begin{equation*}
			\begin{aligned} 
			G(z^{-k}) & = \frac{1}{2\ii\pi} \sum_{m,n\geqslant 0} a_{m,n}\, z^n \int_{\s^1} \zeta^{m-k}\, \dd \zeta \\[2mm]
			& = \frac{1}{2\pi} \sum_{m,n\geqslant 0} a_{m,n}\, z^n \int_{\s^1} \zeta^{m-k+1}\, |\dd \zeta| \\[2mm]
			& = \sum_{n\geqslant 0} a_{k-1,n}\, z^n.
			\end{aligned} 
		\end{equation*}
		Hence, we obtain from \eqref{eq:decompo-Gstar} that $r_{k,m} = \overline{a_{k-1,m}}$.
	\end{proof}
	Coming back to \eqref{eq:der1-logdet} and \eqref{eq:der2-logdet}, we obtain 
	\begin{equation*}
		\begin{aligned} 
		\left. \frac{\dd }{\dd t}\right|_{t=0} \log \det(A_{\vp_t} A^*_{\vp_t}) & = -2\Re\left[ \sum_{k\geqslant 1, m\geqslant 0} \frac{1+m-k}{2}\, c_{m+k+1}\, \overline{a_{k-1,m}} \right] \\[4mm]
		& = \Re \Bigg[ \sum_{p=2}^{\infty} c_p\Bigg( \sum_{\substack{k\geqslant 1, m\geqslant 0\\ m+k+1=p}} (k-m-1)\, \overline{a_{k-1,m}}\Bigg) \Bigg] \\[4mm]
		&  = \Re \Bigg[ \sum_{p=2}^{\infty} c_p\Bigg( \sum_{\substack{\ell\geqslant 0, m\geqslant 0\\ m+\ell+2=p}} (\ell-m)\, \overline{a_{\ell,m}}\Bigg) \Bigg] .
		\end{aligned} 
	\end{equation*}
    Since $K(z,\zeta)=-K(\zeta,z)$, we obtain $a_{\ell,m}=-a_{m,\ell}$ and thus, it holds
    \begin{equation}\label{eq:der-AA}
        \left. \frac{\dd }{\dd t}\right|_{t=0} \log \det(A_{\vp_t} A^*_{\vp_t}) = \Re \Bigg[ \sum_{p=2}^{\infty} c_p\Bigg( \sum_{\substack{\ell + m\geqslant 1\\ m+\ell+2=p}} (m-\ell)\, \overline{a_{m,\ell}}\Bigg) \Bigg] .
    \end{equation}
	In the above sum, the indices are nonnegative. On the other hand, it holds
	\begin{equation*}
		(\dr_z-\dr_{\zeta})K(z,z) = \sum_{n+m\geqslant 1} a_{m,n}\, (n-m)\, z^{n+m-1} = \sum_{p=0}^{\infty}\Bigg( \sum_{\substack{n+m\geqslant 1\\ n+m-1=p}} (n-m)\, a_{m,n} \Bigg)\, z^p .
	\end{equation*}
	In the above sum, the indices are nonnegative. Therefore, we get
	\begin{equation}\label{eq:int_derK}
		\begin{aligned} 
			\int_{\s^1} z^{-3}\, v(z)\, \overline{ (\dr_z-\dr_{\zeta})K(z,z)}\, |\dd z| & = \sum_{k,p\geqslant 0} c_k \Bigg( \sum_{\substack{n+m\geqslant 1\\ n+m-1=p}} (n-m)\, \overline{a_{m,n}} \Bigg) \int_{\s^1} z^{k-p-3}\, |\dd z|\\[2mm]
			& = 2\pi \sum_{k=0}^{\infty} c_k\, \Bigg( \sum_{\substack{n+m\geqslant 1\\ n+m-1=k-3}} (n-m)\, \overline{a_{m,n}} \Bigg) \\[2mm]
			& = 2\pi \sum_{k=2}^{\infty} c_k\, \Bigg( \sum_{\substack{n+m\geqslant 1\\ n+m+2=k}} (n-m)\, \overline{a_{m,n}} \Bigg).
		\end{aligned} 
	\end{equation}
	We use that $c_0=c_1=0$ by Claim \ref{cl:Mobius} in the last line.
	Combining \eqref{eq:der-AA} and \eqref{eq:int_derK} by identifying $\ell=n$, we obtain 
	\begin{equation}\label{eq:der-logA}
		\left. \frac{\dd }{\dd t}\right|_{t=0} \log \det(A_{\vp_t} A^*_{\vp_t}) = - \frac{1}{2 \pi }\, \Re\left[ \int_{\s^1} z^{-3}\, v(z)\, \overline{ (\dr_z-\dr_{\zeta})K(z,z)}\, |\dd z| \right].
	\end{equation}
	\begin{claim}\label{cl:der-K}
		There holds
		\begin{equation*}
			\left(\dr_z-\dr_{\zeta}\right) K(z,z) = \frac{\Sr f(z)}{6}.
		\end{equation*}
	\end{claim}
	\begin{proof} 
	Given $z\in \s^1$, we compute the asymptotic expansion as $\zeta\to z$ in \eqref{eq:kernelK} as follows
	\begin{equation*}
		\begin{aligned} 
			& K(z,\zeta)\\
			 & \ust{\zeta\to z}{=} \frac{1}{\zeta-z} - \frac{\sqrt{f'(z)}\, \sqrt{f'(z) + (\zeta-z)f''(z) + \frac{(\zeta-z)^2}{2} f'''(z) + O(|\zeta-z|^3)}}{(\zeta-z)f'(z) + \frac{(\zeta-z)^2}{2} f''(z) + \frac{(\zeta-z)^3}{6} f'''(z) + O(|\zeta-z|^4)} \\[3mm]
			& \ust{\zeta\to z}{=} \frac{1}{\zeta-z} - \frac{1}{\zeta-z}\, \frac{\sqrt{1 + (\zeta-z)\frac{ f''(z) }{f'(z)} + \frac{(\zeta-z)^2}{2} \frac{f'''(z)}{f'(z)} + O(|\zeta-z|^3)}}{1 + \frac{\zeta-z}{2} \frac{f''(z)}{f'(z)} + \frac{(\zeta-z)^2}{6} \frac{ f'''(z) }{f'(z)}+ O(|\zeta-z|^3)} \\[3mm]
			& \ust{\zeta\to z}{=} \frac{1}{\zeta-z} - \frac{1}{\zeta-z}\, \left( 1 + \frac{\zeta-z}{2}\frac{ f''(z) }{f'(z)} + \frac{(\zeta-z)^2}{4} \frac{f'''(z)}{f'(z)} - \frac{(\zeta-z)^2}{8} \left(\frac{f''(z)}{f'(z)}\right)^2 + O(|\zeta-z|^3) \right) \\[2mm]
			&\qquad \times \left( 1 - \frac{\zeta-z}{2} \frac{f''(z)}{f'(z)} - \frac{(\zeta-z)^2}{6} \frac{ f'''(z) }{f'(z)}+ \frac{(\zeta-z)^2}{4}\, \left(\frac{f''(z)}{f'(z)}\right)^2 +  O(|\zeta-z|^3) \right) \\[3mm]
			& \ust{\zeta\to z}{=} \frac{1}{\zeta-z} \\
            & - \frac{1}{\zeta-z}\, \left( 1 + \left(\frac{1}{4}-\frac{1}{6}\right) \, (\zeta-z)^2\, \frac{ f'''(z) }{f'(z)} +\left( \frac{1}{4}-\frac{1}{4}-\frac{1}{8} \right)\, (\zeta-z)^2\, \left(\frac{f''(z)}{f'(z)}\right)^2 + O(|\zeta-z|^3) \right) \\[3mm]
			& \ust{\zeta\to z}{=}  -\frac{1}{12} \, (\zeta-z)\, \frac{ f'''(z) }{f'(z)} +\frac{1}{8}\, (\zeta-z)\, \left(\frac{f''(z)}{f'(z)}\right)^2 + O(|\zeta-z|^2) \\[3mm]
			& \ust{\zeta\to z}{=}  -\frac{\zeta-z}{12} \, \Sr f(z) + O(|\zeta-z|^2).
		\end{aligned} 
	\end{equation*}
	Therefore, we obtain 
	\begin{equation*}
		\left(\dr_z-\dr_{\zeta}\right) K(z,z) = \frac{\Sr f(z)}{6}.
	\end{equation*}
	\end{proof} 
	Plugging this into \eqref{eq:der-logA}, we obtain 
	\begin{equation}\label{eq:der-logA2}
		\left. \frac{\dd }{\dd t}\right|_{t=0} \log \det(A_{\vp_t} A^*_{\vp_t}) = -\frac{1}{12 \pi }\, \Re\left[ \int_{\s^1} z^{-3}\, v(z)\, \overline{ \mathcal{S}f(z)}\, |\dd z| \right].
	\end{equation}
	On the other hand, the derivative of the Loewner energy is given by the following formula.
	\begin{claim}
		There holds
		\begin{equation*}
				\left. \frac{\dd}{\dd t} \right|_{t=0} \IL(\vp_t)  = \frac{2}{\pi}\, \Re\left(  \int_{\s^1}  z^3\, \overline{v(z)}\, \Sr(f)(z)\,  |\dd z| \right) .
		\end{equation*}
	\end{claim}
	\begin{proof} 
		Since each $u_t\colon\s^1\to\s^1$ is quasisymmetric (conformal weldings of Weil--Petersson curves are quasisymmetric), there exists a smooth family $U_t\colon\D\to\D$ of quasiconformal homeomorphism that extends $u_t$, with $U_0=\id_{\D}$. We denote by $\mu_{U_t^{-1}}$ the Beltrami coefficient of $U_t^{-1}\colon \D\to\D$. We then define the family of Beltrami coefficients
		\begin{equation*}
			\nu_t(w) \coloneqq \begin{cases}
				0 & \text{ if }w\in\Omega^*, \\[2mm]
				\mu_{U_t^{-1}}\big(f^{-1}(w)\big)\, \frac{f'(f^{-1}(w)) }{\overline{ f'(f^{-1}(w)) } } & \text{ if }w\in\Omega.
			\end{cases}
		\end{equation*}
		Since $\|\nu_t\|_{\Ll^{\infty}(\C)}<1$ and using \cite[Theorem 11.1.2]{IM2001}, there exists a map $\omega_t\colon\Chat\to\Chat$ that fixes $-1,-\ii,1$ and is a solution to 
		\begin{equation*}
			\dr_{\bar{z}} \omega_t = \nu_t\, \dr_z \omega_t.
		\end{equation*}
		Since $\nu_0=0$, we have $\omega_0=\id_{\Chat}$. Since $\nu_t|_{\Omega^*}=0$, the map $g_t\coloneqq \omega_t\circ g\colon \D^*\to\Chat$ is biholomorphic. Moreover, the map $f_t\coloneqq \omega_t\circ f\circ U_t\colon\D\to\C$ parametrizes the interior of the curve $\omega_t(\Gamma)$. The Beltrami coefficient of $\omega_t\circ f$ is given for $z\in\D$, by the following formula 
		\begin{equation*}
			\mu_{\omega_t\circ f}(z) = \nu_t(f(z))\, \frac{\overline{f'(z)}}{f'(z)} = \mu_{U_t^{-1}}(z).
		\end{equation*}
		Therefore, we have $\mu_{f_t}=0$ and thus, $f_t$ is holomorphic. Hence the welding of the curve $\omega_t(\Gamma)$ is given by $g_t^{-1}\circ f_t|_{\s^1} = \vp \circ U_t|_{\s^1} = \vp\circ u_t = \vp_t$. We now compute $\dot{\nu} = \dr_t\nu_t|_{t=0}$. Since $U_0^{-1}=\id_{\D}$, it holds 
		\begin{equation*}
			\dot{\nu}(f(z)) = \left. \frac{\dd \nu_t}{\dd t}\right|_{t=0}\left(f(z)\right) = \dr_{\bz} \left(\left. \frac{\dd U_t^{-1}}{\dd t}\right|_{t=0}\right)(z)\ \frac{ f'(z) }{ \overline{f'(z)} }.
		\end{equation*}
		Since $U_t^{-1}\circ U_t(z)=z$ and $U_0(z)=z$, we obtain 
		\begin{equation}\label{eq:nudot}
			\dot{\nu}(f(z)) = -\dr_{\bz} \left(\left. \frac{\dd U_t}{\dd t}\right|_{t=0}\right)(z)\ \frac{ f'(z) }{ \overline{f'(z)} }.
		\end{equation}
	We have from \cite[Theorem 1.4]{W2019-2} and \cite[Corollary II.3.9]{TT2006}, 
	\begin{equation*}
		\left. \frac{\dd}{\dd t} \right|_{t=0} \IL(\vp_t)   = -\frac{4}{\pi}\, \Re\left[ \int_{\Omega} \dot{\nu}(w)\, \Sr(f^{-1})(w)\, |\dd w|^2 \right].
	\end{equation*} 
	By change of variable $w=f(z)$, we obtain that $|\dd w|^2 =  |f'(z)|^2\, |\dd z|^2$ and $\Sr(f^{-1})(w) = -f(z)^{-2}\,\Sr(f)(z)$. Hence, there holds
	\begin{equation*}
		\left. \frac{\dd}{\dd t} \right|_{t=0} \IL(\vp_t) = \frac{4}{\pi}\, \Re\left[ \int_{\Omega} \dot{\nu}\left(f(z)\right)\, \frac{ \overline{f'(z)} }{ f'(z) }\, \Sr(f)(z)\, |\dd z|^2 \right].
	\end{equation*}
	Using \eqref{eq:nudot}, we obtain, with $\dot{U} = \dr_t U_t|_{t=0}$, 
	\begin{equation*}
		\left. \frac{\dd}{\dd t} \right|_{t=0} \IL(\vp_t)  = -\frac{4}{\pi}\, \Re\left[ \int_{\Omega} \dr_{\bz} \dot{U}(z)\, \Sr(f)(z)\, |\dd z|^2 \right].
	\end{equation*}
	By integration by parts as in \cite[Remark 3.4]{GQW2025}, then using \eqref{eq:decompo-u} together with $U_t|_{\s^1}=u_t$, we obtain 
	\begin{equation*}
		\begin{aligned} 
			\left. \frac{\dd}{\dd t} \right|_{t=0} \IL(\vp_t) &= -\frac{2}{\pi}\, \Re\left[ \frac{1}{\ii}\int_{\Omega} \dr_{\bz} \left( \dot{U}(z)\, \Sr(f)(z)\right)\, \dd \bar{z}\wedge \dd z \right] \\[2mm]
            &= -\frac{2}{\pi}\, \Re\left[ \frac{1}{\ii}\int_{\Omega} \dd\left(  \dot{U}(z)\, \Sr(f)(z)\, \dd z\right)  \right] \\[2mm]
            &  = -\frac{2}{\pi}\,\Re\left[ \frac{1}{\ii} \int_{\s^1} \dot{u}(z)\, \Sr(f)(z)\, \dd z \right] \\[2mm]
			& = -\frac{2}{\pi}\,\Re\left[\frac{1}{\ii} \int_{\s^1} \left(v(z) - z^2\, \overline{v(z)}\right)\, \Sr(f)(z)\, \dd z\right].
		\end{aligned} 
	\end{equation*}
	Since both $v$ and $\Sr(f)$ are holomorphic on $\D$, the first integral vanishes and we obtain
	\begin{equation*}
		\begin{aligned} 
			\left. \frac{\dd}{\dd t} \right|_{t=0} \IL(\vp_t) & = \frac{2}{\pi}\, \Re\left[\frac{1}{\ii} \int_{\s^1}  z^2\, \overline{v(z)}\, \Sr(f)(z)\, \dd z \right] \\[2mm]
			& =  \frac{2}{\pi}\, \Re\left[ \frac{1}{\ii} \int_{\s^1}  z^2\, \overline{v(z)}\, \Sr(f)(z)\, \ii z\, |\dd z| \right] \\[2mm]
			& = \frac{2}{\pi}\, \Re\left[  \int_{\s^1}  z^3\, \overline{v(z)}\, \Sr(f)(z)\,  |\dd z| \right] .
		\end{aligned} 
	\end{equation*}
	\end{proof} 
	Using \eqref{eq:der-logA2}, we obtain
	\begin{equation*}
			\left. \frac{\dd}{\dd t} \right|_{t=0} \IL(\vp_t)  = -24\, \left. \frac{\dd }{\dd t}\right|_{t=0} \log \det(A_{\vp_t} A^*_{\vp_t}) .
	\end{equation*}
	Now for the round circle $\s^1$, we have $\IL(\s^1)=0$ and $\vp=\id$. Hence, we have $W_{\s^1}=\id$ and thus it holds $\log\det(A_{\s^1}A^*_{\s^1})=0$. Since the space of orientation-preserving diffeomorphisms of $\s^1$ is path connected, we deduce that for all smooth $\vp$, there holds
	\begin{equation}\label{eq:IL-logdetA}
		\IL(\vp) = -24\log\det(A_{\vp}A_{\vp}^*).
	\end{equation}
	Moreover, the Loewner energy is continuous in the Weil--Petersson topology. Proposition \ref{pr:Mob-Grunsky} together with Lemma~\ref{lem:logdetA} imply that the right-hand side of \eqref{eq:IL-logdetA} is also continuous in the Weil--Petersson topology (see also Theorem 5.1 in \cite{KSS2025}). Therefore the equality \eqref{eq:IL-logdetA} holds also for Weil--Petersson curves and we conclude the proof of Theorem \ref{th:Loewner-Grunsky}, thanks to Lemma~\ref{lem:logdetA}.

	\subsection{Inequalities between the Möbius energy and the Loewner energy.}\label{sec:Mobius-Loewner}
	In this section, we prove the following result.
	\begin{theorem}\label{thm:compare}
		For any closed Weil--Petersson curve $\Gamma\subset \C$, it holds 
		\begin{equation}\label{eq:ineq1}
			\IL(\Gamma) \leqslant \frac{6}{\pi^2}\left( \Mob(\Gamma)-4 \right).
		\end{equation}
		Moreover, if $\Gamma$ is a $K$-quasicircle, then there exists a constant $C_K>0$ depending only on $K$ such that 
		\begin{equation}\label{eq:ineq2}
			\Mob(\Gamma)-4\leqslant C_K\, \IL(\Gamma).
		\end{equation}
	\end{theorem}
	
	\begin{remark}\label{rk:trace}
		The proof shows that the main difference between $\Mob$ and $\IL$ is the optimal constant of the harmonic extension $\W^{1/2,2}(\Gamma)\to \W^{1,2}(\C)$, rather than the quasiconformal factor itself.
	\end{remark}
	
	\begin{proof} 	We first prove \eqref{eq:ineq1}. By Theorem \ref{th:Loewner-Grunsky}, if we let $(\lambda_j)_{j\in\N}$ denote the eigenvalues of $\left[\mathrm{Gr}^{1/2}_{f}\right]^*\, \mathrm{Gr}^{1/2}_{f}$, then we have 
		\begin{equation}\label{eq:LG}
			\IL(\Gamma) =24 \sum_{j\geqslant 0} \log(1+\lambda_j) \leqslant 24\sum_{j\geqslant 0} \lambda_j = 24\left\| \mathrm{Gr}^{1/2}_{f} \right\|_{\rm HS}^2.
		\end{equation}
		Using Proposition \ref{pr:Mob-Grunsky}, we obtain 
		\begin{equation*}
			\IL(\Gamma) \leqslant \frac{6}{\pi^2}(\Mob(\Gamma)-4).
		\end{equation*}
		We now prove \eqref{eq:ineq2}. Up to inversion, we assume that $\Gamma$ contains $\infty$, and thus, is not closed in $\C$.
	By \eqref{eq:Mobius-tangent}, we have 
	\begin{equation*}
		\Mob(\Gamma) \leqslant \frac{1}{2}\, [\tau]_{\W^{1/2,2}(\Gamma)}^2 + 2\|N_{\Gamma}\|_{\rm HS}^2,
	\end{equation*}
    where $N_\Gamma$ is defined in \eqref{eq:TN}.
	We estimate the norm of $N_{\Gamma}$ as follows. Given $x,y\in\Gamma$, we denote $x\triangleright y$ the shortest path along $\Gamma$ from $x$ to $y$. It holds
	\begin{equation*}
		\begin{aligned}
			\int_{\Gamma\times\Gamma} \frac{|\scal{x-y}{n(x)}|^2}{|x-y|^4}\, |\dd x|\, |\dd y|& = \int_{\Gamma\times \Gamma} \left| \int_{x\triangleright y} \scal{\tau(z)}{n(x)}\, |\dd z|\right|^2\, \frac{|\dd x|\, |\dd y|}{|x-y|^4} \\[2mm]
			& = \int_{\Gamma\times \Gamma} \left| \int_{x\triangleright y} \scal{\tau(z)-\tau(x)}{n(x)}\, |\dd z|\right|^2\, \frac{|\dd x|\, |\dd y|}{|x-y|^4} \\[2mm]
			& \leqslant \int_{\Gamma\times \Gamma}  \dist_{\Gamma}(x,y)\, \left(\int_{x\triangleright y} |\tau(z)-\tau(x)|^2\, |\dd z|\right)\, \frac{|\dd x|\, |\dd y|}{|x-y|^4}.
		\end{aligned}
	\end{equation*}
	Since $\Gamma$ is chord-arc, there exists a constant $C_{\rm ca}>0$ such that
	\begin{equation*}
		\begin{aligned}
		\int_{\Gamma\times\Gamma} \frac{|\scal{x-y}{n(x)}|^2}{|x-y|^4}\, |\dd x|\, |\dd y|& \leqslant C_{\rm ca}\, \int_{\Gamma\times \Gamma}  \left( \int_{x\triangleright y} |\tau(z)-\tau(x)|^2\, |\dd z|\right) \, \frac{|\dd x|\, |\dd y|}{|x-y|^3} \\[2mm]
		& \leqslant C_{\rm ca}\, \int_{\Gamma\times \Gamma}  |\tau(z)-\tau(x)|^2\, \left(\int_{\Gamma} \mathbf{1}_{\{z\in x\triangleright y\}} \frac{|\dd y|}{|x-y|^3} \right)\, |\dd z|\, |\dd x| \\[2mm]
		& \leqslant C_{\rm ca}\, \int_{\Gamma\times \Gamma}  |\tau(z)-\tau(x)|^2\, \left(\int_{\dist_{\Gamma}(x,z)}^{\infty}  \frac{\dd r}{r^3} \right)\, |\dd z|\, |\dd x| \\[2mm]
		& \leqslant C_{\rm ca}\, \int_{\Gamma\times \Gamma}  \frac{ |\tau(z)-\tau(x)|^2 }{ \dist_{\Gamma}(x,z)^2 }\, |\dd z|\, |\dd x|.
		\end{aligned} 
	\end{equation*}
	We obtain a constant $C>0$ depending only on the chord-arc constant of $\Gamma$ such that
	\begin{equation*}
		\Mob(\Gamma)\leqslant C\, [\tau]_{\W^{1/2,2}(\Gamma)}^2.
	\end{equation*}
	Since we assumed $\Gamma$ to pass through $\infty$, we can find a lifting of the tangent vector and write $\tau = \ee^{\ii \sigma}$ for $\sigma\colon \Gamma\to \R$. It holds
	\begin{equation*}
		|\ee^{\ii \sigma(x)} - \ee^{\ii \sigma(y)}| =2\, \left|\sin\left(\frac{\sigma(x) - \sigma(y)}{2}\right)\right| \leqslant |\sigma(x)-\sigma(y)|.
	\end{equation*}
	We then obtain 
	\begin{equation*}
		[\tau]_{\W^{1/2,2}(\Gamma)}^2 \leqslant [\sigma]_{\W^{1/2,2}(\Gamma)}^2.
	\end{equation*}
	Moreover, it follows from \cite{WZ2025} that traces of $\W^{1,2}(\C)$ functions coincides with $\W^{1/2,2}(\Gamma)$ functions if and only if $\Gamma$ is a quasiconformal disk. Therefore, if $K$ is the associated quasiconformal factor and we denote $\widetilde{\sigma}$ the harmonic extension of $\sigma$ to $\C$, then it holds, using \cite[Equation (1.4)]{VW2020}
	\begin{equation*}
		[\sigma]_{\W^{1/2,2}(\Gamma)}^2 \approx_K \|\g \widetilde{\sigma}\|_{L^2(\C)}^2 = \pi\, \IL(\Gamma).
	\end{equation*}
	Hence, we obtain
	\begin{equation*}
		\Mob(\Gamma)\aleq_K \IL(\Gamma).\qedhere
	\end{equation*}
	\end{proof} 
	
	\begin{remark}
		One can also have an alternative comparison using \eqref{eq:LG}. Indeed, since the function $\left(x\mapsto x^{-1}\log(1+x)\right)$ is decreasing, we have 
		\begin{equation*}
			\forall j\in\N,\qquad \frac{\log(1+\lambda_j)}{\lambda_j} \geqslant \frac{\log(1+\lambda_0) }{\lambda_0}.
		\end{equation*}
		Hence, we obtain 
		\begin{equation*}
			\frac{1}{4\pi^2}\left(\Mob(\Gamma)-4\right) = \sum_{j\geqslant 0} \lambda_j \leqslant \frac{\lambda_0}{\log(1+\lambda_0)} \sum_{j\geqslant 0} \log(1+\lambda_j) = \frac{\lambda_0}{24\, \log(1+\lambda_0)}\, \IL(\Gamma).
		\end{equation*}
		The number $\lambda_0$ is the operator norm of the fermionic Grunsky operator $\lambda_0 = \left\| \mathrm{Gr}^{1/2}_{f} \right\|_{\rm op}^2$. It would be interesting to relate $\lambda_0$ with the best constant of the trace inequality, see Remark \ref{rk:trace}. If $\Gamma$ is a $K$-quasicircle, it would also be interesting to relate $\lambda_0$ with $K$.
	\end{remark}

    \section{Optimization problems for M\"obius energy}

	\subsection{Existence of minimizers for $\Pc^1_{\Mob}$ and $\Pc^2_{\Mob}$}\label{sec:wp_attainment}
	
	In this section, we prove the existence of minimizers for $\Pc^1_{\Mob}$ and $\Pc^2_{\Mob}$.

	\begin{proposition}\label{pr:arc-existence}
		Let $\gamma_0\colon \T\to \C$ be an embedding of class $\Cd^{1,\beta}$ for some $\beta>1/2$ and $u,v\in \T$. There exists a Weil--Petersson curve achieving the infimum $\Pc^2_{\Mob}(\gamma_0,u,v)$ defined in \eqref{eq:P2Mob}.
	\end{proposition}
	
	\begin{proof}
		Choose a point $p$ in the interior of the prescribed arc $\bar \Gamma_0\coloneqq \gamma_0([u,v])$, and fix one Möbius transformation $M$ satisfying $M(p)=\infty$.
		Let $(\gamma_k)_{k\in\N}$ be a minimizing sequence for \eqref{eq:P2Mob}.  Since $p$ is contained in every curve $\Gamma_k\coloneqq \gamma_k(\T)$, every curve $M(\Gamma_k)$ contains infinity. Up to this Möbius transformation, we can now parametrize every curve $\gamma_k\colon \R\to \C$ by arclength. We denote $\ell_k$ the length of the optimized part $\Gamma_k\setminus \bar\Gamma_0$.
		
		Let $q_-\coloneqq \gamma_0(u)$ and $q_+\coloneqq \gamma_0(v)$.  We choose the arclength parametrization such that $\gamma_k([0,\ell_k])=\Gamma_k\setminus \bar\Gamma_0$ corresponds to the optimized part, with $\gamma_k(0)=q_-$ and $\gamma_k(\ell_k)=q_+$.
		
		By Lemma \ref{lem:distortion} and using the uniform bound $\Mob(\gamma_k)\leqslant C$ (due to the minimizing property), there is a constant $c>0$, independent of $k$, such that
		\begin{equation}\label{eq:unif-chord-arc}
		\forall s,t\in\R,\qquad |\gamma_k(s)-\gamma_k(t)|\geqslant c|s-t|.
		\end{equation}
		In particular, we obtain $c\, \ell_k \leqslant  |q_--q_+| \leqslant \ell_k$. Up to a subsequence, we obtain that $(\ell_k)_{k\in\N}$ converges to some limit $\ell_{\infty}\in(0,\infty)$.
		
		Therefore, $(\gamma_k)_{k\in\N}$ is uniformly bounded and biLipschitz on $[0,\sup_k\ell_k]$.  By Lemma \ref{lem:tangent-control}, the sequence $(\gamma_k')_{k\in\N}$ is uniformly bounded in $\W^{1/2,2}_{\loc}$. We obtain a limiting curve $\gamma_{\infty}\colon \R\to \C$ such that $\gamma_{\infty}(\R\setminus [0,\ell_{\infty}]) = \bar\Gamma_0$ and
		\begin{equation*}
			\begin{cases} 
				\gamma_k \to \gamma_\infty
			&\text{locally uniformly on }\R,\\[2mm]
			\gamma_k' \to \gamma_\infty'
			& \text{strongly in }\Ll^2_{\mathrm{loc}}(\R),\\[2mm]
			\gamma_k' \to \gamma_\infty'
			& \text{weakly in }\W^{1/2,2}_{\mathrm{loc}}(\R).
			\end{cases} 
		\end{equation*}
		The strong $\Ll^2_{\mathrm{loc}}$ convergence and $|\gamma_k'|=1$ imply that $|\gamma_\infty'|=1$ a.e.\ as well.
		The uniform chord-arc inequality \eqref{eq:unif-chord-arc} passes to the limit by uniform convergence:
		\begin{equation}\label{eq:limit-chord-arc}
			|\gamma_\infty(s)-\gamma_\infty(t)|\geqslant c|s-t|.
		\end{equation}
		Thus $\gamma_\infty$ is injective, proper, parametrized by arc-length with tangent vector in $\W^{1/2,2}([0,\ell_{\infty}])$. Hence, it is a Weil--Petersson curve.
		
		For $s\ne t$, local uniform convergence implies that 
		\begin{equation*} 
		\frac{1}{|\gamma_k(s)-\gamma_k(t)|^2}-\frac{1}{|s-t|^2}
		\xrightarrow[k\to \infty]{\text{a.e.}}{}
		\frac{1}{|\gamma_\infty(s)-\gamma_\infty(t)|^2}-\frac{1}{|s-t|^2}.
		\end{equation*} 
		All these functions are nonnegative. We deduce from Fatou's lemma that 
		\begin{equation*} 
			\Mob(\gamma_\infty) \leqslant \liminf_{k\to\infty}\Mob(\gamma_k).
		\end{equation*} 
		Consequently, the limit $\gamma_{\infty}$ is a minimizer for $\Pc^2_{\Mob}(\gamma_0,u,v)$.
	\end{proof}

	We now turn to $\Pc^1_{\Mob}$, for which the proof is very similar. 
	
	\begin{proposition}
		Consider $n\geqslant 2$ distinct points $z_1,\ldots,z_n\in\C$ and let $Z\coloneqq (z_1,\ldots,z_n)$. Let $\tau$ be a relative isotopy class of curves passing through $Z$. Then, there exists a Weil--Petersson curve achieving $\Pc^1_{\Mob}(Z,\tau)$ defined in \eqref{eq:P1Mob}.
	\end{proposition}
	
	\begin{proof}
		If $n\leqslant 3$, a generalized circle contains the prescribed points, which is an absolute minimizer.
		
		Assume now $n\geqslant 4$.  Fix a Möbius transformation $M$ such that $M(z_1)=\infty$ and denote $w_i\coloneqq M(z_i)\in\C$ for $i\in\{2,\ldots,n\}$.
		Let $(\gamma_k)_{k\in\N}$ be a minimizing sequence. Up to exchanging $\gamma_k$ with $M(\gamma_k)$, we can assume that they all pass through $\infty$. Hence, we can consider an arclength parametrization $\gamma_k\colon \R\to\C$. We choose the orientation and the origin of the parametrization so that $\gamma_k(0)=w_2$ and $\gamma_k(s_{3,k})=w_3$ with $s_{3,k}>0$. For every $i\ge2$, let $s_{i,k}$ be the unique parameter satisfying $\gamma_k(s_{i,k})=w_i$.
		
		The energies $\Mob(\gamma_k)$ are uniformly bounded. Hence, Lemma~\ref{lem:distortion} implies that there is a uniform constant $c>0$ such that $|\gamma_k(s)-\gamma_k(t)|\geqslant c|s-t|$. Taking $t=0$ and $s=s_{i,k}$ gives
		\begin{equation*} 
		      c|s_{i,k}| \leqslant |w_i-w_2|\leqslant |s_{i,k}|.
		\end{equation*}
		Up to taking subsequences, we can assume that, for any $i\in\{2,\ldots,n\}$, the sequence $(s_{i,k})_{k\in\N}$ converges to some $s_{i,\infty}\in \R$. Moreover, we have $s_{i,\infty}\neq s_{j,\infty}$ for $i\neq j$. We conclude by following step by step the proof of Proposition \ref{pr:arc-existence}.
	\end{proof}

	\subsection{Regularity of solutions to $\Pc^2_{\Mob}$}\label{sec:reg-P2mob}
	
	In this section, we prove the following result.
	\begin{theorem}\label{th:regu_P2}
        Let $\gamma_0\colon \T\to \C$ be an embedding of class $\Cd^{1,\beta}$ for some $\beta>1/2$ and $u,v\in \T$.
		Consider $\eta$ a solution to $\Pc^2_{\Mob}(\gamma_0,u,v)$ defined in \eqref{eq:P2Mob}. Then it holds $\eta\in \Cd^{1,1/2}(\T_{\ell};\R^2)$, where $\ell$ is the length of $\eta$ and $\eta$ is parameterized by arclength. The optimized part is analytic in its interior. 
	\end{theorem}

    We assumed that $\eta$ does not pass through $\infty$. Otherwise, we can compose $\gamma_0$ and $\eta$ by a Möbius transformation such that $\eta$ does not pass through $\infty$.

    The proof consists of three steps. We recall the Euler--Lagrange equation, proved by Blatt--Reiter--Schikorra \cite{BRS2016} in the Weil--Petersson setting, in Lemma \ref{lm:EL-Mobius}, together with its $\eps$-regularity property in \eqref{eq:eps-reg}. Adapting the arguments from Hildebrandt--Kaul \cite{HK1972}, we prove that $\eta$ is of class $\Cd^1$. This allows us to rewrite the Euler--Lagrange equation in Lemma \ref{lm:EL-Mob}. In a second step, we prove that $\eta$ is of class $\Cd^{1,\alpha}$ for some $\alpha>0$ by a hole-filling trick in Lemma \ref{lm:C1alpha_reg}. Plugging this into the equation, we obtain that $\eta$ is of class $\W^{t,2}$ for some $t>1/2$ in Lemma \ref{lm:Sobolev-reg}. 
    This threshold allows us to initiate the bootstrap argument in \cite{BRS2016}. We obtain finally the $\Cd^{1,1/2}$ regularity in Lemma \ref{lm:Chalf_minus}.
	
	In the remainder of the section, we will keep the notations $\gamma_0,u,v$ of Theorem \ref{th:regu_P2}. 
	\begin{lemma}\label{lm:EL-Mobius}
        With the notations of Theorem \ref{th:regu_P2}, any solution $\eta$ to $\Pc^2_{\Mob}(\gamma_0,u,v)$ satisfies the following Euler--Lagrange equation: for any test function $\vp\in \Cd^{\infty}_c(\T_{\ell}\setminus (u,v))$,
		\begin{equation*}
			0 = \int_{\T_{\ell}\times \T_{\ell}} \left( \scal{\eta'(s)}{\vp'(s)} - \frac{\scal{\eta(s)-\eta(t)}{\vp(s)-\vp(t)}}{|\eta(s)-\eta(t)|^2} \right)\, \frac{\dd s\, \dd t}{|\eta(s)-\eta(t)|^2}.
		\end{equation*}
		In particular, if $\eta$ is regular enough, we have 
		\begin{equation}\label{eq:EL-Mob}
			\forall s\in \T_{\ell}\setminus (u,v),\qquad \int_{\T_{\ell}} \left[ 2\frac{\proj_{\eta'(s)^{\perp}}\left[ \eta(s)-\eta(t) \right] }{|\eta(s)-\eta(t)|^4} +   \frac{\eta''(s)}{|\eta(s)-\eta(t)|^2}\right] \, \dd t =0.
		\end{equation}
	\end{lemma}
	
	\begin{proof}
		The first statement is a direct application of \cite[Theorem II]{BRS2016}.\\
		The second statement is an integration by parts, see for instance \cite[Lemma 3.7.3]{ohara2003}. 
	\end{proof}
	
	Let $\eta$ be a solution to $\Pc^2_{\Mob}(\gamma_0,u,v)$ (with the notations of Theorem \ref{th:regu_P2}) parameterized by arc-length with length $\ell = 2\pi$ (this can always be assumed up to dilation). We consider $\eta$ as a $2\pi$-periodic map $\eta\colon \R\to \C$, with $[u,v]\subset (-\pi,\pi)$. Hence $\eta = \gamma_0$ on $\R\setminus ([u,v]+2\pi\Z)$.
	The arguments from \cite{BRS2016} still hold to prove interior regularity estimates. We conclude that $\eta \in \Cd^{\infty}(\T\setminus \{u,v\};\C)$ and verify the following $\eps$-regularity estimate: there exist $\eps_0>0$, $C>0$ and $\Lambda>1$ such that if $(x-\Lambda r,x+\Lambda r)\Subset (u,v)$ with $[\eta']_{\W^{1/2,2}(x-\Lambda r,x+\Lambda r)}\leqslant \eps_0$, then it holds
	\begin{equation}\label{eq:eps-reg}
		r\, \|\eta''\|_{\Ll^{\infty}\left( (x-r,x+r) \right) } \leqslant C\, \Mob(\eta).
	\end{equation}
    One can also adapt the proof of Corollary 1.3 \cite{BV2019} by introducing cutoff functions, so that the optimized part of $\eta$ is analytic in the interior. We refer to \cite{B2022} for a similar argument. The main idea is to show that the $\Ll^{\infty}$ norm of a $k$-th derivative on an interval of size $\rho$ is controlled by $k!/r^k$ for some $r>0$ depending only on $\rho$. To do so, we need to perform a careful analysis of the estimates of high-order derivatives, using for instance, Faà di Bruno's formula. Since the dominant term can be localized in a similar manner as in \cite{BRS2016}, the approach applies also to our solution $\eta$.
	
   The rest of the section aims to prove the global $\Cd^{1,1/2}$ regularity.  By adapting the arguments of \cite{HK1972} (see also \cite[Section 7.2]{S2011}), we now show that the tangent vector is continuous through the singularities.
	\begin{lemma}\label{lm:C1_reg}
		Let $\eta$ be as in Theorem \ref{th:regu_P2}. Then $\eta$ is globally $\Cd^1$ in arc-length parametrization.
	\end{lemma}
	\begin{proof} 
		Let $\bar{x}\in(u,v)$. Without loss of generality, we can assume that $\tau\coloneqq |\bar{x}-v|\leqslant |\bar{x}-u|$. Then we have $B_{\tau}(\bar{x})\subset (u,v)$.
		Let $\mu\in(0,\Lambda^{-1})$, where $\Lambda$ is defined above \eqref{eq:eps-reg}, and denote $\sigma\coloneqq \mu\, \tau$. We have 
		\begin{equation}\label{eq:est_average}
            \begin{aligned} 
			& \int_{B_{\sigma}(\bar{x})} \left| \eta'(s) - \gamma_0'(2v-s) \right|\, \dd s \\
            & \leqslant \int_{B_{\sigma}(\bar{x})} \left| \eta'(s) - \fint_{B_{2\tau}(v)} \eta' \right|\, \dd s + \int_{B_{\sigma}(\bar{x})} \left| \eta'(2v-s) - \fint_{B_{2\tau}(v)} \eta' \right|\, \dd s  \\[2mm]
			& \leqslant C\, \tau\, \Mr(2\tau), 
            \end{aligned}
		\end{equation}
		where for any $r>0$, we have denoted
		\begin{equation*}
			\Mr(r) \coloneqq \sup_{x\in\R} \sup_{0<\rho<r} \fint_{B_{\rho}(x)} \left| \eta'(s) - \fint_{B_{\rho}(x)} \eta' \right|\, \dd s.
		\end{equation*}
		Since $\eta'\in \VMO(\T)$, we have 
		\begin{equation*}
			\Mr(r) \xrightarrow[r\to 0]{} 0.
		\end{equation*}
		Thanks to \eqref{eq:est_average}, there exists a measurable set $E_{\sigma,\bar{x}}\subset B_{\sigma}(\bar{x})$ of positive measure such that 
		\begin{equation}\label{eq:good_set}
			\forall y\in E_{\sigma,\bar{x}},\qquad \left|\eta'(y) - \gamma_0'(2v-y)\right| \leqslant C\, \frac{\tau}{\sigma}\, \Mr(2\tau) = \frac{C}{\mu}\, \Mr(2\tau).
		\end{equation}
		We obtain for any $y\in E_{\sigma,\bar{x}}$ that
		\begin{equation*}
			\left| \eta'(\bar{x}) - \gamma_0'(v)\right| \leqslant \left| \eta'(\bar{x}) - \eta'(y) \right| + \left| \eta'(y) - \gamma_0'(2v-y) \right| + \left| \gamma_0'(2v-y) - \gamma_0'(v) \right| .
		\end{equation*}
		Thanks to \eqref{eq:eps-reg} applied on the ball $B_{\Lambda^{-1}\tau}(\bar{x})$ and \eqref{eq:good_set}, we obtain 
		\begin{equation}
			\left| \eta'(\bar{x}) - \gamma_0'(v)\right|  \leqslant C\left( \frac{|\bar{x}-y| }{\tau} + \frac{1}{\mu}\, \Mr(2\tau) + \tau \right)  \leqslant C\left( \mu + \frac{1}{\mu}\, \Mr(2\tau) + \tau \right) . \label{eq:pointwise_est}
		\end{equation}
		Choosing first $\mu$ small, then $\tau$ small, we obtain that for any $\eps>0$, there exists $\delta>0$ such that for $\bar{x}\in(v-\delta,v)$, it holds
		\[
		\left| \eta'(\bar{x}) - \gamma_0'(v)\right| \leqslant \eps.
		\]
		In other words, we have $\eta'\in \Cd^0(\T)$. 
	\end{proof} 
    In order to improve the regularity, we first rewrite the Euler--Lagrange equation following \cite{BRS2016}.
    \begin{lemma}\label{lm:EL-Mob}
        Let $\eta$ be as in Theorem \ref{th:regu_P2}. We denote 
	\begin{equation*}
		h(s) \coloneqq \int_{-\pi}^{\pi} \frac{|\eta'(s+w)-\eta'(s)|^2}{|w|^2}\, \dd w.
	\end{equation*}
    Then \eqref{eq:EL-Mob} can be written as follows. For any $\vp\in \Cd^{\infty}_c(\T\setminus (u,v))$,
    \begin{equation*}
        \int_{\T} \scal{|D|^{1/2} \eta' }{ |D|^{1/2}\vp' } = T_1(\eta,\vp) + T_2(\eta,\vp),
    \end{equation*}
    where there are two functions $f_1,f_2\colon \T\to \C$ satisfying the pointwise inequalities $|f_1|+|f_2|\aleq h$ together with 
    \begin{equation*}
        \forall i\in\{1,2\},\qquad T_i(\eta,\vp) = \int_{\T} \scal{f_i}{\vp'}.
    \end{equation*}
    \end{lemma}

    \begin{remark}
        The property $\eta'\in \W^{1/2,2}(\T)$ implies that $h\in \Ll^1(\T)$.
    \end{remark}
    
    \begin{proof}
        We have 
	\begin{equation*}
		\begin{aligned} 
			& \int_{\T} \scal{|D|^{1/2}\eta'}{|D|^{1/2} \vp'} \\[2mm]
            & = \int_{\T\times \T} \scal{\eta'(s)}{\vp'(s)}\left(\frac{1}{|\eta(s)-\eta(t)|^2} - \frac{1}{|s-t|^2}\right)\, \dd s\, \dd t\\[2mm]
			& \quad - \int_{\T\times \T} \scal{\eta(s)-\eta(t)}{\vp(s)-\vp(t)}  \left( \frac{ 1 }{|\eta(s)-\eta(t)|^4} - \frac{1}{|s-t|^4} \right)\, \dd s\, \dd t \\[2mm]
			& = \int_{\T\times \T} \frac{\scal{\eta'(s)}{\vp'(s)}}{|s-t|^2}\left(\frac{|s-t|^2}{|\eta(s)-\eta(t)|^2} - 1\right)\, \dd s\, \dd t\\[2mm]
			& \quad - \int_{\T\times \T} \frac{ \scal{\eta(s)-\eta(t)}{\vp(s)-\vp(t)} }{|s-t|^4} \left( \frac{ |s-t|^4 }{|\eta(s)-\eta(t)|^4} - 1 \right)\, \dd s\, \dd t.
        \end{aligned} 
        \end{equation*}
        We obtain 
        \begin{equation*}
            \begin{aligned} 
            & \int_{\T} \scal{|D|^{1/2}\eta'}{|D|^{1/2} \vp'} \\[2mm]
			& = \int_{[-\pi,\pi]^2} \frac{\scal{\eta'(s)}{\vp'(s)}}{|w|^2}\left(\frac{|w|^2}{|\eta(s)-\eta(s+w)|^2} - 1\right)\, \dd s\, \dd w\\[2mm]
			& \quad - \int_{[-\pi,\pi]^2} \frac{ \scal{\eta(s+w)-\eta(s)}{\vp(s+w)-\vp(s)} }{|w|^4} \left( \frac{ |w|^4 }{|\eta(s)-\eta(s+w)|^4} - 1 \right)\, \dd s\, \dd w.
		\end{aligned} 
	\end{equation*}
    We denote the first term $T_1(\eta,\vp)$ and the second term $T_2(\eta,\vp)$.
	There holds 
	\begin{equation*}
		\begin{aligned}
			\frac{|\eta(s)-\eta(s+w)|^2}{w^2} & = \left|\int_0^1 \eta'(s+xw)\, \dd x\right|^2 \\[2mm]
			& = \int_0^1\int_0^1 \scal{\eta'(s+x_1w)}{\eta'(s+x_2w)}\, \dd x_1\, \dd x_2 \\[2mm]
			& = \frac{1}{2}\int_0^1\int_0^1 \big( 2- |\eta'(s+x_1w) - \eta'(s+x_2w)|^2 \big)\, \dd x_1\, \dd x_2 \\[2mm]
			& = 1- \frac{1}{2}\int_0^1\int_0^1  |\eta'(s+x_1w) - \eta'(s+x_2w)|^2\, \dd x_1\, \dd x_2.
		\end{aligned}
	\end{equation*}
	Therefore, we obtain 
	\begin{equation*}
		\begin{aligned}
			\frac{w^2}{|\eta(s)-\eta(s+w)|^2} - 1 & = \frac{w^2}{2\, |\eta(s)-\eta(s+w)|^2}\int_0^1\int_0^1  |\eta'(s+x_1w) - \eta'(s+x_2w)|^2\, \dd x_1\, \dd x_2
		\end{aligned}
	\end{equation*}
	We also have 
	\begin{equation*}
		\begin{aligned}
			& \frac{w^4}{|\eta(s)-\eta(s+w)|^4} - 1 \\[2mm]
            & = \frac{w^4}{|\eta(s)-\eta(s+w)|^4}\left( 1 - \frac{|\eta(s)-\eta(s+w)|^4}{w^4} \right) \\[2mm]
			& = \frac{w^4}{|\eta(s)-\eta(s+w)|^4}\left( 1 + \frac{|\eta(s)-\eta(s+w)|^2}{w^2} \right) \left( 1 - \frac{|\eta(s)-\eta(s+w)|^2}{w^2} \right) \\[2mm]
			&= \frac{w^4}{2\, |\eta(s)-\eta(s+w)|^4}\left( 1 + \frac{|\eta(s)-\eta(s+w)|^2}{w^2} \right)\int_0^1\int_0^1  |\eta'(s+x_1w) - \eta'(s+x_2w)|^2\, \dd x_1\, \dd x_2.
		\end{aligned}
	\end{equation*}

	We first focus on $T_1$. We have
	\begin{equation*}
		\begin{aligned} 
			T_1(\eta,\vp) = \int_{-\pi}^{\pi} \!\! \scal{\vp'(s)}{ \int_{-\pi}^{\pi} \frac{\eta'(s)}{2\, |\eta(s)-\eta(s+w)|^2}\int_0^1\!\int_0^1 \! |\eta'(s+x_1w) - \eta'(s+x_2w)|^2\, \dd x_1\, \dd x_2\, \dd w }\, \dd s.
		\end{aligned} 
	\end{equation*}
	The second term in the brackets verifies the pointwise estimate
	\begin{equation*}
		\begin{aligned} 
			& \left| \int_{-\pi}^{\pi} \frac{\eta'(s)}{2\, |\eta(s)-\eta(s+w)|^2}\int_0^1\int_0^1  |\eta'(s+x_1w) - \eta'(s+x_2w)|^2\, \dd x_1\, \dd x_2\, \dd w \right| \\[2mm]
			&  \aleq \int_{-\pi}^{\pi} \frac{1}{w^2} \int_{-\pi}^{\pi} |\eta'(s+x w)-\eta'(s)|^2\, \dd x\, \dd w \\[2mm]
			& \aleq \int_{-\pi}^{\pi} \left( \int_{-\pi}^{\pi} \frac{|\eta'(s+xw)-\eta'(s)|^2}{w^2}\, \dd w \right)\, \dd x \\[2mm]
			& \aleq \int_{-\pi}^{\pi} \left( \int_{-\pi}^{\pi} \frac{|\eta'(s+xw)-\eta'(s)|^2}{w^2}\, \dd w \right)\, \dd x.
		\end{aligned} 
	\end{equation*}
	With the change of variable $\widetilde{w} = xw$, $w=\widetilde{w}/x$, we obtain 
	\begin{equation*}
		\begin{aligned}
			& \left| \int_{-\pi}^{\pi} \frac{\eta'(s)}{2\, |\eta(s)-\eta(s+w)|^2}\int_0^1\int_0^1  |\eta'(s+x_1w) - \eta'(s+x_2w)|^2\, \dd x_1\, \dd x_2\, \dd w \right| \\[2mm]
			& \aleq \int_0^1 \left( \int_{-x\pi}^{x\pi} \frac{|\eta'(s+\widetilde{w})-\eta'(s)|^2}{\widetilde{w}^2\, x^{-2}}\, \frac{\dd \widetilde{w}}{x} \right)\, \dd x \\[2mm]
			& \aleq \left( \int_0^1 x\, \dd x\right) \left( \int_{-\pi}^{\pi} \frac{|\eta'(s+\widetilde{w})-\eta'(s)|^2}{\widetilde{w}^2}\, \dd \widetilde{w}\right).
		\end{aligned}
	\end{equation*}
	
	We now focus on $T_2$.	There holds
		\begin{equation*}
			\begin{aligned} 
				& \int_{[-\pi,\pi]^2} \frac{ \scal{\eta(s+w)-\eta(s)}{\vp(s+w)-\vp(s)} }{|w|^4} \left( \frac{ |w|^4 }{|\eta(s)-\eta(s+w)|^4} - 1 \right)\, \dd s\, \dd w \\[2mm]
				& = \int_{[-\pi,\pi]^2} \int_{[0,1]^2} \frac{w^2\, \scal{\eta'(s+t_1w)}{\vp'(s+t_2w)}}{2\, |\eta(s)-\eta(s+w)|^4}\left( 1 + \frac{|\eta(s)-\eta(s+w)|^2}{w^2} \right)\\[2mm]
				&\quad \times \left(\int_0^1\int_0^1  |\eta'(s+x_1w) - \eta'(s+x_2w)|^2\, \dd x_1\, \dd x_2\right) \dd t_1\, \dd t_2\, \dd s\, \dd w.
			\end{aligned} 
		\end{equation*}
		We make the change of variable $u=s+t_2w$, $s=u-t_2w$. We obtain 
		\begin{equation*}
			\begin{aligned} 
				& \int_{[-\pi,\pi]^2} \frac{ \scal{\eta(s+w)-\eta(s)}{\vp(s+w)-\vp(s)} }{|w|^4} \left( \frac{ |w|^4 }{|\eta(s)-\eta(s+w)|^4} - 1 \right)\, \dd s\, \dd w \\[2mm]
				& = \int_{[-\pi,\pi]^2} \int_{[0,1]^2} \frac{w^2\, \scal{\eta'(u+(t_1-t_2)w)}{\vp'(u)}}{2\, |\eta(u-t_2w)-\eta(u+(1-t_2)w)|^4}\left( 1 + \frac{|\eta(u-t_2w)-\eta(u+(1-t_2)w)|^2}{w^2} \right)\\[2mm]
				&\quad \times \left(\int_0^1\int_0^1  |\eta'(u+(x_1-t_2)w) - \eta'(u+(x_2-t_2)w)|^2\, \dd x_1\, \dd x_2\right) \dd t_1\, \dd t_2\, \dd u\, \dd w.
			\end{aligned} 
		\end{equation*}
		We obtain the pointwise estimate in a similar manner as for $T_1$, using the following two estimates
        \begin{equation*}
            \frac{|\eta(u-t_2w)-\eta(u+(1-t_2)w)|^2}{w^2} \aleq 1,
            \qquad \text{ and } \qquad
            \frac{w^2}{|\eta(u-t_2w)-\eta(u+(1-t_2)w)|^4} \aleq w^{-2}.
        \end{equation*}
    \end{proof}
    We now improve the $\Cd^1$ regularity to $\Cd^{1,\alpha}$ regularity for some $\alpha>0$.
    \begin{lemma}\label{lm:C1alpha_reg}
		Let $\eta$ be as in Theorem \ref{th:regu_P2}. Then $\eta$ is globally $\Cd^{1,\alpha}$ in arc-length parametrization for some $\alpha>0$.
	\end{lemma}
    \begin{proof}
    Since $\gamma_0\in \Cd^{1,\beta}$ for some $\beta>1/2$, there exists  $\sigma>0$ such that $\gamma_0\in \W^{3/2+\sigma,2}(\T)$. Hence, we have 
	\begin{equation*}
		\begin{aligned} 
		\int_{\T} \scal{|D|^{1/2}(\eta-\gamma_0)'}{|D|^{1/2} \vp'} & = \int_{\T}\scal{\vp'}{f_1+f_2} - \int_{\T} \scal{|D|^{1/2}\gamma_0'}{|D|^{1/2} \vp'} \\[2mm]
		& = \int_{\T}\scal{\vp'}{f_1+f_2} - \int_{\T} \scal{|D|^{1/2+\sigma}\gamma_0' }{ |D|^{1/2-\sigma} \vp'}.
		\end{aligned} 
	\end{equation*}
	We now take a cutoff function $\chi$ equal to $1$ in $B_{\rho}(u)$ and $0$ on $\T\setminus B_{2 \rho}(u)$, with $\rho^i |\g^i \chi|\aleq 1$. We test with $\vp = \chi^2(\eta-\gamma_0)$. There holds for $0<s<1$
	\begin{equation*}
		\begin{aligned} 
			|D|^s \vp' & = |D|^s\big( 2\chi\chi'(\eta-\gamma_0) + \chi^2\, (\eta-\gamma_0)' \big) \\[2mm]
			& = \left(|D|^s (2\chi\chi')\right)(\eta-\gamma_0) + 2\chi\chi'\left( |D|^s (\eta-\gamma_0)\right) + H_s\left(2\chi\chi',\eta-\gamma_0\right) \\[2mm]
			&\quad + \left(|D|^s (\chi^2)\right)(\eta-\gamma_0)' + \chi^2 \left( |D|^s (\eta-\gamma_0)'\right) + H_s \left(\chi^2,(\eta-\gamma_0)'\right).
		\end{aligned} 
	\end{equation*}
    In the above formula, we have denoted
    \begin{equation*}
        H_s(a,b) \coloneqq |D|^s(ab) - a\, |D|^s b - b\, |D|^s a.
    \end{equation*}
	Using \cite[Lemma A.7]{BRS2016} together with $\left\| |D|^{1/2-\sigma}(2\chi\chi')\right\|_{\Ll^{\infty}(\T)} \aleq \rho^{\sigma-3/2}$ obtained by interpolation inequalities, we obtain  
	\begin{equation*}
		\left\| |D|^{1/2-\sigma} \vp' \right\|_{\Ll^2(\T)} \aleq_{\Mob(\eta),\|(\eta-\gamma_0)'\|_{\Ll^{\infty}} } \rho^{\sigma}. 
	\end{equation*}
	We obtain
	\begin{equation}\label{eq:hole-filling}
		\begin{aligned} 
			&\int_{\T} \scal{ |D|^{1/2}(\eta-\gamma_0)' }{ |D|^{1/2}\left[ \chi^2 (\eta-\gamma_0)' \right] } \\[2mm]
            & \leqslant C\, \left\| (\eta-\gamma_0)' \right\|_{\Ll^{\infty}(B_{2\rho}(u))} \int_{B_{2 \rho}(u)} |h(s)|\,\dd s + C\rho^{\sigma}  \\[2mm]
            &\quad + \int_{\T} \scal{ |D|^{1/2}(\eta-\gamma_0)' }{ |D|^{1/2}\left[ 2\chi\chi' (\eta-\gamma_0) \right] } .
		\end{aligned} 
	\end{equation}
    We first estimate the left-hand side of \eqref{eq:hole-filling}.
	\begin{claim}
		There holds
		\begin{align*}
		  \left[ (\eta-\gamma_0)' \right]_{\W^{1/2,2}(B_{\rho}(u))}^2 & \aleq \int_{\T} \scal{ |D|^{1/2}(\eta-\gamma_0)' }{ |D|^{1/2}\left[ \chi^2 (\eta-\gamma_0)' \right] } \\[2mm]
            & \quad + \rho^{-2}\int_{\T\setminus B_{\rho}(u)} |(\eta-\gamma_0)'|^2.  
		\end{align*}
	\end{claim}
	\begin{proof}
		We know from \cite[Theorem 1.5]{RS2016} that there exists $K\colon \T\times\T\to (0,\infty)$ such that $K(x,y)\approx |x-y|^{-2}$ and for any $a,b\colon \T\to\R$, there holds
		\begin{equation*}
			\int_{\T} \left( |D|^{1/2} a\right)\left( |D|^{1/2} b\right) = \int_{\T} (a(x)-a(y))\, (b(x)-b(y))\, K(x,y)\, |\dd x|\, |\dd y|.
		\end{equation*}
		Hence, we have
		\begin{equation*}
			\begin{aligned} 
				& \int_{\T} \scal{ |D|^{1/2}(\eta-\gamma_0)' }{ |D|^{1/2}\left[ \chi^2 (\eta-\gamma_0)' \right] } \\[2mm]
				& = \int_{\T\times\T} \scal{ (\eta-\gamma_0)'(x) - (\eta-\gamma_0)'(y) }{ \chi(x)^2 (\eta-\gamma_0)'(x) - \chi(y)^2\, (\eta-\gamma_0)'(y)}\, K(x,y)\, |\dd x|\, |\dd y|\\[2mm]
				& = \int_{\T\times\T} \Big( \chi(x)^2 |(\eta-\gamma_0)'(x)|^2 +\chi(y)^2 |(\eta-\gamma_0)'(y)|^2 \\[2mm]
                &\qquad - \left( \chi(x)^2 + \chi(y)^2\right) \scal{ (\eta-\gamma_0)'(x) }{ (\eta-\gamma_0)'(y) } \Big)\, K(x,y)\, |\dd x|\, |\dd y|\\[2mm]
				& = \int_{\T\times\T} \Big( \left| \chi(x)(\eta-\gamma_0)'(x)-\chi(y)(\eta-\gamma_0)'(y) \right|^2  \\[2mm]
                &\qquad - \left( \chi(x)^2 -2\chi(x)\chi(y) + \chi(y)^2\right) \scal{ (\eta-\gamma_0)'(x) }{ (\eta-\gamma_0)'(y) } \Big)\, K(x,y)\, |\dd x|\, |\dd y|\\[2mm]
				& = \int_{\T\times\T} \Big( \left| \chi(x)(\eta-\gamma_0)'(x)-\chi(y)(\eta-\gamma_0)'(y) \right|^2 \\[2mm]
                &\qquad - \left( \chi(x)-\chi(y)\right)^2 \scal{ (\eta-\gamma_0)'(x) }{ (\eta-\gamma_0)'(y) } \Big)\, K(x,y)\, |\dd x|\, |\dd y|.
			\end{aligned} 
		\end{equation*}
		The first term satisfies
		\begin{equation*}
			\begin{aligned} 
				& \int_{\T\times\T}  \left| \chi(x)(\eta-\gamma_0)'(x)-\chi(y)(\eta-\gamma_0)'(y) \right|^2\, K(x,y)\, |\dd x|\, |\dd y|\\[2mm]
				& \geqslant \int_{ B_{\rho}(u)\times B_{\rho}(u)}  \left| (\eta-\gamma_0)'(x)-(\eta-\gamma_0)'(y) \right|^2\, K(x,y)\, |\dd x|\, |\dd y|\\[2mm]
				& \gtrsim \int_{ B_{\rho}(u)\times B_{\rho}(u)} \frac{ \left| (\eta-\gamma_0)'(x)-(\eta-\gamma_0)'(y) \right|^2 }{ |x-y|^2} \, |\dd x|\, |\dd y|= \left[ (\eta-\gamma_0)' \right]_{\W^{1/2,2}(B_{\rho}(u))}^2.
			\end{aligned} 
		\end{equation*}
		The second term verifies 
		\begin{equation*}
			\begin{aligned}
				& \int_{\T\times \T} \left( \chi(x)-\chi(y)\right)^2 \scal{ (\eta-\gamma_0)'(x) }{ (\eta-\gamma_0)'(y) } \, K(x,y)\, \dd x\, |\dd y|\\[2mm]
				& \leqslant\int_{\left( \T\times \T\right)\setminus \left(B_{\rho}(u) \times B_{\rho}(u)\right)} \left( \chi(x)-\chi(y)\right)^2 \left| (\eta-\gamma_0)'(x) \right|\, \left|(\eta-\gamma_0)'(y)\right| \, K(x,y)\, \dd x\, |\dd y|\\[2mm]
				& \aleq \int_{\left( \T\times \T\right)\setminus \left(B_{\rho}(u)\times B_{\rho}(u)\right)} \frac{ \left( \chi(x)-\chi(y)\right)^2 }{|x-y|^2} \left| (\eta-\gamma_0)'(x) \right|\, \left|(\eta-\gamma_0)'(y)\right| \, \dd x\, |\dd y|.
            \end{aligned} 
        \end{equation*}
        Using that $|\chi'|\aleq \rho^{-1}$ and Hölder inequality, we obtain 
        \begin{equation*}
            \begin{aligned} 
                & \int_{\T\times \T} \left( \chi(x)-\chi(y)\right)^2 \scal{ (\eta-\gamma_0)'(x) }{ (\eta-\gamma_0)'(y) } \, K(x,y)\, \dd x\, |\dd y|\\[2mm]
				& \aleq \rho^{-2}\int_{\T\setminus B_{\rho}(u) }\left|(\eta-\gamma_0)'\right|^2.
			\end{aligned}
		\end{equation*}
	\end{proof}
	We now estimate the right-hand side of \eqref{eq:hole-filling}.
	\begin{claim}
		There holds
		\begin{align*}
			&\left|\int_{\T} \scal{ |D|^{1/2}(\eta-\gamma_0)' }{ |D|^{1/2}\left[ 2\chi\chi' (\eta-\gamma_0) \right] }  \right| 
            \\
            & \aleq \rho^{-1}\, \int_{B_{3 \rho}(u)} |(\eta-\gamma_0)'|^2  + \rho \left( \int_{\T\setminus B_{2\rho}(u)} \frac{|(\eta-\gamma_0)'(y)|}{|y-u|^{3/2}}|\dd y|\right)^2.
		\end{align*}
	\end{claim}
	\begin{proof} 
		We introduce another cutoff function $\zeta$ with $\zeta=1$ in $B_{2 \rho}(u)$ and $\zeta=0$ in $\T\setminus B_{3 \rho}(u)$. We decompose
		\begin{equation*}
			(\eta-\gamma_0)' = \zeta(\eta-\gamma_0)' + (1-\zeta) (\eta-\gamma_0)' \eqqcolon v_0+v_1.
		\end{equation*}
		We integrate by parts:
		\begin{equation*}
			\begin{aligned} 
				\int_{\T} \scal{ |D|^{1/2}(\eta-\gamma_0)' }{ |D|^{1/2}\left[ 2\chi\chi' (\eta-\gamma_0) \right] } 
				& = \int_{\T} \scal{ (\eta-\gamma_0)' }{ |D|^1\left[ 2\chi\chi' (\eta-\gamma_0) \right] } \\[2mm]
				& = \int_{\T} \scal{ v_0+v_1 }{ |D|^1\left[ 2\chi\chi' (\eta-\gamma_0) \right] }.
			\end{aligned} 
		\end{equation*}
		Concerning the first term, we have 
		\begin{equation}\label{eq:der-chi-1}
			\begin{aligned} 
				\left| \int_{\T} \scal{ v_0 }{ |D|^1\left(2\chi\chi'(\eta-\gamma_0)\right) }\right|
				& \leqslant\|v_0\|_{\Ll^2(\T)} \left\| |D|^1\left(2\chi\chi'(\eta-\gamma_0)\right) \right\|_{\Ll^2(\T)}\\[2mm]
				& \aleq \|v_0\|_{\Ll^2(\T)} \left\| \left( 2\chi\chi'(\eta-\gamma_0)\right)' \right\|_{\Ll^2(\T)}\\[2mm]
				&\aleq \rho^{-1} \|v_0\|_{\Ll^2(\T)}^2  + \rho^{-2}\|\eta-\gamma_0\|_{\Ll^2(B_{2 \rho}(u))}\\[2mm]
				& \aleq \rho^{-1}\|(\eta-\gamma_0)'\|_{\Ll^2(B_{3\rho}(u))}^2.
			\end{aligned} 
		\end{equation}
		For the last estimate, we used that $(\eta-\gamma_0)'=v_0$ on $\supp(\chi)$.
	
	For $x\in B_{\rho}(u)$ and $y\notin B_{2\rho}(u)$, we have $|x-y|\geqslant|y-u|/2$.
	Since $v_1=0$ in $B_{2\rho}(u)$, we have the pointwise estimates
	\begin{equation}
		\left| |D|^1v_1(x) \right| \aleq \int_{\T\setminus B_{2\rho}(u)} \frac{|v_1(y)|}{|x-y|^2}|\dd y|\aleq  \frac{1}{\rho^{1/2 }} \int_{\T\setminus B_{2\rho}(u)} \frac{|v_1(y)|}{|y-u|^{3/2}}\, |\dd y|. \label{eq:point-one-v1}
	\end{equation}
	Therefore, we have 
	\begin{equation}\label{eq:der-chi-2}
		\begin{aligned} 
			\left| \int_{\T} \scal{ v_1 }{ |D|^1\left(2\chi\chi'(\eta-\gamma_0)\right)}	\right| & = \left| \int_{\T} \scal{ |D|^1 v_1 }{ 2\chi\chi'(\eta-\gamma_0)}	\right| \\[2mm]
			&\aleq \rho^{-1/2} \left\|2\chi\chi'(\eta-\gamma_0) \right\|_{\Ll^1(\T)} \int_{\T\setminus B_{2\rho}(u)} \frac{|v_1(x)|}{|x-u|^{3/2}}\, |\dd x|\\[2mm]
			&\aleq \|(\eta-\gamma_0)'\|_{\Ll^2(B_{ 2\rho}(u))} \int_{\T\setminus B_{2\rho}(u)} \frac{|(\eta-\gamma_0)'(x)|}{|x-u|^{3/2}}\, |\dd x|.
		\end{aligned} 
	\end{equation}
	\end{proof}
	
	We obtain that for any $\rho>0$, there holds 
	\begin{equation*}
		\begin{aligned} 
			\left[ (\eta-\gamma_0)' \right]_{\W^{1/2,2}(B_{\rho}(u))}^2  & \aleq  \|\eta-\gamma_0\|_{\W^{1,\infty}(B_{2 \rho}(u))} \int_{B_{2 \rho}(u)} |h| + \rho^{\sigma}  \\
			&\quad + \rho^{-1}\, \int_{B_{6 \rho}(u)} |(\eta-\gamma_0)'|^2  + \rho \left( \int_{\T\setminus B_{4\rho}(u)} \frac{|(\eta-\gamma_0)'(y)|}{|y-u|^{3/2}}|\dd y|\right)^2.
		\end{aligned} 
	\end{equation*}
	Given $x\in\T$ and $\rho>0$, we denote 
	\begin{equation*}
		\begin{aligned} 
			E(x,\rho) & \coloneqq \left[ (\eta-\gamma_0)' \right]_{\W^{1/2,2}(B_{\rho}(x))}^2  = \int_{B_{\rho}(x)\times B_{\rho}(x)} \frac{ |(\eta-\gamma_0)'(x)-(\eta-\gamma_0)'(y)|^2 }{|x-y|^2}\, |\dd x|\, |\dd y|.
		\end{aligned} 
	\end{equation*}
	Since $(\eta-\gamma_0)'=0$ on half of the ball $B_{\rho}(u)$, we have 
	\begin{equation*}
		\begin{aligned} 
			\rho^{-1}\, \int_{B_{6 \rho}(u)} |(\eta-\gamma_0)'|^2 
			& \leqslant\frac{1}{2\rho^2} \int_{(B_{8\rho}(u)\setminus B_{7 \rho}(u))\times B_{6\rho}(u)} \left|(\eta-\gamma_0)'(x)-(\eta-\gamma_0)'(y) \right|^2\, |\dd x|\, |\dd y|\\[2mm]
			& \aleq  \int_{(B_{8\rho}(u)\setminus B_{7\rho}(u))\times B_{6\rho}(u)} \frac{ \left|(\eta-\gamma_0)'(x)-(\eta-\gamma_0)'(y) \right|^2 }{ |x-y|^2 }\, |\dd x|\, |\dd y|\\[2mm]
			& \aleq E(u,8\rho) - E(u,6\rho).
		\end{aligned} 
	\end{equation*}
	Similarly, we have 
	\begin{equation*}
		\begin{aligned} 
			\rho & \left( \int_{\T\setminus B_{4\rho}(u)} \frac{|(\eta-\gamma_0)'(y)|}{|y-u|^{3/2}}|\dd y|\right)^2 \\[2mm]
			& \aleq \rho  \left( \sum_{j\geqslant2} \int_{B_{2^{j+1}\rho}(u)\setminus B_{2^j\rho}(u)} \frac{|(\eta-\gamma_0)'(y)|}{|y-u|^{3/2}}|\dd y|\right)^2 \\[2mm]
			& \aleq \rho^{-1}  \left( \sum_{j\geqslant2} 2^{-3j/2} \int_{B_{2^{j+1}\rho}(u)\setminus B_{2^j\rho}(u)} |(\eta-\gamma_0)'| \right)^2 \\[2mm]
			& \aleq \rho^{-1}  \left( \sum_{j\geqslant0} 2^{-3j/2} (2^j\rho)^{1/2}  \|(\eta-\gamma_0)'(y)\|_{\Ll^2\left( B_{2^{j+1}\rho}(u)\setminus B_{2^j\rho}(u) \right)} \right)^2 \\[2mm]
			& \aleq   \left( \sum_{j\geqslant0} 2^{-j}  \left( E(u,2^{j+1} \rho) - E(u,2^j\rho)\right)^{1/2} \right)^2 \\[2mm]
			  & \aleq \sum_{j\geqslant0} 2^{-j}  \left( E(u,2^{j+1} \rho) - E(u,2^j\rho)\right).
		\end{aligned} 
	\end{equation*}	
	Hence, we obtain, for any $\eps>0$, the existence of $\rho_0>0$ such that for any $\rho_0>\rho>0$ it holds
	\begin{equation*}
    \begin{aligned} 
		E(u,\rho) & \leqslant\eps\, E(u,2 \rho)+ C\rho^{\sigma} + C \big(E(u,8\rho) - E(u,6 \rho)\big) +C\sum_{j\geqslant0} 2^{-j}  \left( E(u,2^{j+1} \rho) - E(u,2^j\rho)\right) \\[2mm]
        & \leqslant \eps\, E(u,2 \rho)+ C\rho^{\sigma}  +C\sum_{j\geqslant0} 2^{-j}  \left( E(u,2^{j+1} \rho) - E(u,2^j\rho)\right).
    \end{aligned} 
	\end{equation*}
	For $\eps>0$ small enough, we can apply \cite[Lemma 8]{BRS2016} and we obtain  $\alpha>0$ verifying 
	\begin{equation*}
		E(u,\rho)\aleq \rho^{\alpha}.
	\end{equation*}
    We now prove that this rate of decay also holds nearby $0$. If $x\in B_1(u)\setminus (u,v)$ and $1/10>r>0$, there are two possibilities. If $x>\sqrt{r}$, then we have from \eqref{eq:eps-reg} that $E(x,r)\aleq \frac{r}{|x|-r}\aleq \sqrt{r}$. If $x<\sqrt{r}$, then we have $E(x,r)\leqslant E(0,r+x)\aleq (r+x)^{\alpha}\aleq r^{\alpha/2}$. Using Sobolev embeddings in \cite[Equation (8.4)]{DiNPV2012} together with \cite[Theorem 3.1]{HL2011}, we obtain $(\eta-\gamma_0)'\in \Cd^{0,\tilde{\alpha}}$ for some $\tilde{\alpha}>0$.
    \end{proof}

    As a corollary, we can improve the Sobolev regularity of $\eta$.

    \begin{lemma}\label{lm:Sobolev-reg}
        Let $\eta$ be as in Theorem \ref{th:regu_P2}. Then there is $1>t>1/2$ such that $\eta$ is globally $\W^{t,2}$ in arc-length parameterization.
    \end{lemma} 
    \begin{proof} 
    We denote $k\coloneqq \eta-\gamma_0$. Since $\gamma_0\in \Cd^{1,\beta}(\T)$ for some $\beta>1/2$, there exists $1>\sigma>0$ such that $\gamma_0'\in \W^{(1+\sigma)/2,2}(\T)$. Thanks to Lemma \ref{lm:EL-Mob}, $k$ verifies, in $\T\setminus [u,v]$,
    \begin{equation*}
        \left( |D|^1 k' \right)' = f_1' + f_2' - \left( |D|^1 \gamma_0' \right)'.
    \end{equation*}
    We define $v\coloneqq \Hh(k') \in \Cd^{0,\alpha}(\T)$, where $0<\alpha<\beta$ is given by  Lemma \ref{lm:C1alpha_reg}. Then,  there holds, in $\T\setminus [u,v]$
    \begin{equation*}
        \left( v + |D|^1\gamma_0\right)' = |D|^1(k') +|D|^1(\gamma_0') \in \Ll^1(\T\setminus [u,v]).
    \end{equation*}
    Therefore, we have $v + |D|^1\gamma_0\in \W^{1,1}(\T\setminus [u,v])$. Hence, we obtain $v\in \W^{(1+\sigma)/2,2}(\T\setminus [u,v])$. Given $x\in [u,v]$, there holds
    \begin{equation*}
        |v'(x)| = \left| |D|^1 k'(x) \right| \aleq \int_{\T\setminus [u,v]} \frac{ | k'(y)-k'(x)| }{|x-y|^2}\, |\dd y| \aleq \int_{\T\setminus [u,v]} \frac{ |\dd y|  }{|x-y|^{2-\alpha}} \aleq |x-y|^{1-\alpha}.
    \end{equation*}
    Hence, we have $v\in \W^{1,p}([u,v])$ for $1<p<1/(1-\alpha)$. Since $\W^{1,p}([u,v])\hookrightarrow \W^{3/2-1/p,2}([u,v])$, we obtain that $v\in \W^{1/2+\delta,2}(\T)$ for $\delta = \min\left\{\sigma/2, 1-1/p \right\}$. Hence, there holds $k'\in\W^{1/2+\delta,2}(\T)$.
     \end{proof}

	We now conclude to the global regularity in the general case.
	\begin{lemma}\label{lm:Chalf_minus}
		Let $\eta$ be as in Theorem \ref{th:regu_P2}. 
		Then $\eta$ is globally of class $\Cd^{1,1/2}$.
	\end{lemma}
	\begin{proof} 		
        Thanks to Lemma \ref{lm:Sobolev-reg}, we can initiate the bootstrap argument in  Lemma 4.5 of \cite{BRS2016}. We choose $1/p = 1 - 2\left(t-1/2\right)+\eps$ for $\eps>0$ small enough to apply Corollary 4.4 in \cite{BRS2016} (as in the proof of Lemma 4.5 in \cite{BRS2016}): for any $\vp\in \Cd^{\infty}_c(\T\setminus [u,v])$, it holds
		\begin{equation*}
			|T_1(\eta,\vp)| + |T_2(\eta,\vp)| \aleq  \| \vp\|_{\W^{1,p'}(\T)} .
		\end{equation*}
		We denote $k\coloneqq \eta-\gamma_0$. Since $\gamma_0\in \Cd^{1,\beta}(\T)$ for some $\beta>1/2$, there exists $1>\sigma>0$ such that $\gamma_0'\in \W^{(1+\sigma)/2,2}(\T)$. We deduce that, for any $\vp\in \Cd^{\infty}_c(\T\setminus [u,v])$,
		\begin{equation*}
        \begin{aligned} 
			\left| \int_{\T} \scal{|D|^{1/2} k'}{|D|^{1/2} \vp'} \right| & \aleq  \| \vp\|_{\W^{1,p'}(\T)} + \left\| |D|^{(1+\sigma)/2}\gamma_0' \right\|_{\Ll^2(\T)}\, \left\| |D|^{(1-\sigma)/2}  \vp \right\|_{\Ll^2(\T)} \\[2mm]
            & \aleq  \| \vp\|_{\W^{1,p'}(\T)} +  \| \vp \|_{\W^{(1-\sigma)/2,2}(\T)}.
        \end{aligned} 
		\end{equation*}
		In other words, it holds
		\begin{equation*}
			\left((-\lap)^{1/2} k'\right)' \in \left(\W^{1,p'}_0\cap \W_0^{(1-\sigma)/2,2}(\T\setminus [u,v])\right)^* = \left(\W^{-1,p} + \W^{(\sigma-1)/2,2}\right)(\T\setminus [u,v]).
		\end{equation*}
		Thus, $k$ verifies the following system:
		\begin{equation}\label{eq:syst_k}
			\begin{cases}
				(-\lap)^{1/2} k' \in \left(\Ll^{p}+\W^{\sigma/2,2}\right)(\T\setminus [u,v]), \\[2mm]
				k'=0 \quad \text{ in }[u,v]. 
			\end{cases}
		\end{equation}
		By \Cref{th:global_reg_Lpdata}, we obtain $k'\in \W^{s,r}_0(\R)$ for any $1/2<s<1$ and any $r>1$ such that
        \begin{equation*} 
            \frac{1}{r}-s > \frac{1}{p}-1 = 2\left(\frac{1}{2}-t\right)+\eps.
        \end{equation*}
        We can then bootstrap this argument, similarly to Lemma 4.5 in \cite{BRS2016}, until we reach the integrability exponent $(-\lap)^{1/2} k' \in \Ll^q$ for some $q>2$. We deduce $k'\in \Cd^{0,1/2}(\T)$ by Lemma \ref{lm:C0half}.
	\end{proof}

	\subsection{Regularity of solutions to $\Pc^1_{\Mob}$}\label{sec:reg-P1mob}
	The goal of this section is to prove Theorem \ref{th:reg-P1Mob} below. That is to say, solutions to $\Pc^1_{\Mob}$ are log-Lipschitz. To this end, we first compute the associated Euler--Lagrange equation.
	\begin{lemma}\label{lm:EL}
		Let $Z=(z_1,\ldots,z_n)\in \C^n$ with $z_i\neq z_j$ for $i\neq j$. Let $\tau$ be a relative isotopy class of curves passing through $z_1,\ldots,z_n$. 
        
        A solution $\gamma\colon I\to \Chat$, with $I=\T_L$ for some $L>0$ or $I=\R$, to $\Pc^1_{\Mob}(Z,\tau)$, defined in \eqref{eq:P1Mob}, verifies the following system. We denote $a_1,\ldots,a_n\in I$ the points such that $\gamma(a_i)=z_i$ and let $\nu_i\coloneqq \nu(a_i)$ with $\nu$ the normal vector of $\gamma$.
		There are coefficients $\alpha_1,\ldots,\alpha_n\in\R$ such that 
		\begin{equation}\label{eq:Dirac_mass}
			\forall \vp\in \Cd^{\infty}_c(I;\C),\qquad \delta\Mob(\gamma)[\vp] = \sum_{i=1}^n \alpha_i\, \scal{\vp(a_i)}{\nu_i}.
		\end{equation}
		In other words, it holds
		\begin{equation*}
			\delta\Mob(\gamma) = \sum_{i=1}^n \alpha_i\, \delta_{a_i}\, \nu_i.
		\end{equation*}
	\end{lemma}
	
	\begin{proof}
		We first show that the quantity $\delta\Mob(\gamma)[\vp]$ depends only on $\scal{\vp}{\nu}$. 
		We denote $\tau = -\ii \nu$ the tangent vector of $\gamma$. Given $\vp\in \Cd^{\infty}_c(I;\C)$, we decompose $\vp = \vp^{\top}\tau + \vp^{\perp}\nu$ for some $\vp^{\top},\vp^{\perp}\in\Cd^{\infty}_c(I;\R)$. We consider $(\alpha_t)_{t\in(-1,1)}\subset \Diff(I)$ a family of reparametrization of $I$ such that $\dr_t\alpha_t|_{t=0} = -\vp^{\top}$ and $\alpha_0=\id_I$. Then we have 
		\begin{equation*}
			\begin{aligned} 
			\delta\Mob(\gamma)[\vp] & = \left. \frac{\dd}{\dd t} \right|_{t=0} \Mob(\gamma+t\vp) \\[2mm]
			& = \left. \frac{\dd}{\dd t} \right|_{t=0} \Mob(\gamma\circ\alpha_t+t\vp\circ\alpha_t) \\[2mm]
			& = \delta\Mob(\gamma)[ \vp^{\top} \tau +\vp ] \\[2mm]
			& = \delta\Mob(\gamma)[\vp^{\perp}\nu].
			\end{aligned} 
		\end{equation*}
		Hence the admissible variations for $\Pc^1_{\Mob}(Z,\tau)$ are in one-to-one correspondence with elements of 
		\begin{equation*}
			V\coloneqq \left\{ \vp\in \Cd^{\infty}_c(I;\C) : \scal{\vp(a_1)}{\nu_1}=\cdots = \scal{\vp(a_n)}{\nu_n}=0 \right\}.
		\end{equation*}
		Therefore, given two maps $\vp,\psi\in \Cd^{\infty}_c(I;\C)$, if $\scal{\vp(a_1)}{\nu_1}=\scal{\psi(a_1)}{\nu_1}$,\ldots, $\scal{\vp(a_n)}{\nu_n}= \scal{\psi(a_n)}{\nu_n}$, then it holds
		\begin{equation*}
			\delta\Mob(\gamma)[\vp-\psi]=0.
		\end{equation*}
		This implies that $\delta\Mob(\gamma)[\vp]$ is a linear map in the values $\vp(a_1),\ldots,\vp(a_n)$. This implies \eqref{eq:Dirac_mass} by the representation of linear forms in $\R^n$.
	\end{proof}

	\begin{theorem}\label{th:reg-P1Mob}
    Let $Z=(z_1,\ldots,z_n)\in \C^n$ with $z_i\neq z_j$ for $i\neq j$. Let $\tau$ be a relative isotopy class of curves passing through $z_1,\ldots,z_n$.
		Let $\gamma\colon I\to \C$ be a solution to $\Pc^1_{\Mob}(Z,\tau)$ parameterized by arc-length. Then it holds $\gamma''\in\BMO(I)$. In particular, the tangent vector $\gamma'$ is log-Lipschitz.
	\end{theorem}
	
	\begin{remark}\label{rk:log-singularities}
		The proof actually shows that the singularities of $\gamma''$ (and thus, of the curvature) at the points $a_i$ is of the form $\log|\cdot-a_i|$.
	\end{remark}
	
	\begin{proof}
        Up to a Möbius transformation, we can assume that $I=\T$.
        In dimension 1, the Dirac masses are obtained as derivatives of functions with jumps: we have $\delta_{a_i} \in \W^{-1,\infty}(\T)$ since $\W^{1,1}(\T)\hookrightarrow \Cd^0(\T)$. Thus, the Euler--Lagrange equation as written in Lemma \ref{lm:EL-Mob} now becomes the following: for any $\vp\in \Cd^{\infty}(\T;\C)$,
		\begin{equation*}
			\int_{\T} \scal{|D|^{1/2} \gamma'}{ |D|^{1/2} \vp'} = T_1(\gamma,\vp)-T_2(\gamma,\vp) - \sum_{i=1}^n \alpha_i\, \scal{ \scal{\vp}{\delta_{a_i}}_{\Dr',\Dr}}{ \nu_i}_{\C}.
		\end{equation*}
        Following the proof of Theorem III in \cite{BRS2016} step by step only requires to check Lemma 3.3 in \cite{BRS2016}. It is straightforward to check that the proof still holds. Hence, there exists $4>q>2$ such that $\gamma\in \W^{3/2,q}(\T;\C)$. We can then apply Corollary 4.4 in \cite{BRS2016} with $\beta_0=\beta=0$ and $1/p=2/q+\eps\in(1/2,1)$ for $\eps>0$ small enough. We obtain 
		\begin{equation*}
			\left| \int_{\T} \scal{|D|^{1/2} \gamma'}{ |D|^{1/2} \vp'} \right| \aleq \|\vp\|_{\W^{1,p'}(\T)} + \|\vp'\|_{\Ll^1(\T)} \aleq \|\vp\|_{\W^{1,p'}(\T)} \aleq \|\vp\|_{\W^{3/2-1/p',2}(\T)}.
		\end{equation*}
        We deduce that 
        \begin{equation*} 
            \left( (-\lap)^{1/2}\gamma'\right)'\in \left(\W^{3/2-1/p',2}(\T)\right)^* = \W^{1/p'-3/2,2}(\T). 
        \end{equation*}
        We obtain $\gamma'\in \W^{1/2+1/p',2}(\T)$. We then follow the bootstrap argument in Lemma 4.5 in \cite{BRS2016}: there holds for  $\frac{1}{r}=1-2/p'+\eps$ for $\eps>0$ small enough and any $\vp\in \Cd^{\infty}(\T)$,
        \begin{equation*}
			\left| \int_{\T} \scal{|D|^{1/2} \gamma'}{ |D|^{1/2} \vp'} \right| \aleq \|\vp\|_{\W^{1,r'}(\T)} + \|\vp'\|_{\Ll^1(\T)} \aleq \|\vp\|_{\W^{3/2-1/r',2}(\T)}.
		\end{equation*}
        We obtain $\gamma'\in \W^{1/2+1/r',2}(\T)=\W^{1/2+2/p'-\eps,2}(\T)$. The bootstrap argument in Lemma 4.5 of \cite{BRS2016} leads to $\gamma'\in \W^{s,2}(\T)$ for some $s>1$ and the bootstrap stops at the following step
        \begin{equation*} 
            \left( (-\lap)^{1/2}\gamma'\right)'\in \left(\W^{1,1}(\T)\right)^* = \W^{-1,\infty}(\T). 
        \end{equation*}
        Thus, by \eqref{eq:Hilbert-Laplace}, we have
		\begin{equation*}
			\Hh[\gamma''] = (-\lap)^{1/2}\gamma' \in \Ll^{\infty}(\T).
		\end{equation*}
		Hence, we have $\gamma'' \in \Hh^{-1}\left( \Ll^{\infty}(I) \right)\subset \BMO(I)$. We conclude the proof of Theorem \ref{th:reg-P1Mob} thanks to Lemma \ref{lm:logLip} below. 
	\end{proof}
	
	\begin{lemma}\label{lm:logLip}
		A differentiable function $u\colon [0,1]\to \R$ verifies for any $x,y\in[0,1]$
		\begin{equation*}
			|u(x)-u(y)| \aleq \left( \|u'\|_{\BMO([0,1])} + \|u'\|_{\Ll^1([0,1])} \right)\, |x-y|\, \Big(1 +\log\left(\frac{1}{|x-y|}\Big)\right).
		\end{equation*}
	\end{lemma}
	
	\begin{proof}
        We denote 
        \begin{equation*}
            (u')_{B_r(x)} = \fint_{B_r(x)} u' = \frac{1}{|B_r(x)|}\int_{B_r(x)} u'.
        \end{equation*}
		For $0<r<1/2$ and $B_r(x)\subset [0,1]$, we have 
		\begin{equation*}
			\left| (u')_{B_r(x)} - (u')_{B_{r/2}(x)} \right| \leqslant \fint_{B_{r/2}(x)} \left| (u')_{B_r(x)} - u' \right| \leqslant 2\, \|u'\|_{\BMO(B_r(x))}.
		\end{equation*}
		If $r\approx 2^k$, then we obtain 
		\begin{equation*}
			\left| (u')_{B_{2^k}(x)} \right| \aleq k\, \|u'\|_{\BMO([0,1])} + \|u'\|_{\Ll^1([0,1])}.
		\end{equation*}
		Hence, we have
		\begin{equation*}
			\left| (u')_{B_r(x)} \right| \aleq \log\left(\frac{1}{r}\right)\, \|u'\|_{\BMO([0,1])} + \|u'\|_{\Ll^1([0,1])}.
		\end{equation*}
		Therefore, we obtain 
		\begin{equation*}
			\begin{aligned} 
				\left| u(x) - u(y) \right| & = \left|\int_y^x u'(z)\,\dd z\right| \\[2mm]
				& = |x-y|\, \left|(u')_{B_{|x-y|/2}(\frac{x+y}{2})}\right|  \\[2mm]
				& \aleq |x-y|\, \left[1 + \log\left(\frac{1}{|x-y|}\right)\right] \, \left( \|u'\|_{\BMO([0,1])} + \|u'\|_{\Ll^1([0,1])} \right).
			\end{aligned} 
		\end{equation*}
	\end{proof}
	
	\section{\texorpdfstring{Uniqueness in relative isotopy classes of solutions to $\Pc^1_{\Mob}$}{Uniqueness in relative isotopy classes of solutions to P1Mob}}\label{sec:uniqueness}
	
	The goal of this section is to prove the following result.
	\begin{proposition}
		Let $\tau$ be a relative isotopy class of Weil--Petersson curves passing through distinct points $z_1,\ldots,z_n\in\C$ with $n\geqslant 3$. Assume that there are $\gamma_{-1},\gamma_1\in\tau$ that are both minimizing, that can be chosen with compact image up to inversion. Then there holds $\gamma_{-1}(\T)=\gamma_1(\T)$.
	\end{proposition}
	
	\begin{proof} 
	Given $i\in\{-1,1\}$, we consider $(\gamma_{i,k})_{k\in\N}$ a minimizing sequence of smooth Jordan curves passing through $z_1,\ldots,z_n$ and converging to $\gamma_i$ in the $\Cd^{1,\alpha}$ topology for $1/2<\alpha<1$.
	
	Let $k\in\N$. Let $c_k\in\Cd^{\infty}([-1,1]\times\T;\C)$ be a regular isotopy in $\tau$ such that $c_k(-1,\cdot)=\gamma_{-1,k}(\cdot)$ and $c_k(1,\cdot)=\gamma_{1,k}(\cdot)$. We can choose the sequence $(c_k)_{k\in\N}$ to converge to a $\Cd^{1,\alpha}$ isotopy between $\gamma_{-1}$ and $\gamma_1$ as $k\to \infty$. Therefore, the following number is finite
	\begin{equation}\label{eq:sup-Mob}
		M\coloneqq \sup_{k\in\N} \max_{-1\leqslant t\leqslant 1} \Mob(c_k(t))<\infty.
	\end{equation}
	We denote 
    \begin{equation}\label{eq:value-inf}
        m\coloneqq \Mob(\gamma_1)=\Mob(\gamma_{-1}) = \inf\left\{ \Mob(\gamma) : \gamma\in\tau \right\} . 
    \end{equation} 
    Let $\eps>0$. Up to extracting a subsequence, we can assume that
	\begin{equation}\label{eq:almost-minimizing}
		\sup_{k\in\N} \Mob(\gamma_{i,k}) \leqslant m+\eps.
	\end{equation}
	
	Let $0<\theta<1/100$ and consider $r_{\theta}\in \Cd^{\infty}(\s^1;[-1,1])$ such that 
	\begin{enumerate}
		\item $r_{\theta}(\ee^{\ii\vp})=1$ for $\vp\in[0,\frac{\pi}{2}-\theta]$,
		\item $r_{\theta}(\ee^{\ii\vp})=-1$ for $\vp\in[\frac{\pi}{2}+\theta,\pi]$,
		\item $r_{\theta}(\ee^{-\ii\vp})=r_{\theta}(\ee^{\ii\vp})$ for $\vp\in[0,\pi]$,
		\item $r_{\theta}(\s^1)=[-1,1]$.
	\end{enumerate}
	We consider $a_{j,k}\in\T$ such that 
    \begin{equation}\label{eq:def-ajk}
        \gamma_{1,k}(a_{j,k})=z_j.
    \end{equation}
    We will drop the subscript $\theta$ in $r_\theta$ for a moment, and will only take the limit $\theta \to 0$ at the end.
	\begin{claim}\label{cl:holomorphic-motion}
		There is $g_k\in \Cd^{\infty}(\overline{\D}\times\T;\C)$ with $g_k(\cdot,s)$ holomorphic for any $s\in\T$, such that
		\begin{equation}\label{eq:hm}
			\begin{cases}
				g_k(1,s)=\gamma_{1,k}(s), \\[2mm]
				\forall j\in\{1,\ldots,n\},\ \forall z\in\D,\ g_k(z,a_{j,k})=z_j, \\[2mm]
				\forall \zeta\in\s^1,\ 
				g_k(\zeta,\T)=c_k(r(\zeta),\T) \text{ as Jordan curves} .
			\end{cases}
		\end{equation}
		Every $\gamma_{z,k}=g_k(z,\cdot)\colon\T\to\C$ for $z\in\overline{\D}$ is a regular embedding.
	\end{claim} 
	\begin{proof} 
		We drop the index $k$ for simplicity.
		
		We first claim that up to a translation by a holomorphic moving point, we can assume that all the curves enclose a fixed neighborhood of the origin.
		
		Let $t\in[0,1]$ and $a_t\in\T$ be such that $c(t,a_t)=z_1$. We denote $n_t\in\s^1$ the normal of $c(t,\cdot)$ at $a_t$ pointing inward the bounded domain $\Omega_t$ enclosed by $c(t,\T)$. The map $s\in[0,1]\mapsto n_{s r(\zeta)}$ defines a homotopy between $n_{r(\zeta)}$ and the constant $n_0$. Hence, the map $n_{r(\cdot)}\colon\s^1\to \s^1$ has degree $0$ and there exists $\Psi\colon \s^1\to \R$ such that $n_{r(\zeta)}=\ee^{\ii \Psi(\zeta)}$.
		Let $\Hh$ be the Hilbert transform on $\s^1$. We denote
		\begin{equation*}
			\forall \zeta\in\s^1,\qquad G(\zeta) \coloneqq  \ee^{-\Hh \Psi(\zeta)+\ii\Psi(\zeta)}= \ee^{-\Hh\Psi(\zeta)}\, n_{r(\zeta)}.
		\end{equation*}
		Then $G$ extends to a smooth holomorphic function on $\overline{\D}$ since the function $-\Hh \Psi(\zeta) + \ii\Psi$ extends to a holomorphic function on $\overline{\D}$.
        Moreover, for $\zeta \in \s^1$, $G(\zeta)$ points inward $\Omega_{r(\zeta)}$ at $a_{r(\zeta)}$.
		Choose $\delta>0$ small enough so that, if we denote the holomorphic function 
        $$b \coloneqq z_1+\delta G \colon \overline{\D} \to \C,$$ then $b(\zeta)\in \Omega_{r(\zeta)}$ holds for all $\zeta\in\s^1$. 
	   Hence, there exists $\rho>0$ such that for all $\zeta\in\s^1$, the closed disk $\overline{\D_{\rho}}$ is contained in the interior of $\Gamma_\zeta \coloneqq c(r(\zeta),\T)-b(\zeta)$.

		We now apply \cite[Corollary 2]{F1988} (see also \cite[Filling Theorem p.7]{BF2016} or \cite[Theorem 4.1]{AM2001}) to $(\Gamma_\zeta)_{\zeta \in \s^1}$, 
        we obtain $\widetilde g\in \Cd^{\infty}(\overline{\D}\times\T;\C)$ such that for any $s\in\T$, the function $\widetilde g(\cdot,s)$ satisfies:
		\begin{equation}\label{eq:cond-g}
			\begin{cases} 
           \widetilde g(\cdot,s) \text{ is holomorphic  in } \D, \\[2mm]
				\widetilde g(\cdot,s)\in \Cd^0\left( \overline{\D};\C \right), \\[2mm]
				\widetilde g(1,s) = \gamma_1(s) - b(1), \\[2mm]
				\forall \zeta\in\s^1,\ \widetilde g(\zeta,s) \in \Gamma_\zeta, \\[2mm]
				\forall (z,s)\in\overline{\D}\times\T,\ \widetilde g(z,s) \neq 0.
			\end{cases}
		\end{equation}
This function $\widetilde g(\cdot, s)$ is unique, by Theorem~3.v of \cite{F1988}, or Theorem 4.1 of \cite{AM2001}.

Moreover, $\widetilde g(z,\cdot)\colon\T\to\C$  is a regular embedding for any $z\in\overline{\D}$, see for instance Equation (4.8) or the proof of Theorem 3 on p.885 in \cite{F1988}. Together with the second last property in \eqref{eq:cond-g}, we obtain that for $\zeta\in\s^1$, the set $\widetilde g(\zeta,\T)$ is a Jordan curve included in $\Gamma_\zeta$. Hence, there holds 
$$
    \widetilde g(\zeta,\T)=\Gamma_\zeta, \qquad \forall \zeta\in\s^1.
$$ 

        Now let $g (z, s) \coloneqq \widetilde{g} (z,s) +b(z)$. We have $g (\zeta, \T) = c (r(\zeta),\T)$.

		We now check the second condition in \eqref{eq:hm}
		\begin{equation*} 
			\forall j\in\{1,\ldots,n\},\ \forall z\in\D,\ g(z,a_{j,k})=z_j.
		\end{equation*} 
For this, by possibly reducing $\delta$, we assume further that
\begin{equation*}
			\delta\, \|G\|_{\Ll^{\infty}(\D)}< \frac{1}{2}\, \min_{j\geqslant 2}|z_j-z_1|.
		\end{equation*}
		Then $z_1-b=-\delta\, G$ and $|z_j-b|\geqslant |z_j-z_1|-\delta\|G\|_{\Ll^{\infty}(\D)}>0$ for $j\geqslant 2$.
        
		Thus, for $j \geqslant 2$, the functions $\widehat g_j \colon z \mapsto z_j-b (z)$ have no zero on $\overline{\D}$.
 We check easily that $\widehat g_j$ also satisfies  all other conditions in \eqref{eq:cond-g} for $s=a_{j,k}$.  Therefore, the uniqueness implies $\widetilde g(z,a_{j,k}) = \widehat{g}(z)$ which means $g(z,a_{j,k}) \equiv z_j$.
	\end{proof} 
	
	Using Claim \ref{cl:holomorphic-motion}, we obtain a family of curves $\gamma_{z,k}\coloneqq g_k(z,\cdot)$ for $z\in\D$ parameterizing Jordan curves of class $\Cd^{\infty}$. For $s,t\in\T$, we denote 
	\begin{equation*}
		q_{\gamma_{z,k}}(s,t) = \frac{\gamma_{z,k}'(s)\, \gamma_{z,k}'(t)}{ (\gamma_{z,k}(s)-\gamma_{z,k}(t))^2 }.
	\end{equation*}
	For any fixed $s\in[0,2\pi)$, the function $z\mapsto \gamma_{z,k}'(s)$ is holomorphic and has no zero in $\overline{\D}$. Hence, it has a square root $\sqrt{\gamma_{z,k}'(s)}$. Moreover, we can choose it so that it is continuous in $s$ and holomorphic in $z$. We obtain a square root of $q_{\gamma_{z,k}}$ that we denote
	\begin{equation*}
		k_{\gamma_{z,k}}(s,t) \coloneqq \frac{\sqrt{\gamma_{z,k}'(s)}\, \sqrt{\gamma_{z,k}'(t)} }{ \gamma_{z,k}(s)-\gamma_{z,k}(t) } .
	\end{equation*}
	We have $k_{\gamma_{z,k}}(s+2\pi,t)=k_{\gamma_{z,k}}(s,t+2\pi) = -k_{\gamma_{z,k}}(s,t)$.
	
	\begin{claim}
		There holds 
		\begin{equation}\label{eq:contour}
			\fint_{\s^1} \Mob(\gamma_{z,k})\, |\dd z| - \Mob(\gamma_{0,k}) = \fint_{\s^1} \int_{[0,2\pi)\times[0,2\pi)} |k_{\gamma_{z,k}}(s,t)-k_{\gamma_{0,k}}(s,t)|^2\, \dd s\, \dd t\, |\dd z|.
		\end{equation}
	\end{claim}
	\begin{proof} 
    To simplify the notation, we drop the index $k$ in this proof.
    
		We denote $\gamma'_z(s)= |\gamma'_z(s)|\, \ee^{\ii\theta_z(s)}$ and $\gamma_z(s)-\gamma_z(t) = |\gamma_z(s)-\gamma_z(t)|\, \ee^{\ii \alpha_z(s,t)}$. Then it holds
		\begin{equation*}
			\begin{aligned} 
			\frac{\gamma'_z(s)\, \gamma'_z(t)}{ (\gamma_z(s)-\gamma_z(t))^2 } & = \frac{|\gamma_z'(s)|\, |\gamma'_z(t)|}{|\gamma_z(s)-\gamma_z(t)|^2}\,  \ee^{\ii\left[\theta_z(s)+\theta_z(t)-2\alpha_z(s,t) \right]} \\[2mm]
			 & = \left( \cos\theta_{\gamma_z}(s,t) + \ii\sin\theta_{\gamma_z}(s,t)\right)\, \frac{|\gamma_z'(s)|\, |\gamma'_z(t)|}{|\gamma_z(s)-\gamma_z(t)|^2}.
			\end{aligned} 
		\end{equation*}
		In the above formula, the term $\theta_{\gamma_z}$ denotes the conformal angle of $\Gamma_z$ defined in Definition \ref{def:conformal-angle} and we used the equality \eqref{eq:theta_tau}.
		Therefore, we have by \eqref{eq:cos-law}, 
		\begin{equation*}
			\begin{aligned} 
				\Mob(\gamma_z)-4 & = \int_{[0,2\pi)\times[0,2\pi)} \left( |q_{\gamma_z}(s,t)| - \Re[q_{\gamma_z}(s,t)]\right)\, \dd s\, \dd t \\[2mm]
				& = \int_{[0,2\pi)\times[0,2\pi)} \left( |k_{\gamma_z}(s,t)|^2 - \Re[k_{\gamma_z}(s,t)^2]\right)\, \dd s\, \dd t.
			\end{aligned} 
		\end{equation*}
		By holomorphicity in $z$, we have for any $s\neq t$
		\begin{equation*}
			q_{\gamma_0}(s,t) = \fint_{\s^1} q_{\gamma_z}(s,t)\, |\dd z|,
            \qquad \text{ and }\qquad 
            k_{\gamma_0}(s,t) = \fint_{\s^1} k_{\gamma_z}(s,t)\, |\dd z|.
		\end{equation*}
		Hence, it holds
		\begin{equation*}
			\begin{aligned} 
				& \fint_{\s^1}\left( |q_{\gamma_z}(s,t)| - \Re[q_{\gamma_z}(s,t)] \right)\, |\dd z| \\[2mm]
                & =  \fint_{\s^1}|k_{\gamma_z}(s,t)|^2 |\dd z|- \Re[q_{\gamma_0}(s,t)] \\[2mm]
				& = \fint_{\s^1} \left( |k_{\gamma_z}(s,t) - k_{\gamma_0}(s,t)|^2 + 2\Re\left( k_{\gamma_z}(s,t)\, \overline{k_{\gamma_0}(s,t)} \right) - |k_{\gamma_0}(s,t)|^2\right)|\dd z| - \Re[q_{\gamma_0}(s,t)] \\[2mm]
				& = \fint_{\s^1}  |k_{\gamma_z}(s,t) - k_{\gamma_0}(s,t)|^2 \, |\dd z| + 2\, |k_{\gamma_0}(s,t)|^2 - |k_{\gamma_0}(s,t)|^2- \Re[q_{\gamma_0}(s,t)] \\[2mm]
				& = \fint_{\s^1}  |k_{\gamma_z}(s,t) - k_{\gamma_0}(s,t)|^2 \, |\dd z| +  |q_{\gamma_0}(s,t)|- \Re[q_{\gamma_0}(s,t)].
			\end{aligned} 
		\end{equation*}
		We integrate over the region $I_{\eps}\coloneqq  \{ (s,t) \in [0,2\pi)\times[0,2\pi) : |s-t|\geqslant \eps \mod 2\pi \}$ and obtain 
        \begin{equation*}
        \begin{aligned} 
            & \fint_{\s^1} \left[ \int_{I_{\eps}} \left( |q_{\gamma_z}(s,t)| - \Re[q_{\gamma_z}(s,t)] \right)\dd s\, \dd t\right]\, |\dd z|  - \int_{I_{\eps}} \left( |q_{\gamma_0}(s,t)|- \Re[q_{\gamma_0}(s,t)] \right) \dd s\, \dd t \\[2mm]
            & = \fint_{\s^1} \left[ \int_{I_{\eps}} |k_{\gamma_z}(s,t) - k_{\gamma_0}(s,t)|^2\, \dd s\, \dd t\right]\, |\dd z|.
        \end{aligned} 
        \end{equation*}
        As $\eps\to 0$, the right-hand side converges by montone convergence theorem. The convergence of the left-hand side follows from dominated convergence theorem together with the uniform $\Cd^2$-regularity of the family $(\gamma_z)_{z\in\overline{\D}}$ (see for instance the proof of Proposition 3.1.1 in \cite{ohara2003}).
        This provides the formula \eqref{eq:contour}.
	\end{proof}

		By \eqref{eq:sup-Mob}, \eqref{eq:value-inf}, \eqref{eq:almost-minimizing} and \eqref{eq:contour}, we have 
		\begin{equation*}
			\begin{aligned} 
				\fint_{\s^1} \int_{[0,2\pi)\times[0,2\pi)} |k_{\gamma_{\zeta,k}}(s,t)-k_{\gamma_{0,k}}(s,t)|^2\, \dd s\, \dd t\, |\dd \zeta|& = \fint_{\s^1}\Mob(\gamma_{z,k})\, |\dd z| - \Mob(\gamma_{0,k})\\[2mm]
				& \leqslant   (1-\theta)(m+\eps) + \theta M - m \\[2mm] 
				& \leqslant \theta(M-m) +\eps,
			\end{aligned} 
		\end{equation*}
        where the first inequality follows from the fact that $\gamma_{0,k}$ is in the isotopy class $\tau$. 
        
		Choosing $\theta>0$ small enough, this implies that there are two points $\zeta_1,\zeta_2\in\s^1$ with $r(\zeta_1)=-1$ and $r(\zeta_2)=1$ such that 
		\begin{equation*}
			\int_{
            [0,2\pi)\times[0,2\pi)} \hspace{-5pt}\Big( |k_{\gamma_{\zeta_1,k}}(s,t)-k_{\gamma_{0,k}}(s,t)|^2+ |k_{\gamma_{\zeta_2,k}}(s,t)-k_{\gamma_{0,k}}(s,t)|^2\Big)\, \dd s\, \dd t\leqslant 3\theta(M-m) + 3\eps.
		\end{equation*}
		By construction, $\gamma_{\zeta_1,k}$ and $\gamma_{\zeta_2,k}$ are reparametrizations of $\gamma_{\pm1,k}$. Hence, there are two diffeomorphisms $\alpha_{-1,k}$ and $\alpha_{1,k}$ of $\T$ such that 
		\begin{equation}\label{eq:slices}
			\begin{aligned} 
				3\theta(M-m)+3\eps & \geqslant \int_{[0,2\pi)\times[0,2\pi)} |k_{\gamma_{-1,k}\circ\alpha_{-1,k}}(s,t)-k_{\gamma_{1,k}\circ\alpha_{1,k}}(s,t)|^2\, \dd s\, \dd t \\[2mm]
				&
                = \int_{[0,2\pi)\times[0,2\pi)} |k_{\gamma_{-1,k}\circ\alpha_{-1,k}\circ\alpha_{1,k}^{-1}}(s,t)-k_{\gamma_{1,k}}(s,t)|^2\, \dd s\, \dd t.
			\end{aligned} 
		\end{equation}
		Moreover, we have by Claim \ref{cl:holomorphic-motion} that $\gamma_{\pm 1,k}\circ\alpha_{\pm 1,k}(a_{j,k})=z_j$. By definition of $a_{j,k}$ in \eqref{eq:def-ajk}, we obtain $\alpha_{1,k}(a_{j,k})=a_{j,k}$ and hence, we have $\gamma_{-1,k}\circ\alpha_{-1,k}\circ\alpha_{1,k}^{-1}(a_{j,k})=z_j$.
		
		Thanks to \eqref{eq:slices}, we have on a domain $R\subset [0,2\pi)\times[0,2\pi)$ that does not intersect the diagonal $\Delta = \{(x,x) : x\in [0,2\pi)\}$,
		\begin{equation*}
			\begin{aligned}
				& \int_{R} \left| q_{\gamma_{1,k}}(s,t)-q_{\gamma_{-1,k}\circ\alpha_{-1,k}\circ\alpha_{1,k}^{-1}}(s,t) \right|\, \dd s\, \dd t \\[2mm]
				& = \int_{R} \bigg|2k_{\gamma_{1,k}}(s,t)\left( k_{\gamma_{-1,k}\circ\alpha_{-1,k}\circ\alpha_{1,k}^{-1}}(s,t) - k_{\gamma_{1,k}}(s,t) \right) \\[2mm]
                &\quad +\left( k_{\gamma_{-1,k}\circ\alpha_{-1,k}\circ\alpha_{1,k}^{-1}}(s,t) - k_{\gamma_{1,k}}(s,t) \right)^2 \bigg|\, \dd s\, \dd t .
			\end{aligned}
		\end{equation*}
        By Cauchy--Schwarz inequality, we obtain a constant $C>0$ depending only on $m,M,\|q_{\gamma_{1,k}}\|_{\Ll^1(R)}$ such that
        \begin{equation}\label{eq:rectangle}
            \int_{R} \left| q_{\gamma_{1,k}}(s,t)-q_{\gamma_{-1,k}\circ\alpha_{-1,k}\circ\alpha_{1,k}^{-1}}(s,t) \right|\, \dd s\, \dd t  \leqslant C\, \left( \sqrt{\theta} + \sqrt{\eps}\right).
        \end{equation}
		We consider $R$ of the form $R = [a_1,a_2]\times [a_3,x]$ (where we dropped the index $k$ in $a_{j,k}$). Then we have, with $\widetilde{\gamma}\coloneqq \gamma_{-1,k}\circ\alpha_{-1,k}\circ\alpha_{1,k}^{-1}$,
		\begin{equation}\label{eq:der-log}
			\dr_s\dr_t \log\left(\widetilde{\gamma}(s)-\widetilde{\gamma}(t) \right) = q_{\widetilde{\gamma}}(s,t).
		\end{equation}
		Since we restrict the above computations to $(s,t)\in R$, the above logarithm is well-defined. Integrating twice \eqref{eq:der-log}, we obtain 
		\begin{equation*}
			\begin{aligned} 
				\int_R q_{\widetilde{\gamma}}(s,t)\, \dd s\, \dd t & = \int_{a_3}^x \dr_t \log\left( \frac{ \widetilde{\gamma}(a_2)-\widetilde{\gamma}(t) }{\widetilde{\gamma}(a_1)-\widetilde{\gamma}(t)}  \right)\, \dd t \\[2mm]
				& = \log\left( \frac{ \widetilde{\gamma}(a_2)-\widetilde{\gamma}(x) }{\widetilde{\gamma}(a_1)-\widetilde{\gamma}(x)}\, \frac{ \widetilde{\gamma}(a_1)-\widetilde{\gamma}(a_3) }{\widetilde{\gamma}(a_2)-\widetilde{\gamma}(a_3)} \right)
			\end{aligned} 
		\end{equation*}
		Therefore, it holds
		\begin{equation*}
			\exp\left( \int_R q_{\widetilde{\gamma}}(s,t) \, \dd s\, \dd t\right) = \frac{ \left( \widetilde{\gamma}(a_2)-\widetilde{\gamma}(x) \right)  \left( \widetilde{\gamma}(a_1)-\widetilde{\gamma}(a_3) \right) }{\left( \widetilde{\gamma}(a_1)-\widetilde{\gamma}(x) \right)\left( \widetilde{\gamma}(a_2)-\widetilde{\gamma}(a_3) \right)  } =  \frac{ \left( z_2-\widetilde{\gamma}(x) \right)  \left( z_1-z_3 \right) }{\left( z_1-\widetilde{\gamma}(x) \right)\left( z_2-z_3 \right)  }.
		\end{equation*}
		Therefore, we obtain from \eqref{eq:rectangle}
		\begin{equation*}
			\begin{aligned}
				\left| \frac{z_2-\widetilde{\gamma}(x)}{z_1-\widetilde{\gamma}(x)} - \frac{z_2-\gamma_{1,k}(x)}{z_1-\gamma_{1,k}(x)} \right| & = \left| \frac{z_2-z_3}{z_1-z_3}  \right|\, \left| \exp\left( \int_R q_{\widetilde{\gamma}} \right) - \exp\left( \int_R q_{\gamma_{1,k}}\right) \right| \\[4mm]
				& \leqslant C(\|q_{\gamma_{1,k}}\|_{\Ll^1(R)})\, \left| \frac{z_2-z_3}{z_1-z_3}  \right|\, \left| \ee^{C \left( m,M,\left\| q_{\gamma_{1,k}} \right\|_{\Ll^1(R)} \right)\, \left( \sqrt{\theta} + \sqrt{\eps} \right) } -1 \right|.
			\end{aligned} 
		\end{equation*}
		For $\theta$ small enough, we obtain for any $x\in [0,2\pi)\setminus \left( [a_1,a_2]\cup\{a_3\}\right)$, that
		\begin{equation*}
			\left| \frac{z_2-\widetilde{\gamma}(x)}{z_1-\widetilde{\gamma}(x)} - \frac{z_2-\gamma_{1,k}(x)}{z_1-\gamma_{1,k}(x)} \right| \leqslant C\, \left( \sqrt{\theta} + \sqrt{\eps}\right),
		\end{equation*}
        where $C$ depends only on $m,M,z_1,z_2,z_3,\|q_{\gamma_{1,k}}\|_{\Ll^1(R)}$.
		On the other hand, it holds
		\begin{equation*}
			\begin{aligned} 
				\left| \frac{z_2-\widetilde{\gamma}(x)}{z_1-\widetilde{\gamma}(x)} - \frac{z_2-\gamma_{1,k}(x)}{z_1-\gamma_{1,k}(x)} \right| & = \left| \frac{ \left(z_2-\widetilde{\gamma}(x)\right)\, \left(z_1-\gamma_{1,k}(x)\right) -\left(z_1-\widetilde{\gamma}(x)\right)\, \left(z_2-\gamma_{1,k}(x) \right) }{ \left(z_1-\widetilde{\gamma}(x)\right)\, \left(z_1-\gamma_{1,k}(x)\right) } \right| \\[4mm]
				& = \left| \frac{ 
                    -z_2\, \gamma_{1,k}(x) - z_1\, \widetilde{\gamma}(x) + z_2\, \widetilde{\gamma}(x)+z_1\, \gamma_{1,k}(x) 
                        }{ 
                    \left(z_1-\widetilde{\gamma}(x)\right)\, \left(z_1-\gamma_{1,k}(x)\right) 
                    } \right| \\[4mm]
				& = |z_2-z_1|\, \left| \frac{\gamma_{1,k}(x)-\widetilde{\gamma}(x)}{\left(z_1-\widetilde{\gamma}(x)\right)\, \left(z_1-\gamma_{1,k}(x)\right)} \right|.
			\end{aligned} 
		\end{equation*}
		Consequently, there is a constant $C>0$ depending only on $m,M,z_1,z_2,z_3,\|q_{\gamma_{1,k}}\|_{\Ll^1(R)}$ such that
		\begin{equation*}
			|\gamma_{1,k}(x)-\widetilde{\gamma}(x)| \leqslant C\, \left| \left(z_1-\widetilde{\gamma}(x)\right)\, \left(z_1-\gamma_{1,k}(x)\right) \right|\,\left( \sqrt{\theta} + \sqrt{\eps} \right).
		\end{equation*}
		There exists $\rho>1$, independent of $k$, such that both $\gamma_{1,k}(\T)$ and $\gamma_{-1,k}(\T)$ are included in $\D_{\rho}$. We obtain a constant $C>0$ depending only on $m,M,z_1,z_2,z_3,\|q_{\gamma_{1,k}}\|_{\Ll^1(R)},\rho$ such that 
		\begin{equation*}
			|\gamma_{1,k}(x)-\widetilde{\gamma}(x)| \leqslant C\, \left( \sqrt{\theta} + \sqrt{\eps} \right).
		\end{equation*}
		Up to exchanging the roles of $z_1,z_2$ and $z_3$, and thus, changing the domain $R$, we obtain that $x$ can cover $\T$. Therefore, we obtain that for all $0<\theta<1$, there holds
		\begin{equation}\label{eq:Linfty-distance}
			\sup_{x\in\T}|\gamma_{1,k}(x)-\widetilde{\gamma}(x)| \leqslant C\, \left( \sqrt{\theta} + \sqrt{\eps}\right),
		\end{equation}
		for some constant $C>0$ depending only on $m,M,z_1,z_2,z_3,\rho$ and $\|q_{\gamma_{1,k}}\|_{\Ll^1\left((\T\times\T)\setminus B_{\delta}(\Delta)\right)}$ for some $\delta>0$ depending only on $\gamma_1,z_1,\ldots,z_n$.
		
		Although the map $\widetilde{\gamma}$ depends on $\theta$, it is a reparameterization of $\gamma_{-1,k}(\T)$. Hence the estimate \eqref{eq:Linfty-distance} implies that the Hausdorff distance between the curves $\gamma_{-1,k}(\T)$ and $\gamma_{1,k}(\T)$ satisfies
		\begin{equation*}
			\dist_{\rm Hausdorff}(\gamma_{1,k}(\T),\gamma_{-1,k}(\T))  \leqslant C\, \left( \sqrt{\theta} + \sqrt{\eps} \right).
		\end{equation*}
		Letting first $k\to \infty$ then $\theta\to 0$ then $\eps\to 0$, we obtain that $\gamma_1(\T)=\gamma_{-1}(\T)$.
		\end{proof} 

    \begin{remark}
        The same proof applies to $\Pc^2_{\Mob}(\gamma_0,u,v)$ in \eqref{eq:P2Mob}. In this case, there is only one relative isotopy class. Hence, there is uniqueness of the solution among all closed extensions of $\gamma_0|_{(u,v)}$.
    \end{remark}
    
	\section{\texorpdfstring{A balancing law for the singularities of $\Pc^1_{\Mob}$}{A balancing law for the singularities of Pc1Mob}}\label{sec:balance}
	
	In this section, we prove that the Möbius invariance of $\Mob$ imposes some relations between the singularities in Lemma \ref{lm:EL}.
	
	\begin{proposition}\label{pr:balancing_law}
        We use the notations of Lemma \ref{lm:EL} and denote $r_i\coloneqq \overline{\alpha_i\, \nu_i} \in\C$. Then it holds
		\begin{equation}\label{eq:balancing-law}
			0 = \sum_{i=1}^n r_i = \sum_{i=1}^n r_i z_i = \sum_{i=1}^n r_i z_i^2.
		\end{equation}
	\end{proposition}

    \begin{remark}
        This is reminiscent of Theorem 4.3 in \cite{MRW2026}.
    \end{remark}
	
	\begin{proof}
		We test \eqref{eq:Dirac_mass} against infinitesimal M\"obius transformations. Since the energy is invariant by Möbius transformation \cite[Theorem 2.1(i)]{FHW1994}, the variation vanishes completely.
		
		\underline{Translations.} For the variation $\gamma_t\coloneqq \gamma +ta$ with $a\in\C$, we have $\dot{\gamma}=a$. Hence
		\begin{equation*}
			\forall a\in\C,\qquad \scal{a}{\sum_{i=1}^n \overline{r_i} }=0.
		\end{equation*}
		Since this is true for every $a\in\C$, we get $\sum_i r_i=0$.
		
		\underline{Dilations and rotations.} For the variation $\gamma_t\coloneqq (1+t)\gamma$, we have $\dot{\gamma}=\gamma$. Therefore
		\begin{equation*}
			0=\sum_{i=1}^n \scal{z_i}{\overline{r_i}} =\Re\left(\sum_{i=1}^n z_i\, r_i\right).
		\end{equation*}
		For the infinitesimal rotation $\dot{\gamma}=\ii\gamma$, we similarly get
		\begin{equation*}
			0=\sum_{i=1}^n \scal{\ii z_i}{\overline{r_i} } = \Re\left(\ii\sum_{i=1}^n z_i\, r_i\right).
		\end{equation*}
		We obtain the second relation in \eqref{eq:balancing-law}.
		
		\underline{Special conformal transformations.} Define
		$\iota(z)=\frac{z}{|z|^2} = \frac{1}{\bar{z}}$ for $z\ne0$. Given $a\in\C$, we consider, for
		$|t|$ sufficiently small, and denoting $T_{-ta}$ the translation by the vectore $-ta$, the family
		\begin{equation*}
			\gamma_t = \iota\circ T_{-ta}\circ\iota \circ\gamma = \frac{1}{\frac{1}{\gamma}-t\bar{a}} = \frac{\gamma }{1-t\, \bar{a}\, \gamma}.
		\end{equation*}
		Differentiating at $t=0$, we obtain
		\begin{equation*}
			\left. \frac{\dd \gamma_t}{\dd t} \right|_{t=0} = \gamma^2 \bar{a}.
		\end{equation*}
		Thus, for every $a\in\C$, we obtain 
		\begin{equation*}
			\Re\left( \sum_{i=1}^n \overline{r_i}\, z_i^2\, \bar{a} \right) =0.
		\end{equation*}
		Taking $a=1$ and then $a=\ii$, we obtain the last identity in \eqref{eq:balancing-law}.
	\end{proof}
	
	As a consequence, we obtain the following result.
	
	\begin{lemma}\label{lm:polynomial}
		Assume $n\geqslant 4$ and that $z_1,\ldots,z_n$ are pairwise distinct. Define
		\begin{equation*}
			P(z) = \prod_{i=1}^{n} (z-z_i).
		\end{equation*}
		There exists a polynomial $Q\in\C[z]$ with $\deg Q\leqslant n-4$ and such that
		\begin{equation}\label{eq:polynomial}
			r_i = \frac{Q(z_i)}{P'(z_i)}.
		\end{equation}
	\end{lemma}
	
	\begin{remark}
		For $n\leqslant 3$, there always exists a circle passing through the points $z_1$, $z_2$ and $z_3$. Hence, we must have $Q=0$.
		For $n=4$, $Q$ is a constant. If $Q=0$, then we get $r_i=0$ and the right-hand side of the Euler--Lagrange equation \eqref{eq:Dirac_mass} vanishes. We obtain a planar critical point of $\Mob$, which has to be a round circle by \cite[Theorem 4.7]{B2020}.
	\end{remark}
	
	\begin{proof}
		We define
		\begin{equation*}
			R(z) = \sum_{i=1}^n	\frac{r_i}{z-z_i} = \sum_{i=1}^n \frac{r_i}{z}\, \frac{1}{1-\frac{z_i}{z}}.
		\end{equation*}
		If $|z|> \max_i |z_i|$, we obtain
		\begin{equation*}
			R(z) = \sum_{i=1}^n \left( \frac{r_i}{z} \sum_{k=0}^{\infty}  \left(\frac{z_i}{z}\right)^k \right) = \sum_{k=0}^{\infty} \frac{1}{z^{k+1}} \left(\sum_{i=1}^n r_i z_i^k \right).
		\end{equation*}
		By Proposition \ref{pr:balancing_law}, the terms corresponding to $k=0,1,2$ vanish. Therefore
		\begin{equation*}
			R(z) = \sum_{k=3}^{\infty} \frac{1}{z^{k+1}} \left(\sum_{i=1}^n r_i z_i^k \right),
		\end{equation*}
		and hence $R(z)= O(|z|^{-4})$ as $|z|\to \infty$.
		
		We now set $Q(z)\coloneqq P(z)R(z)$. 
		The poles of $R$ are simple poles at $z_1,\ldots,z_n$, and they are cancelled by the corresponding simple zeros of $P$. Thus $Q$ extends to a polynomial. Since $P$ has degree $n$ and $R(z)=O(|z|^{-4})$, we deduce that $Q(z)=O(|z|^{n-4})$ as $|z|\to\infty$. Hence, we obtain $\deg Q\leqslant n-4$.
		
		It remains to compute the value of $Q$ at the marked points. Near $z_i$, we have
		\begin{equation*}
			\begin{cases} 
				P(z)=P'(z_i)(z-z_i)+O(|z-z_i|^2), \\[2mm]
				\displaystyle R(z)=\frac{r_i}{z-z_i}+O(1).
			\end{cases} 
		\end{equation*}
		Multiplying the two expansions, we obtain
		\begin{align*}
			Q(z)
			&=\Big(P'(z_i)(z-z_i)+O(|z-z_i|^2)\Big)\left(\frac{r_i}{z-z_i}+O(1)\right)\\
			&=r_iP'(z_i)+o(1).
		\end{align*}
		Letting $z\to z_i$, we obtain $Q(z_i)=r_iP'(z_i)$, which is \eqref{eq:polynomial}.
	\end{proof}

	\section{Round-arc unique continuation}\label{sec:Round_pieces}

    In this section we show that the Loewner energy and the M\"obius energy exhibit similar rigidity. 
    
\subsection{For Loewner energy}\label{sec:unique_cont_IL}
\begin{theorem}\label{thm:unique-continuation-IL}
		Let $\eta$ be a solution to $\Pc^2_{\IL}(\gamma_0,u,v)$. If a nonempty open subarc of the optimized part of $\eta|_{\T\setminus (u,v)}$ is round, then the whole curve $\eta(\T)$ is a circle. In particular, $\gamma_0|_{(u,v)}$ is also a piece of a circle.
	\end{theorem}
    \begin{proof}
        By Lemma~\ref{lem:minI}, $\Gamma_1 := \eta|_{\T\setminus (u,v)}$ is the hyperbolic geodesic in $\Chat \setminus \Gamma_0$ connecting the two end points of $\Gamma_0$, where $\Gamma_0$ is the image of $\gamma_0|_{(u,v)}$. 
        Let $A$ be a nonempty open subarc of $\Gamma_1$ such that $A$ is a circular arc by assumption. As a subarc of a hyperbolic geodesic, $A$ is also a hyperbolic geodesic in the simply connected domain $D : = \Chat \setminus (\Gamma_0 \cup \Gamma_1 \setminus A) = \Omega \cup \Omega^* \cup A$, where $\Omega$ and $\Omega^*$ are the connected components of $\Chat \setminus (\Gamma_0 \cup \Gamma_1)$. 

        Therefore, there exists an anti-holomorphic involution $\iota$ of $D$ where $A$ is the set of fixed points, sending $\Omega$ to $\Omega^*$ and vice versa. (If $D$ is the upper half-plane $\mathbb H$ and the geodesic is $\ii \mathbb R_+$, then $\iota$ is the involution $x + \ii y \mapsto -x + \ii y$. For a general simply connected domain $D$, $\iota$ is obtained from a conjugation by a conformal map $D \to \mathbb H$.) 
        On the other hand, the reflection $R : \Chat \to \Chat$ with respect to the circle $\mathcal C$ containing $A$ coincides with $\iota$ along $A$. 
        By the unique analytic continuation, $\iota=R|_D$.  We claim that it implies that $D = \Chat \setminus (\mathcal C \setminus A)$, i.e., $\Gamma_0 \cup \Gamma_1 = \mathcal C$.

        For this, after a conjugation by a M\"obius map, we assume that $\mathcal C = \widehat {\mathbb R}$, then $\iota$ and $R$ are given by the complex conjugation, and $A$ is an interval of $\widehat {\mathbb R}$.  For all interval $I \subset \widehat {\mathbb R} \setminus A$, if $I \subset \Omega$, then $I = \iota (I) \subset \Omega^*$ which would imply that $\Omega$ and $\Omega^*$ are not disjoint. We obtain that $(\widehat {\mathbb R} \setminus A) \cap \Omega = \emptyset$, similarly  $(\widehat {\mathbb R} \setminus A) \cap \Omega^* = \emptyset$. This concludes the proof of $D = \Chat \setminus (\mathcal C \setminus A)$.
    \end{proof}

	\subsection{For M\"obius energy} \label{sec:unique_mob}
	\begin{theorem}\label{thm:unique-continuation}
		Let $\eta$ be a solution to $\Pc^2_{\Mob}(\gamma_0,u,v)$. If a nonempty open subarc of the optimized part of $\eta|_{\T\setminus (u,v)}$ is round, then the whole curve $\eta(\T)$ is a circle. In particular, $\gamma_0|_{(u,v)}$ is also a piece of a circle.
	\end{theorem}
	
	\begin{proof}
		The optimized part is real-analytic by Theorem \ref{th:regu_P2}. Hence, the whole optimized piece $\eta|_{\T\setminus (u,v)}$ is parameterizing an arc of a generalized circle. Up to applying a Möbius transformation, we can assume that this generalized circle is the axis $\R\times \{0\}$. Parameterizing the curve by arclength $\eta(s) = (x(s),y(s))$ for $s\in\R$ with $\eta(s)=(s,0)$ for $s\leqslant 0$ and for $s\geqslant L>0$, the Euler--Lagrange equation \eqref{eq:EL-Mob} is given by
		\begin{equation}\label{eq:EL-simple}
			\forall p\in \R\setminus (0,L),\qquad \int_{\R}\frac{y(s)}
			{\bigl((x(s)-p)^2+y(s)^2\bigr)^2}\,\dd s=0.
		\end{equation}
		We define the following subset of $\R^6$
		\begin{equation*}
		\Sigma\coloneqq  \overline{ \big\{(x(s),|y(s)|\omega): s\in [0,L],\ y(s)\ne0,\ \omega\in\s^4\big\} } \subset\R\times\R^5.
		\end{equation*}
		It is a compact subset because $y$ is supported in $[0,L]$ and the map $(s,\omega)\mapsto(x(s),|y(s)|\omega)$ is continuous.
		For \(X=(a,Y)\in\R^6\setminus\Sigma\), we denote
		\begin{equation}\label{eq:potential}
			U(a,Y) \coloneqq \frac{1}{|\s^4|}\int_{\R} y(s) \left(\int_{\s^4} \frac{\dd\omega}{ \left((a-x(s))^2+|Y-|y(s)|\omega|^2\right)^2} \right) \dd s.
		\end{equation}
		The function $\Phi(W)\coloneqq |W|^{-4}$ on $\R^6\setminus\{0\}$ is proportional to the Green kernel on $\R^6$ and thus is harmonic. Moreover, it holds
		\begin{equation*}
			U(a,Y) = \frac{1}{|\s^4|}\int_{\s^4\times [0,L]} y(s)\ \Phi\begin{pmatrix}
				a-x(s)\\ Y - |y(s)|\, \omega
			\end{pmatrix}\, \dd \omega\, \dd s.
		\end{equation*}
		Therefore, we get $U\in \Cd^{\infty}(\R^6\setminus \Sigma)$ with $\lap_{\R^6} U=0$ in $\R^6\setminus \Sigma$. Moreover, $U$ is also invariant under
		rotations of $Y$.  For $Y=0$, we get 
		\begin{equation*}
		U(p,0)=\int_{\R}
		\frac{y(s)}{\left( (p-x(s))^2+y(s)^2 \right)^2}\,\dd s.
		\end{equation*}
		Thus \eqref{eq:EL-simple}, implies that $U(p,0)=0$ for $p\in \R\setminus [0,L]$.
		
		\begin{claim}
			Let $\Omega_\infty$ be the unbounded connected component of $\R^6\setminus\Sigma$. It holds
			\begin{equation}\label{eq:U-zero}
				U=0 \qquad\text{on }\Omega_\infty.
			\end{equation}
		\end{claim}
		\begin{proof} 
			Fix $p_0\in I\Subset  \R\setminus [0,L]$ with $I$ being an open interval. Then we have $\dist(\overline I\times\{0\},\Sigma)>0$. 		
			Since $(p_0,0)\notin\Sigma$, the harmonic function $U$ is real-analytic in a Euclidean ball centered at $(p_0,0)$.  After decreasing $\rho>0$, we may arrange simultaneously that
			\begin{enumerate}
				\item we have $(p_0-2\rho,p_0+2\rho)\subset I$,
				
				\item it holds $\Cr_{2\rho}\coloneqq
				\big\{(p,Y):|p-p_0|<2\rho,\ |Y|<2\rho\big\}
				\Subset\R^6\setminus\Sigma$,
				
				\item the Taylor series of $U$ at $(p_0,0)$ converges on a neighbourhood of $\overline{\Cr_{2\rho}}$.
			\end{enumerate}
			
			Let $e_1=(1,0,0,0,0)\in\R^5$, and consider the restriction $v(p,r)\coloneqq U(p,r e_1)$.
			It is real-analytic for $|p-p_0|<2\rho$ and $|r|<2\rho$.  Moreover,
			$v(p,r)=v(p,-r)$ by definition of $U$ in \eqref{eq:potential}.
			Hence, possibly decreasing $\rho$ once more, there are real-analytic functions $a_k\colon (p_0-\rho,p_0+\rho)\to\R$ such that for $|p-p_0|<\rho$ and $|r|<\rho$, we have
			\begin{equation}\label{eq:Taylor-U}
				v(p,r)=\sum_{k=0}^\infty a_k(p)r^{2k}.
			\end{equation}		
			Using the rotational invariance of $U$, we obtain the following expansion
			\begin{equation}\label{eq:Taylor-U2}
				U(p,Y)=\sum_{k=0}^\infty a_k(p)|Y|^{2k}.
			\end{equation}
			We now use the equation $\lap_{\R^6} U=0$. It holds for any integer $j\geqslant 1$,
			\begin{equation*}
			\lap_{\R^5} |Y|^{2j} = \left( \frac{\dr^2}{\dr |Y|^2} + \frac{4}{r}\frac{\dr}{\dr|Y|} \right)|Y|^{2j}=2j(2j+3)|Y|^{2j-2}.
			\end{equation*}
			Hence, we obtain
			\begin{equation*}
			\begin{aligned} 
				0& = \sum_{k=0}^{\infty}  a_k''(p) |Y|^{2k} + \sum_{k=1}^{\infty} a_k(p)\, 2k(2k+3)|Y|^{2k-2}  \\[2mm]
				& = \sum_{k=0}^\infty \left(a_k''(p)+2(k+1)(2k+5)a_{k+1}(p)\right) |Y|^{2k}.
			\end{aligned} 
			\end{equation*}
			We obtain for all $k\geqslant 0$ and $|p-p_0|<\rho$,
			\begin{equation}\label{eq:induction-ak}
				a_k''(p)+2(k+1)(2k+5)a_{k+1}(p)=0.
			\end{equation}
			For $k=0$, we have $a_0(p)=U(p,0)=0$ on $(p_0-\rho,p_0+\rho)$.  Thus $a_0''=0$, and the relation \eqref{eq:induction-ak} with $k=0$ implies that $a_1=0$. Inductively, if $a_k=0$, then it holds $a_k''=0$, and the relation \eqref{eq:induction-ak} implies that $a_{k+1}=0$. Therefore every coefficient $a_k$ vanishes on $(p_0-\rho,p_0+\rho)$, and we obtain $U=0$ on $\{(p,Y):|p-p_0|<\rho,\ |Y|<\rho\}$. We conclude by unique continuation property.
		\end{proof}

		Up to translating the parametrization of $\eta$, we can assume that $\inf\{s\in[0,L] : y(s)\neq 0\}=0$. By the inverse function theorem in $\Cd^1$, there exists $s_1>0$ such that $\eta([0,s_1])$ is a graph 
		\begin{equation}\label{eq:C1_departure_graph}
			\{(t,f(t)) : t\in[0,t_0]\} \qquad \text{ with  }f(0)=f'(0)=0.
		\end{equation}
		Let $\omega_0\in\s^4$ and denote
		\begin{equation*}
			E_{\omega_0}\coloneqq \left\{\omega\in\s^4: |\omega-\omega_0|\leqslant \frac{1}{4} \right\}, \qquad \text{ and }\qquad \beta\coloneqq |E_{\omega_0}|.
		\end{equation*}
		By rotational invariance of $\s^4$, the number $\beta$ is independent of $\omega_0$.  Let $c_0\coloneqq \frac{3\beta}{8 |\s^4|}>0$ and choose $\eps>0$ such that 
		\begin{equation}\label{eq:choice-eps}
			0<\eps\le\frac{1}{4}, \qquad \text{ and }\qquad 2\eps^3\le\frac{c_0}{2}.
		\end{equation} 
		Up to reducing $t_0$, we can assume that 
		\begin{equation}\label{eq:small-slope}
			\|f'\|_{\Ll^{\infty}([0,t_0])}\leqslant \eps .
		\end{equation}
		Let $K\coloneqq \supp(y)\subset [0,L]$. let $\psi\colon [0,x(s_1)]\to [0,s_1]$ be the inverse of $ x\colon [0,s_1]\longrightarrow[0,x(s_1)]$. Let  $J\coloneqq \psi([0,x(s_1)])=[0,s_1]$.
		With the convention $\inf\emptyset=+\infty$, we denote
		\begin{equation}\label{eq:def-d}
			d\coloneqq  \min\left\{ 1,\, \inf_{s\in K\setminus J} |\eta(s)| \right\}.
		\end{equation}
		Then we have $d>0$ since $\eta$ is a chord-arc curve and $\eta(0)=0$. Therefore, for any $s\in K\setminus J$ and any $\omega\in\s^4$, one has
		\begin{equation}\label{eq:norm-d}
			\big|(x(s),|y(s)|\omega)\big| = \sqrt{x(s)^2+y(s)^2} = |\eta(s)| \geqslant d.
		\end{equation}
		Consequently, it holds
		\begin{equation}\label{eq:inclusion-norm-d}
			\overline{\big\{(x(s),|y(s)|\omega): s\in K\setminus J,\ \omega\in\s^4\big\}} 
			\subset \left\{ X\in\R^6:|X|\geqslant d \right\}.
		\end{equation}
		By contradiction, we assume that $f$ is not identically zero on any right neighbourhood of $0$.  Hence, there exists a sequence $t_j\to 0^+$ such that $f(t_j)\neq 0$.

		For $j$ large enough depending on $x(s_1)$ and $\eps$ chosen in \eqref{eq:choice-eps}, we can assume that
		\begin{equation}\label{eq:small-t1}
			t_j <\frac{x(s_1)}{4}, 
            \qquad \text{ and }\qquad 
            \sqrt{1+4\eps^2}\,t_j<\frac{d}{2}.
		\end{equation}
		We define
		\begin{equation}\label{eq:defX}
			X_j\coloneqq \left(t_j,\, \frac{17\, |f(t_j)|}{16}\, \omega_0\right) \in \R^6.
		\end{equation}
		Since $|f(t_j)|\le\eps\, t_j$ by \eqref{eq:small-slope} and \eqref{eq:C1_departure_graph}, we deduce from \eqref{eq:small-t1} that
		\begin{equation}\label{eq:normX0}
			|X_j|< t_1\, \sqrt{1 + \left( \frac{17\, \eps}{16}\right)^2} < \frac{d}{2}.
		\end{equation}
		
		\begin{claim} \label{claim:path}
			We have
			\begin{equation}\label{eq:X0-out}
				X_j \in\Omega_\infty,
			\end{equation}
			and consequently
			\begin{equation}\label{eq:Ux0-0}
				U(X_j)=0.
			\end{equation}
		\end{claim} 

\begin{figure}[ht]
    \begin{center}    \includegraphics[width=\textwidth]{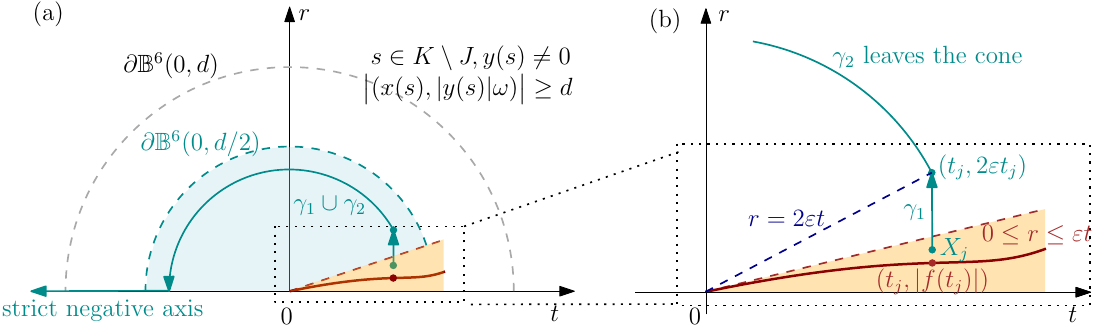}
    \end{center}
    \caption{Illustration of the paths in the proof of Claim~\ref{claim:path}.}
    \label{fig:claim_3_2_path}
\end{figure}

		\begin{proof} 
			From \eqref{eq:small-slope} and \eqref{eq:C1_departure_graph}, we have $|f(t)|\leqslant \eps t$ for all $0\leqslant t\leqslant t_0$.
			We define
			\begin{equation*}
			S_J\coloneqq \{(t,|f(t)|\omega): 0\leqslant t\leqslant x(s_1),\ \omega\in\s^4\}.
			\end{equation*}
			Then $S_J$ is compact and its image under the map $(a,Y)\mapsto(a,|Y|)$ is contained in the cone
			\begin{equation}\label{eq:cone}
				\{(t,r):0\leqslant t\leqslant x(s_1),\ 0\leqslant r\le\eps t\}.
			\end{equation}
			As a consequence of \eqref{eq:inclusion-norm-d}, we have 
			\begin{equation}\label{eq:Sigma-SJ}
				\Sigma\cap \B^6(0,d)\subset S_J.
			\end{equation}
			We denote
			\begin{equation*}
			\rho \coloneqq t_j\sqrt{1+4\eps^2},\qquad \text{ and }\qquad \theta_0 \coloneqq \arctan(2\eps ).
			\end{equation*}
			We join \(X_0\) to the strict negative axis $\{p\in\R^6 : p^1<0,\ p^2=\cdots=p^6=0\}$ by the two paths
			\begin{equation*}
				\begin{cases} 
					\gamma_1(r)=(t_j,r\omega_0), & \frac{17}{16}|f(t_j)|\leqslant r\leqslant 2 \eps t_j,\\[2mm]
					\gamma_2(\theta)=(\rho\cos(\theta), \rho\sin(\theta)\,\omega_0), &\theta_0\le\theta\le\pi.
				\end{cases} 
			\end{equation*}
			They meet because $\gamma_1(2\eps t_j)=\gamma_2(\theta_0)$, and $\gamma_2(\pi)=(-\rho,0)$.  The path $\gamma_1$ does not meet $S_J$ since the radial coordinate verifies $r>|f(t_j)|$.
			The path $\gamma_2$ does not meet $S_J$ either since we have for $\rho\, \cos(\theta)>0$ and $\theta\in [\theta_0,\pi]$, that
			\begin{equation}
				\frac{\rho\, \sin(\theta)}{\rho\, \cos(\theta)} = \tan(\theta) \geqslant \tan(\theta_0) = 2\eps>\eps.
			\end{equation}
			Therefore, $\gamma_2$ never meets the cone \eqref{eq:cone}.
			
			Moreover, we have $|\gamma_1(r)|\leq\rho$ and $|\gamma_2(\theta)|=\rho$.  By \eqref{eq:small-t1} and choice of $\rho$, both paths lie in $\B^6(0,\frac{d}{2})$.
			Hence neither path meets $\Sigma$, meaning that $X_j\in\Omega_{\infty}$.
			We conclude by \eqref{eq:U-zero}.
		\end{proof} 
		
		\begin{claim} 
			It holds
			\begin{equation}\label{eq:lower-bound}
				\sgn\left[f(t_j)\right]\,  U(X_j)\geqslant \frac{ c_0-2\eps^3 }{|f(t_j)|^2}-C_0, \qquad \text{ with }\quad C_0\coloneqq \left(\frac{2}{d}\right)^4 \int_{\R}|y(s)|\,\dd s.
			\end{equation}
		\end{claim} 
		
		\begin{proof} 
			We split the integral \eqref{eq:potential} as follows
			\begin{equation*}
				\begin{aligned}
					U(X_j) & = \frac{1}{|\s^4|} \int_{J\times \s^4} \frac{ y(s) }{ \left((t_j - x(s))^2 + \left| |y(s)|\, \omega - \frac{17}{16}|f(t_j)|\, \omega_0 \right|^2\right)^2 }\, \dd s\, \dd\omega \\[2mm]
					& \qquad + \frac{1}{|\s^4|} \int_{(K\setminus J)\times \s^4} \frac{ y(s) }{ \left((t_j - x(s))^2 + \left| |y(s)|\, \omega - \frac{17}{16}|f(t_j)|\, \omega_0 \right|^2\right)^2 }\, \dd s\, \dd\omega \\[2mm]
					& \eqqcolon U_J(X_j) + U_{K\setminus J}(X_j).
				\end{aligned}
			\end{equation*}
			Thanks to \eqref{eq:inclusion-norm-d} and \eqref{eq:normX0},
			we estimate the part on $K\setminus J$ as follows 
			\begin{equation}\label{eq:U-t2}
				|U_{K\setminus J}(X_j)| \leqslant \left(\frac2d\right)^4\int_{\R}|y(s)|\,\dd s \eqqcolon C_0.
			\end{equation}
			Concerning $U_J$, we change of coordinate $\dd s=\sqrt{1+f'(t)^2}\,\dd t$ and we obtain 
			\begin{equation}\label{eq:UJ}
				U_J(X_j) =\frac{1}{|\s^4|}\int_0^{x(s_1)} f(t)\sqrt{1+f'(t)^2} \left(\int_{\s^4} \frac{ \dd\omega}{|X_j-(t,|f(t)|\omega)|^4} \right)\dd t.
			\end{equation}
			We denote
			\begin{equation*}
				\begin{cases} 
					A_+ \coloneqq \{t\in[0,x(s_1)]:\sgn[f(t_j)] f(t)>0\}, \\[2mm]
					A_-\coloneqq \{t\in[0,x(s_1)]:\sgn[f(t_j)] f(t)<0\}.
				\end{cases} 
			\end{equation*}
			We split the first integral in \eqref{eq:UJ} in $A_+$ and $A_-$ and denote the corresponding two integrals by $U_+(X_j)$ and $U_-(X_j)$.  
			
			For $|t-t_j|\leqslant |f(t_j)|/4$, the estimate \eqref{eq:small-slope} implies that
			\begin{equation}\label{eq:IT}
				\big|\, |f(t)|-|f(t_j)|\, \big|\leqslant |f(t) - f(t_j)| = \big| \sgn[f(t_j)]\, f(t) - |f(t_j)|\, \big| \leqslant \frac{\eps}{4}\, |f(t_j)|.
			\end{equation}
			Hence, we have 
			\begin{equation}\label{eq:onAplus}
				\sgn[f(t_j)]\, f(t)\geqslant |f(t_j)|-\eps|t-t_j|\geqslant \frac{3}{4}\, |f(t_j)|.
			\end{equation}
			Because $|f(t_j)|\leq\eps t_j$, $\eps\leqslant \frac{1}{4}$, and $t_j<x(s_1)/4$, one has
			\begin{equation*}
				0<t_j-\frac{|f(t_j)|}{4} < t_j+\frac{|f(t_j)|}{4}< x(s_1).
			\end{equation*}
			Hence the entire interval $\left[t_j-\frac{|f(t_j)|}{4} , t_j+\frac{|f(t_j)|}{4}  \right]$ lies in $A_+$.  If $\omega\in E_{\omega_0}$, then we have by \eqref{eq:IT}
			\begin{equation*}
				\begin{aligned} 
				\left|\frac{ 17\, |f(t_j)|}{16}\, \omega_0-|f(t)|\, \omega\right|
				&\leqslant \frac{17\, |f(t_j)|}{16}\, |\omega_0-\omega| +\left|\frac{17\, |f(t_j)|}{16}\,-|f(t)|\, \right|\\[2mm]
				&\leqslant \frac{17\, |f(t_j)|}{64}+\frac {|f(t_j)|}{16}+\frac{\eps |f(t_j)|}{4} \\[2mm]
				& \leqslant \frac{25\, |f(t_j)|}{64}  < \frac{|f(t_j)|}{2}.
				\end{aligned} 
			\end{equation*}
			Since we have by assumption that $|t-t_j|\leqslant |f(t_j)|/4$, we obtain 
			\begin{equation}\label{eq:ineq_Xj-Aplus}
				|X_j-(t,|f(t)|\,\omega)|^2 \leqslant \frac{|f(t_j)|^2}{16} + \frac{|f(t_j)|^2}{4} < |f(t_j)|^2.
			\end{equation}
			We now estimate $U_+(X_j)$ as follows
			\begin{equation*}
				\begin{aligned} 
					& \sgn[f(t_j)]\, U_+(X_j) \\[2mm]
					& = \frac{1}{|\s^4|}\int_{A_+} \sgn[f(t_j)]\, f(t)\sqrt{1+f'(t)^2} \left(\int_{\s^4} \frac{ \dd\omega}{|X_j-(t,|f(t)|\omega)|^4} \right)\dd t \\[2mm]
					& \geqslant \frac{1}{|\s^4|}\int_{t_j-|f(t_j)|/4}^{t_j+|f(t_j)|/4} \sgn[f(t_j)]\, f(t)\sqrt{1+f'(t)^2} \left(\int_{E_{\omega_0}} \frac{ \dd\omega}{|X_j-(t,|f(t)|\omega)|^4} \right)\dd t.
				\end{aligned} 
			\end{equation*}
			Using \eqref{eq:onAplus} and \eqref{eq:ineq_Xj-Aplus}, we obtain 
			\begin{equation}\label{eq:same_sign_lower_bound}
				\begin{aligned} 
					\sgn[f(t_j)]\, U_+(X_j) & \geqslant \frac{1}{|\s^4|}\int_{t_j-|f(t_j)|/4}^{t_j+|f(t_j)|/4} \frac{3\, |f(t_j)|}{4}
					\left(\int_{E_{\omega_0}} \frac{ \dd\omega}{|f(t_j)|^4}
					\right)\dd t \\[2mm]
					& \geqslant \frac{|f(t_j)|\,\beta}{2\, |\s^4|}\, \frac{3}{4\, |f(t_j)|^3} \\[2mm]
					& \geqslant \frac{c_0}{|f(t_j)|^2}.
				\end{aligned} 
			\end{equation}
			We now turn to $U_-(X_j)$. If $t\in A_-$, then it holds by \eqref{eq:small-slope} and a computation similar to \eqref{eq:IT}
			\begin{equation*}
				|f(t_j)|+|f(t)|=|f(t_j)-f(t)| \leqslant \eps\, |t-t_j|.
			\end{equation*}
			Hence, we obtain the following two inequalities for any $t\in A_-$
			\begin{equation}\label{eq:ineq-Aminus}
				|t-t_j|\geqslant \frac{|f(t_j)|}{\eps}, \qquad \text{ and }\qquad |f(t)|\leq\eps|t-t_j|. 
			\end{equation} 
			Using \eqref{eq:small-slope}, we obtain 
			\begin{equation*}
				\begin{aligned} 
				\sgn[f(t_j)]\, U_-(X_j) & = \frac{1}{|\s^4|}\int_{A_-} \sgn[f(t_j)]\, f(t)\, \sqrt{1+f'(t)^2} \left(\int_{\s^4} \frac{ \dd\omega}{|X_j-(t,|f(t)|\omega)|^4} \right)\dd t \\[2mm]
				& =- \frac{1}{|\s^4|}\int_{A_-} |f(t)|\, \sqrt{1+f'(t)^2} \left(\int_{\s^4} \frac{ \dd\omega}{|X_j-(t,|f(t)|\omega)|^4} \right)\dd t.
				\end{aligned} 
			\end{equation*} 
			Therefore, we have 
			\begin{equation*} 
				\begin{aligned} 
				\sgn[f(t_j)]\, U_-(X_j) &\geqslant - \frac{1}{|\s^4|}
				\int_{ \left\{ |t-t_j|\geqslant \frac{ |f(t_j)| }{\eps} \right\} } 2\, \eps\, |t-t_j|\, \left(\int_{\s^4} \frac{\dd\omega}{|t-t_j|^4}\right) \,\dd t.
				\end{aligned} 
			\end{equation*} 
			By change of variable $u=t-t_j$, we obtain 
			\begin{equation}\label{eq:U-}
				\sgn[f(t_j)]\, U_-(X_j) \geqslant - 2\eps  \int_{ \left\{ |u|\geqslant \frac{|f(t_j)|}{\eps} \right\} } \frac{\dd u}{|u|^3}  
				= - \frac{2\eps^3}{|f(t_j)|^2}. 
			\end{equation}
			Since $U_J=U_+ +U_-$, the estimates \eqref{eq:U-t2}, \eqref{eq:same_sign_lower_bound}, and \eqref{eq:U-} imply that
			\begin{equation*}
				\sgn[f(t_j)]\, U(X_j) \geqslant \frac{ c_0-2\eps^3}{ |f(t_j)| }-C_0.
			\end{equation*}
			This is \eqref{eq:lower-bound}.
		\end{proof} 
		
		By \eqref{eq:Ux0-0} and \eqref{eq:choice-eps}, we have
		\begin{equation}\label{eq:exterior_point_contradiction}
			0=\sgn[f(t_j)]\, U(X_0) \geqslant \frac{c_0}{2\, |f(t_j)|^2}-C_0.
		\end{equation}
		Letting $j\to\infty$, we have $f(t_j)\to 0$ and we obtain a contradiction.
	\end{proof}

	\begin{remark}
		The map $U$ may seem surprising at first sight. However, there is a natural way to this potential. On the one hand, the formula \eqref{eq:EL-simple} recalls the Green kernel of the Laplacian in dimension 6. On the other hand, it also recalls the Laplacian of the Green kernel of the bi-Laplacian in dimension 6. We consider the map (this is a 6-dimensional extension of the potential already introduced in \cite{FHW1994})
		\begin{equation*}
			\forall X\in\R^6,\qquad \Phi_{\Gamma}(X) \coloneqq \frac{1}{16\pi^3}\int_{\Gamma} \frac{\dd z}{|X-z|^2}.
		\end{equation*}
		Then $\Phi_{\Gamma}$ satisfies  $(\lap_{\R^6})^2\Phi_{\Gamma} = \delta_{\Gamma}$. Up to the normalisations used at the beginning of the proof of Theorem \ref{thm:unique-continuation}, restricting $\Phi$ to $X=(a,b,0,0,0,0)$, the Euler--Lagrange equation \eqref{eq:EL-simple} reduces to $\dr_b\Phi|_{b=0}=0$ for $|a|$ large enough. The map $U$ is the harmonic extension of $\dr_b\Phi|_{b=0}$.
	\end{remark}

    \section{\texorpdfstring{Solutions to $\Pc^1_{\IL}$ and $\Pc^1_{\Mob}$ are distinct}{Solutions to P1IL and P1Mob are distinct}}\label{sec:distinct}

    In this section, we prove the following result.
    \begin{proposition}\label{pr:distinct}
        Consider $n$ distinct points $Z\coloneqq (z_1,\ldots,z_n)\in\C^n$ and $\tau$ a relative isotopy class of curves passing through $Z$. Assume that the solutions $\gamma_{\IL}$ to $\Pc^1_{\IL}$ and $\gamma_{\Mob}$ to $\Pc^1_{\Mob}$ parametrize the same Jordan curve $\Gamma\subset \C$. Then $\Gamma$ is a round circle.
    \end{proposition}

    For the rest of the section, we keep the notations of Proposition \ref{pr:distinct}. We fix 
    \begin{equation}\label{eq:regu-Gamma}
        1/2<\alpha<1 \text{ and we know  } \Gamma \text{ is of class }\Cd^{1,\alpha} \text{ by Theorem~\ref{th:main-Mob}}.
    \end{equation}
    We let $f\in \Cd^{1,\alpha} (\D; \Omega)$ be a biholomorphic map with $\Omega$ being the bounded connected component of $\C\setminus \Gamma$. We consider the following expansion of the kernel $(2\ii \pi)^{-1} K$ of the fermionic Grunsky operator defined in \eqref{eq:def-Grunsky}-\eqref{eq:kernelK}, for $z,w\in\D$,
    \begin{equation}\label{eq:decompo-K}
        (2\ii \pi)^{-1} K(z,w) = \sum_{n,m\geqslant 0} K^f_{nm}\, z^n\, w^m.
    \end{equation}
    In order to shorten the notation, we denote
    \begin{equation*}
        G_f \coloneqq {\rm Gr}_f^{1/2}.
    \end{equation*}
    Let $(f_{\eps})_{\eps\in(-1,1)}$  be a smooth family of conformal maps with $Z\subset f_{\eps}(\s^1)$. 
    
    By Proposition \ref{pr:Mob-Grunsky}, the Euler--Lagrange equation associated to $\Pc^1_{\Mob}$ can be written as 
    \begin{equation}\label{eq:EL-Mob-op}
        0 = \Re\tr\left[ G_f^*\, \left.\frac{\dd}{\dd \eps}\right|_{\eps=0} G_{f_{\eps}} \right] .
    \end{equation}

    On the other hand, we have from Theorem \ref{th:Loewner-Grunsky}, that the Euler--Lagrange equation associated to $\Pc^1_{\IL}$ is given by 
    \begin{equation}\label{eq:EL-Loewner-op}
        0= \Re\tr\left[ (I+G_f^*G_f)^{-1}G_f^*\, \left.\frac{\dd}{\dd \eps}\right|_{\eps=0} G_{f_{\eps}} \right].
    \end{equation}
    Taking the difference between \eqref{eq:EL-Mob-op} and \eqref{eq:EL-Loewner-op}, we obtain 
    \begin{equation*}
        0 = \Re\tr\left[ \left(I - (I+G_f^*G_f)^{-1}\right)G_f^* \, \left.\frac{\dd}{\dd \eps}\right|_{\eps=0} G_{f_{\eps}} \right] = \Re\tr\left[G_f^*G_f(I+G_f^*G_f)^{-1}G_f^*\, \left.\frac{\dd}{\dd \eps}\right|_{\eps=0} G_{f_{\eps}} \right].
    \end{equation*}
    Dropping the condition $Z\subset f_{\eps}(\s^1)$, we obtain coefficients $\beta_1,\ldots,\beta_n\in\C$ such that, for $p_j\in\s^1$ such that $f(p_j)=z_j$, and for \emph{any} smooth variation $(f_{\eps})_{\eps\in(-1,1)}$ of conformal maps
    \begin{equation}\label{eq:new-EL}
        \sum_{j=1}^n \Re\left( \beta_j\, \left. \frac{\dd f_{\eps}}{\dd \eps} \right|_{\eps=0}(p_j) \right) = \Re\tr\left[ G_f^*G_f(I+G_f^*G_f)^{-1}G_f^*\, \left.\frac{\dd}{\dd \eps}\right|_{\eps=0} G_{f_{\eps}}   \right].
    \end{equation}
    We will prove that all the coefficients $\beta_j$ are equal to $0$ if $\Gamma$ solves both $\Pc^1_{\IL} (Z,\tau)$ and $\Pc^1_{\Mob} (Z,\tau)$. To this end, we first show that the kernel of the operator $G_f^*G_f(I+G_f^*G_f)^{-1}G_f^*$ is Hölder continuous. 
    \begin{claim}\label{cl:reg-kernel}
        Let $\alpha$ be as in \eqref{eq:regu-Gamma}. The kernel $a\colon\s^1\times\s^1\to\C$ of the operator 
        \begin{equation}\label{eq:def-A}
            A \coloneqq G_f^*G_f(I+G_f^*G_f)^{-1}G_f^*
        \end{equation} satisfies, for any $x,y\in\s^1$,
        \begin{subnumcases}{}
            a(x,y) = -a(y,x), \label{eq:asym-a}\\
            |a(x,y_1)-a(x,y_2)|\lesssim |y_1-y_2|^{\alpha-\frac{1}{4}}. \label{eq:holder-a}
        \end{subnumcases}
    \end{claim}
    \begin{proof}
        The symmetry of $a$ follows from the antisymmetry of $K$, which implies that $G_f^{\top} = -G_f$. Moreover, we have $G_f(I+G_f^*G_f)=(I+G_f G_f^*)G_f$. Hence, there holds
    \begin{equation*}
    \begin{aligned} 
        \left[ G_f^*G_f(I+G_f^*G_f)^{-1}G_f^*\right]^{\top} & = (-G_f)^*(I+G_f^*G_f)^{-1} (-G_f)(-G_f)^* \\[2mm]
        & = -G_f^*G_f(I- G_f^*G_f)^{-1} G_f^*.
    \end{aligned} 
    \end{equation*} 
    This implies \eqref{eq:asym-a}.
    
        By Theorem \ref{th:main-Mob}, the curve $\Gamma$ is of class $\Cd^{1,\alpha}$. Hence the biholomorphism $f\colon\D\to\Omega$ is of class $\Cd^{1,\alpha}$ up to the boundary and its inverse $f^{-1}\colon \Omega\to \D$ is also $\Cd^{1,\alpha}$ up to the boundary. For each $z,y\in\s^1$, there holds 
        \begin{equation*}
        \begin{aligned} 
            \left| K(z,y) \right|
            & =\left| \frac{  \sqrt{ f'(z)}\, \sqrt{f'(y)}  }{ f(z)-f(y) } -\frac{1}{z-y} \right| \\[2mm]
            & = \frac{1}{|f(z)-f(y)|}\left| \sqrt{ f'(z)}\, \sqrt{f'(y)}  - \frac{f(z)-f(y)}{z-y}\right| \\[2mm]
            & = \frac{1}{|f(z)-f(y)|}\left| \sqrt{ f'(z)}\, \left( \sqrt{f'(y)} - \sqrt{f'(z)}\right) + f'(z)  - \frac{f(z)-f(y_j)}{z-y}\right|.
        \end{aligned} 
        \end{equation*}
        There holds 
        $$
            |f(z)-f(y)| \geqslant \left\| (f^{-1})' \right\|_{\Ll^{\infty}(\Omega)}^{-1} |z-y|. 
        $$
        We obtain the upper bound 
        \begin{equation}\label{eq:zy-small}
            \left| K(z,y) \right| \aleq \frac{1}{ |z-y| } \left( |z-y|^{\alpha} + |z - y|^{\alpha} \right) \aleq  |z-y|^{\alpha-1}.
        \end{equation}
        Since $\alpha-1>-1/2$ by assumption, we have that $\sup_{z\in\s^1} \|K(z,\cdot)\|_{\Ll^2(\s^1)}<\infty$. Hence the kernel $b$ of the operator $G_f^*G_f(I+G_f^*G_f)^{-1}=G_f^*(I+G_f^*G_f)^{-1}G_f$ satisfies, by Cauchy--Schwarz
        \begin{equation*}
        \begin{aligned} 
            & \|b\|_{\Ll^{\infty}(\s^1\times\s^1)} \\[2mm]
            & \leqslant 4\pi^2 \left(\sup_{z\in \s^1} \|K(z,\cdot)\|_{\Ll^2(\s^1)}\right)\, \left(\sup_{w\in \s^1} \left\| (I+G_f^*G_f)^{-1}K(\cdot,w) \right\|_{\Ll^2(\s^1)}\right) \\[2mm]
            & \leqslant 4\pi^2 \left(\sup_{z\in \s^1} \|K(z,\cdot)\|_{\Ll^2(\s^1)}\right)\,  \left\| (I+G_f^*G_f)^{-1} \right\|_{\Ll^2(\s^1)\to\Ll^2(\s^1)}\, \left(\sup_{w\in \s^1} \|K(\cdot,w)\|_{\Ll^2(\s^1)}\right) \\[2mm]
            & <\infty.
        \end{aligned} 
        \end{equation*}
        Thus, we have
        \begin{equation}\label{eq:a-diff}
        \begin{aligned} 
            |a(x,y_1) - a(x,y_2)| & =  \frac{1}{2\pi} \left| \int_{\s^1} b(x,z)\, \left[ \overline{K(z,y_1)} -\overline{K(z,y_2)}\right]\, |\dd z| \right| \\[2mm]
            & \aleq \int_{\{|z-y_2|<2|y_1-y_2|\}} \left| K(z,y_1)-K(z,y_2) \right|\, |\dd z| \\[2mm]
            &\quad + \int_{\{|z-y_2|>2|y_1-y_2|\}} \left| K(z,y_1)-K(z,y_2) \right|\, |\dd z|.
        \end{aligned} 
        \end{equation}
        We deduce from \eqref{eq:zy-small} that 
        \begin{equation}\label{eq:a-diff1}
        \begin{aligned} 
            & \int_{\{|z-y_2|<2|y_1-y_2|\}} \left| K(z,y_1)-K(z,y_2) \right|\, |\dd z| \\[2mm]
            & \aleq \int_{\{|z-y_2|<2|y_1-y_2|\}} \left( |z-y_1|^{\alpha-1} + |z-y_2|^{\alpha-1} \right)\, |\dd z| \\[2mm]
            & \aleq |y_1-y_2|^{\alpha}.
        \end{aligned} 
        \end{equation}
        If now $|z-y_2|>2|y_1-y_2|$, then we also have $|z-y_1|>|y_1-y_2|$. Hence, 
        \begin{align} 
            & \left| K(z,y_1)-K(z,y_2) \right| \nonumber \\[2mm]
            & = \left| \frac{  \sqrt{ f'(z)}\, \sqrt{f'(y_1)}  }{ f(z)-f(y_1) } -\frac{1}{z-y_1} -\frac{  \sqrt{ f'(z)}\, \sqrt{f'(y_2)}  }{ f(z)-f(y_2) } + \frac{1}{z-y_2} \right| \nonumber \\[2mm]
            & \leqslant |f'(z)|^{1/2} \left| \frac{ \sqrt{f'(y_1)} - \sqrt{f'(y_2)} }{ f(z)-f(y_1) } \right| \label{eq:t1}\\[2mm]
            &\quad +  \left| \sqrt{f'(z)}\, \sqrt{f'(y_2)}\left(   \frac{1}{f(z)-f(y_1)} - \frac{1}{f(z)-f(y_2)} \right) + \frac{ y_2-y_1 }{ (z-y_2)(z-y_1) } \right|. \label{eq:t2}
        \end{align}
        We estimate \eqref{eq:t1} as
        \begin{equation}\label{eq:t21}
            |f'(z)|^{1/2} \left| \frac{ \sqrt{f'(y_1)} - \sqrt{f'(y_2)} }{ f(z)-f(y_1) } \right| \aleq \frac{|y_1-y_2|^{\alpha}}{ |z-y_1| }  .
        \end{equation}
        We now rewrite \eqref{eq:t2} as follows 
        \begin{equation*}
            \begin{aligned}
                & \sqrt{f'(z)}\, \sqrt{f'(y_2)}\left(   \frac{1}{f(z)-f(y_1)} - \frac{1}{f(z)-f(y_2)} \right) + \frac{ y_2-y_1 }{ (z-y_2)(z-y_1) } \\[2mm]
                & =\sqrt{f'(z)}\, \sqrt{f'(y_2)} \, \frac{f(y_1)-f(y_2)}{(f(z)-f(y_2))(f(z)-f(y_1))}  + \frac{ y_2-y_1 }{ (z-y_2)(z-y_1) } \\[2mm]
                & = \left( K(z,y_2) + \frac{1}{z-y_2} \right) \, \frac{f(y_1)-f(y_2)}{f(z)-f(y_1)}  + \frac{ y_2-y_1 }{ (z-y_2)(z-y_1) } \\[2mm]
                & = K(z,y_2)\, \frac{f(y_1)-f(y_2)}{f(z)-f(y_1)} + \frac{1}{z-y_2}\left[ \frac{f(y_1)-f(y_2)}{f(z)-f(y_1)} + \frac{y_2-y_1}{z-y_1} \right] \\[2mm]
                & = K(z,y_2)\, \frac{f(y_1)-f(y_2)}{f(z)-f(y_1)} + \frac{y_2-y_1}{(z-y_2)(f(z)-f(y_1))}\left[ \frac{f(y_1)-f(y_2)}{y_2-y_1} + \frac{f(z)-f(y_1)}{z-y_1} \right] \\[2mm]
                & = K(z,y_2)\, \frac{f(y_1)-f(y_2)}{f(z)-f(y_1)} \\[2mm]
                &\quad + \frac{y_2-y_1}{(z-y_2)(f(z)-f(y_1))}\int_0^1 \Big( f'((1-t)y_1+tz)-f'((1-t)y_1+ty_2)\Big)\dd t.
            \end{aligned}
        \end{equation*}
        Therefore, we obtain using \eqref{eq:zy-small}
        \begin{equation*}
            \begin{aligned}
                & \left| \sqrt{f'(z)}\, \sqrt{f'(y_2)}\left(   \frac{1}{f(z)-f(y_1)} - \frac{1}{f(z)-f(y_2)} \right) + \frac{ y_2-y_1 }{ (z-y_2)(z-y_1) } \right| \\[2mm]
                & \aleq |z-y_2|^{\alpha-1} \frac{|y_1-y_2|}{|z-y_1|} + \frac{|y_1-y_2|}{|z-y_2|\, |z-y_1|}\, |z-y_2|^{\alpha}.
            \end{aligned}
        \end{equation*}
        We estimate the denominator of each term using the inequalities 
        $$
            |z-y_1|>|z-y_2| - |y_1-y_2|\geqslant \frac{1}{2} |z-y_2|.
        $$
        We obtain 
        \begin{equation}\label{eq:t22}
            \begin{aligned} 
                & \left| \sqrt{f'(z)}\, \sqrt{f'(y_2)}\left(   \frac{1}{f(z)-f(y_1)} - \frac{1}{f(z)-f(y_2)} \right) + \frac{ y_2-y_1 }{ (z-y_2)(z-y_1) } \right| \\[2mm]
                &\aleq |z-y_2|^{\alpha-2}\, |y_1-y_2|.
            \end{aligned}
        \end{equation}
        Since $\alpha-2<-1$, we obtain the following estimate by combining \eqref{eq:t21}, \eqref{eq:t22}, 
        \begin{equation*}
        \begin{aligned} 
            & \int_{\{|z-y_2|>2|y_1-y_2|\}} \left| K(z,y_1)-K(z,y_2) \right|\, |\dd z| \\[2mm]
            & \aleq |y_1-y_2|^{\alpha}\,  |\log|y_1-y_2|| + |y_1-y_2|^{\alpha} \\[2mm]
            & \aleq |y_1-y_2|^{\alpha-\frac{1}{4}}.
        \end{aligned} 
        \end{equation*}
        We conclude the proof of \eqref{eq:holder-a} using the above estimate with \eqref{eq:a-diff} and \eqref{eq:a-diff1}.
    \end{proof}

    We now prove that all the coefficients $\beta_j$ in \eqref{eq:new-EL} are actually $0$.

    \begin{claim}
        Let $j_0\in\{1,\ldots,n\}$. Then we have $\beta_{j_0}=0$.
    \end{claim}
    \begin{proof} 
    Let $\lambda\in\Omega^*$, the unbounded connected component of $\C\setminus \Gamma$. Let $n_{j_0}$ be the outward-pointing normal vector to $\Gamma$ at the point $z_{j_0}$. We choose $\lambda$ of the form $\lambda=z_{j_0}+ t n_{j_0}$ for $t>0$ that we will send to $0$ eventually. We consider the variation 
    \begin{equation*}
        f_{\eps}(z) = f(z) + \frac{\eps}{f(z)-\lambda}.
    \end{equation*}
    We then have 
    \begin{equation*}
    \begin{aligned} 
        & \left.\frac{\dd}{\dd \eps}\right|_{\eps=0} \frac{ \sqrt{f_{\eps}'(z)}\, \sqrt{f_{\eps}'(w)} }{f_{\eps}(z) - f_{\eps}(w) } \\[2mm]
        & = \frac{ \sqrt{f'(z)}\, \sqrt{f'(w)} }{f(z) - f(w) } \left[ \frac{\dot{f}'(z)}{2\, f'(z)} + \frac{\dot{f}'(w)}{ 2\, f'(w)} - \frac{\dot{f}(z) - \dot{f}(w)}{f(z)-f(w)} \right] \\[2mm]
        & = \frac{ \sqrt{f'(z)}\, \sqrt{f'(w)} }{f(z) - f(w) } \left[ \frac{-1}{2\, (f(z)-\lambda)^2} + \frac{-1}{ 2\, (f(w)-\lambda)^2 } + \frac{1}{(f(z)-\lambda)(f(w)-\lambda) } \right] \\[2mm]
        & = \frac{ \sqrt{f'(z)}\, \sqrt{f'(w)} }{f(z) - f(w) }\  \frac{ (f(z)-\lambda)(f(w)-\lambda) - \frac{1}{2}(f(w)-\lambda)^2 - \frac{1}{2}(f(z)-\lambda)^2 }{(f(z)-\lambda)^2 (f(w)-\lambda)^2 } \\[2mm]
        & = \frac{ \sqrt{f'(z)}\, \sqrt{f'(w)} }{f(z) - f(w) }\  \frac{  - \frac{1}{2}(f(z)-f(w))^2 }{(f(z)-\lambda)^2 (f(w)-\lambda)^2 } \\[2mm]
        & =  \sqrt{f'(z)}\, \sqrt{f'(w)} \,  \frac{  - \frac{1}{2}(f(z)-f(w)) }{(f(z)-\lambda)^2 (f(w)-\lambda)^2 } \\[2mm]
        & = -\frac{1}{2} \left[ \frac{\sqrt{f'(z)}\, \sqrt{f'(w)} }{ (f(z)-\lambda)(f(w)-\lambda)^2 } - \frac{ \sqrt{f'(z)}\, \sqrt{f'(w)} }{ (f(z)-\lambda)^2(f(w)-\lambda) } \right].
    \end{aligned} 
    \end{equation*}
    Plugging this into \eqref{eq:new-EL}, we obtain thanks to \eqref{eq:asym-a}
    \begin{equation*}
    \begin{aligned} 
        \sum_{j=1}^n\Re\left( \frac{\beta_j}{z_j-\lambda} \right) 
        & = \Re \left( \frac{1}{2\ii\pi} \int_{\s^1\times\s^1} \overline{a(x,y)}\, \frac{\sqrt{f'(x)}\, \sqrt{f'(y)}}{ (f(x)-\lambda)^2(f(y)-\lambda) }\, |\dd x|\, |\dd y| \right) \\[2mm]
        & = \Re \left( \frac{1}{2\ii \pi} \int_{\s^1} \frac{\sqrt{f'(y)} }{ f(y)-\lambda } \left( \int_{\s^1} \frac{ \overline{a(x,y)} }{\ii x \sqrt{f'(x)}}\, \frac{f'(x)}{(f(x)-\lambda)^2}\, \dd x \right)\, |\dd y| \right) .
    \end{aligned} 
    \end{equation*}
    Moreover, there holds
    \begin{equation*}
        \int_{\s^1} \frac{f'(x)}{(f(x)-\lambda)^2}\, \dd x  = \int_{\s^1} \dd \left( \frac{1}{f-\lambda} \right) =0.
    \end{equation*}
    Therefore, we have 
    \begin{equation*}
        \begin{aligned} 
        & \sum_{j=1}^n\Re\left( \frac{\beta_j}{z_j-\lambda} \right) \\[2mm]
        & = \Re\int_{\s^1} \frac{\sqrt{f'(y)} }{ f(y)-\lambda } \left( \int_{\s^1} \left[ \frac{ \overline{a(x,y)} }{\ii x \sqrt{f'(x)}} - \frac{ \overline{a(p_{j_0},y)} }{\ii p_{j_0} \sqrt{f'(p_{j_0})}} \right]\, \frac{f'(x)}{(f(x)-\lambda)^2}\, \dd x \right)\, |\dd y|.
        \end{aligned} 
    \end{equation*}
    We consider the change of variable $x=p_{j_0}\ee^{\ii \theta}$ and $y=p_{j_0} \ee^{\ii\vp}$. For $|\theta|$ small enough, there holds
    \begin{equation*}
        \begin{aligned} 
            |f(x)-\lambda| & = \left| f\left( p_{j_0}\ee^{\ii\theta} \right) - f(p_{j_0}) - tn_{j_0} \right| \\[2mm]
            & = \left| \int_0^{\theta} \ii p_{j_0}\ee^{\ii\tilde{\theta}} f'\left( p_{j_0}\ee^{\ii\tilde{\theta}}\right)\dd \tilde{\theta} - tn_{j_0} \right| \\[2mm]
            & = \left| \int_0^{\theta} \ii p_{j_0} \left[ \ee^{\ii\tilde{\theta}} f'\left( p_{j_0} \ee^{\ii\tilde{\theta}}\right) -f'(p_{j_0})\right]\dd \tilde{\theta} +\ii p_{j_0}\, \theta\, f'(p_{j_0})- tn_{j_0} \right|.
        \end{aligned} 
    \end{equation*}
    For $\theta$ small enough, the dominant terms are the last two terms, since $\ii p_{j_0} f'(p_{j_0})$ is orthogonal to $n_{j_0}$. Therefore, we have 
    \begin{equation*}
        |f(x)-\lambda| \geq \sqrt{\theta^2+t^2} - \theta^{1+\alpha}\, \|f'\|_{\Cd^{0,\alpha}}.
    \end{equation*}
    Together with Claim \ref{cl:reg-kernel}, we obtain 
    \begin{equation*}
        \begin{aligned}
             \left|\sum_{j=1}^n\Re\left( \frac{\beta_j}{z_j-\lambda} \right) \right| & \aleq \int_{[-\pi,\pi]^2} \frac{ |p_{j_0}\ee^{\ii\theta} -p_{j_0} |^{\alpha-\frac{1}{4}} }{ |f(p_{j_0}\ee^{\ii\vp})- \lambda|\, |f(p_{j_0}\ee^{\ii\theta})- \lambda|^2 }\, \dd\theta \, \dd \vp  \\[2mm]
             & \aleq \int_{[-\pi,\pi]^2} \frac{ |\ee^{\ii\theta} -1 |^{\alpha-\frac{1}{4}} }{ |f(p_{j_0}\ee^{\ii\vp})- \lambda|\, |f(p_{j_0}\ee^{\ii\theta})- \lambda|^2 }\, \dd\theta \, \dd \vp  \\[2mm]
             & \aleq \left(1+\int_{-1}^1 \frac{\dd\vp}{|\vp|+t}\right)\left( 1+ \int_{-1}^1 \frac{ |\theta|^{\alpha} }{ (|\theta|+t)^2 } \dd \theta\right) \\[2mm]
             & \aleq \left( 1 + \log(1+t)\right) \left( 1 + t^{\alpha-1} \right).
        \end{aligned}
    \end{equation*}
    Multiplying by $t$ and letting $t\to0$, we obtain 
    \begin{equation*}
        \Re\left(\frac{\beta_{j_0}}{n_{j_0}} \right) =0.
    \end{equation*}
    We can perform the same computation with $\eps/(f-\lambda)$ replaced by $\ii\eps/(f-\lambda)$ and we end up with $\Re(\beta_{j_0}(\ii n_{j_0})^{-1})=0$. Hence $\beta_{j_0}$ is orthogonal to both $1/n_{j_0}$ and $\ii/n_{j_0}$, which form a basis of $\C$. We deduce that $\beta_{j_0}=0$.
    \end{proof}
    
    Coming back to \eqref{eq:new-EL}, we obtain for all variations $(f_{\eps})_{\eps}$, that
    \begin{equation}\label{eq:new-EL2}
        0= \Re\tr\left[ 
        A \left.\frac{\dd}{\dd \eps}\right|_{\eps=0} G_{f_{\eps}}   \right]
    \end{equation}
    where  $A= G_f^*G_f(I+G_f^*G_f)^{-1}G_f^*$ is defined in \eqref{eq:def-A}.
    
    Let $r<1$ and $f_r(z) = r^{-1}f(rz)$ for $z\in\D$. We can now take the variation $f_{\lambda,r}(z)=f(z)+\lambda(f_r(z)-f(z))$, for $\lambda\in\C$ with $|\lambda|\aleq (1-r)^{-\alpha}$ for any $1/2<\alpha<1$ as in \eqref{eq:regu-Gamma} and $r$ close enough to $1$. Indeed, we have 
    \begin{equation*}
    \begin{aligned} 
        |f_{\lambda,r}(z)-f_{\lambda,r}(w)| & \geqslant |f(z)-f(w)| - \lambda \|f_r'-f'\|_{\Ll^{\infty}(\D)}|z-w| \\[2mm] 
        & \geqslant \left( \left\|(f^{-1})'\right\|_{\Ll^{\infty}(\D)}- \lambda \|f_r'-f'\|_{\Ll^{\infty}(\D)}  \right)|z-w|.
        \end{aligned} 
    \end{equation*}
    Moreover, there holds by Theorem \ref{th:main-Mob}
    \begin{equation*}
        \|f_r'-f\|_{\Ll^{\infty}(\D)} = \sup_{z\in\D} |f'(rz)-f'(z)| \aleq (1-r)^{\alpha}.
    \end{equation*}
    Hence, for $|\lambda|\leqslant t_0\coloneqq C(1-r)^{-\alpha}$ for some $C>0$ depending only on $f$, the map $f_{\lambda,r}$ is holomorphic and bi-Lipschitz, i.e., it is a biholomorphism on its image.
    
    We then have by Claim \ref{cl:reg-kernel} and the uniform estimate $|K_{f_{\lambda,r}}(x,y)| \aleq |x-y|^{-1}$ for $|\lambda|\leqslant t_0$,
    \begin{equation*}
    \begin{aligned}
        \left| \tr(A\, G_{f_{\lambda,r}}) \right| & \aleq \int_{\s^1\times \s^1} \frac{|a(x,y)|}{|x-y|}\, |\dd x|\, |\dd y| \\[2mm]
        & = \int_{\s^1\times \s^1} \frac{|a(x,y)-a(x,x)|}{|x-y|}\, |\dd x|\, |\dd y| \\[2mm]
        & \aleq \int_{\s^1\times \s^1}\frac{|\dd x|\, |\dd y|}{|x-y|^{5/4-\alpha}} <\infty. 
    \end{aligned} 
    \end{equation*}
    Moreover, the map $\big(\vp\colon  \lambda\mapsto \tr(A G_{f_{\lambda,r}}) \big)$ is analytic in the disk of radius $t_0$ with a critical point at $\lambda=0$. Hence, there holds
    \begin{equation*}
    \begin{aligned} 
        |\lambda|^{-2} \left| \tr(A G_{f_{\lambda,r}}) - \tr(A G_{f})\right| & = \left| \frac{\vp(\lambda) - \vp(0)}{\lambda^2} \right|  \\[2mm]
        & \aleq \sup_{|\lambda|=t_0} \frac{|\vp(\lambda)-\vp(0)|}{|\lambda|^2} \\[2mm]
        & \aleq \frac{1}{t_0^2} \left(\sup_{|\lambda|\leqslant t_0} |\tr(A G_{f_{\lambda,r}})| \right) \\[2mm]
        & \aleq (1-r)^{2\alpha}.
    \end{aligned} 
    \end{equation*}
    Taking $\lambda=1$, we obtain 
    \begin{equation}\label{eq:upper1}
        \left| \tr(A (G_{f_r}- G_f))\right| \aleq (1-r)^{2\alpha}.
    \end{equation}

    \begin{claim}\label{cl:G-zero}
        There holds 
        $$
            \|G_f \|_{\rm HS}^2 - \log\det(1+G_f^*G_f) \aleq (1-r)^{2\alpha-1}.
        $$
    \end{claim}
    \begin{proof}
        The coefficients in \eqref{eq:decompo-K} are given by
    \begin{equation*}
        K^{f_{r}}_{nm} = r^{1+n+m} K^f_{nm}.
    \end{equation*}
        We define the linear operator $D_r\colon \Ll^2(\s^1)\to\Ll^2(\s^1)$ by its action on each $e_n\colon z\mapsto z^n$ for $n\in\Z$ as
        \begin{equation*}
            D_r (e_n) \coloneqq r^n\, e_n.
        \end{equation*}
        We then have $\|D_r\|_{\rm op} \leqslant 1$ and $G_{f_r} = rD_r G_f D_r$. Let $j\in\N$. We denote $P_j\colon \Ll^2(\s^1)\to \Ll^2(\s^1)$ the projection on the $j-1$ first eigenfunctions of $G_f^*G_f$. Then the $j^{th}$ eigenvalue $s_j(G_{f_r}^*G_{f_r})$ satisfies
        \begin{equation*}
            s_j(G_{f_r}^*G_{f_r}) = r^2\|D_r( G_f-P_j)D_r\|^2_{\rm op} \leqslant r^2\|G_f-P_j\|_{\rm op}^2 = r^2\, s_j(G_f^*G_f).
        \end{equation*}
        We consider the function $\psi(x)\coloneqq x-\log(1+x)$ for $x\geqslant 0$. Then, there holds, $\psi'(x) = 1-1/(1+x)$ and $\psi''(x)=1/(1+x)^2>0$. Hence $\psi$ is convex and we obtain, for $1>r>0$, 
        \begin{equation*}
            \psi(r^2x) = \psi(r^2x+(1-r^2)x) \leqslant r^2 \psi(x) + (1-r^2)\psi(0) = r^2\psi(x).
        \end{equation*}
        Thus, we obtain 
        \begin{equation*}
        \begin{aligned} 
            g(r) & \coloneqq \|G_{f_r}\|_{\rm HS}^2 - \log\det(1+G_{f_r}^*G_{f_r}) \\[2mm]
            & = \sum_{j\geqslant 0} \psi\big( s_j(G_{f_r}^*G_{f_r}) \big) \\[2mm]
            &\leqslant r^2\left( \|G_f \|_{\rm HS}^2 - \log\det(1+G_f^*G_f)\right)  =r^2 g(1).
        \end{aligned} 
        \end{equation*}
        Hence, we have 
        \begin{equation}\label{eq:upper2}
            (1-r^2)\, g(1) \leqslant g(1)-g(r)
        \end{equation}
        Moreover, the function $F\colon X\mapsto \left( \|X\|_{\rm HS}^2 - \log\det(1+X^*X) \right) $ is convex. Its derivative at the point $X=G_f$ in a direction $H$ is computed similarly as in \eqref{eq:derivative-logdet}, \eqref{eq:EL-Mob-op} and \eqref{eq:EL-Loewner-op} (the coefficient $2$ has been omitted in \eqref{eq:EL-Mob-op} and \eqref{eq:EL-Loewner-op} for simplicity) in terms of the operator $A$ defined in \eqref{eq:def-A}:
        \begin{equation*}
            \dd F(G_f)\cdot H = 2\Re\tr\left[ A\, H\right].
        \end{equation*}
        Since the graph of $F$ lies above its tangent at the point $X=G_f$, there holds
        \begin{equation*}
            g(r) \geqslant g(1)+ 2\Re\tr\left[ A (G_{f_r}-G_f)\right].
        \end{equation*}
        Thanks to \eqref{eq:upper1}, we obtain 
        \begin{equation*}
            g(1)-g(r) \leqslant 2\Re\tr\left[ A (G_f-G_{f_r})\right] \aleq (1-r)^{2\alpha}.
        \end{equation*}
        Then \eqref{eq:upper2} implies that 
        $g(1) \aleq (1-r)^{2\alpha-1}$ which completes the proof. 
    \end{proof}
    
    We let $r\to 1$ in Claim \ref{cl:G-zero} and obtain 
    \begin{equation*}
        \|G_f \|_{\rm HS}^2 - \log\det(1+G_f^*G_f) \leqslant 0.
    \end{equation*}
    If $(\lambda_j)_{j\geqslant 0}$ denotes the eigenvalues of $G_f^*G_f$ and $\psi(x)\coloneqq x-\log(1+x)$, we obtain 
    \begin{equation*}
        \sum_{j\geqslant 0} \psi(\lambda_j) \leqslant 0.
    \end{equation*}
    Since $\psi(x)>0$ for $x>0$, we deduce that $\lambda_j=0$ for any $j\geqslant0$. Hence, we have $G_f=0$.    
    By Claim \ref{cl:der-K}, this implies that $\Sr f=0$ and hence, $f$ is a Möbius transformation. Thus, $\Gamma$ is a round circle. This concludes the proof of Proposition \ref{pr:distinct}.

	\appendix
	
	\section{\texorpdfstring{Dirichlet problem for the \(\frac{1}{2}\)-Laplacian in dimension $1$}{Dirichlet problem for the half-Laplacian in dimension 1}}\label{sec:CZ}
	
	In this section, we develop the regularity up to the boundary for the Dirichlet problem with $\frac{1}{2}$-Laplacian. The Green kernel for the associated Dirichlet problem on $(-1,1)$ is introduced in \Cref{sec:Green}. We prove some pointwise estimates on its Green kernel and its gradient in \Cref{sec:pointwise_Green} and \Cref{sec:pointwise_dGreen}. We deduce some regularity properties in \Cref{sec:Lap_Green} on the fractional Laplacian of the Green kernel. We obtain global regularity estimates for the Dirichlet problem in Section \ref{sec:global_reg}.
	
	\subsection{The Green kernel on $(-1,1)$.}\label{sec:Green}
	Given $z\in(-1,1)$, we denote
	\begin{equation}
		\delta(z)\coloneqq \dist(z,\{-1,1\}) = \min(|z-1|,|z+1|) = 1-|z|.
	\end{equation}
	By \cite[Theorem 3.1]{B2016} (see also \cite[Corollary 4]{BGR1961} or \cite[Definition 1.1 and Theorem 1.4]{CS1998}), it holds for $u=-1$ and $v=1$
	\begin{equation*}
		f(x) = \int_{-1}^1 G(x,y)\, (-\lap)^{1/2} f(y)\ \dd y,
	\end{equation*}
	with
	\begin{equation}\label{eq:Green}
		\begin{aligned} 
			G(x,y) & = \frac{1}{\pi}\, \log\left( \frac{1-x\, y + \sqrt{(1-x^2)\, (1-y^2)} }{|x-y|} \right) \\[2mm]
			&   = \frac{1}{2\pi}\, \log\left( \frac{1-x\, y + \sqrt{(1-x^2)\, (1-y^2)} }{ 1-x\, y - \sqrt{(1-x^2)\, (1-y^2) } } \right) .
		\end{aligned} 
	\end{equation} 
	\begin{proof}[Proof of the second equality]
		We denote $A \coloneqq 1-x\, y$ and $B\coloneqq  \sqrt{(1-x^2)\, (1-y^2)}$. We have 
		\begin{equation}\label{eq:A-B}
		\begin{aligned}
			A^2 - B^2 & = 1-2\, x\, y + x^2\, y^2 - (1-x^2)\, (1-y^2) \\[1mm]
			& = 1-2\, x\, y + x^2\, y^2 - (1-x^2 - y^2 + x^2\, y^2) \\[1mm]
			& = (x-y)^2.
		\end{aligned}
		\end{equation} 
		Hence, we obtain 
		\begin{equation*}
			(A-B)(A+B) = |x-y|^2.
		\end{equation*}
		That is to say, it holds
		\begin{equation*}
			\frac{A+B}{|x-y|} = \frac{|x-y|}{A-B} = \sqrt{\frac{A+B}{A-B}}.\qedhere
		\end{equation*}
	\end{proof}
	
	\subsection{Pointwise estimate of the Green kernel.}\label{sec:pointwise_Green}
	
	The formula \eqref{eq:Green} is not necessarily simple to work with. We provide some pointwise estimates in this section.
	
	\begin{lemma}\label{lm:log-G}
		For every $x,y\in(-1,1)$ with $x\neq y$, it holds
		\begin{equation*}
		G(x,y) \leqslant \frac{1}{\pi} \log\left( 1+4\frac{\sqrt{\delta(x)\delta(y)}}{|x-y|} \right).
		\end{equation*}
	\end{lemma}
	
	\begin{proof}
		We denote
		\begin{equation*}
		A\coloneqq 1-xy, \qquad \text{ and } \qquad 
		B \coloneqq \sqrt{(1-x^2)(1-y^2)},
		\qquad \text{ and } \qquad  d\coloneqq |x-y|.
		\end{equation*}
		Since $A>0$, it follows from \eqref{eq:A-B} that $A=\sqrt{d^2+B^2}$. 
		By \eqref{eq:Green}, we have
		\begin{equation}\label{eq:est-G1} 
			G(x,y) =\frac{1}{\pi} \log\left(\frac{A+B}{d}\right) =\frac{1}{\pi} \log\left( \frac{\sqrt{d^2+B^2}+B}{d} \right) =\frac{1}{\pi} \log\left( \sqrt{1+\left(\frac{B}{d}\right)^2} +\frac{B}{d} \right).
		\end{equation} 
		For $s\geqslant 0$, we have $\sqrt{1+s^2}\leqslant 1+s$. Applying this with $s=\frac{B}{d}$, we obtain
		\begin{equation*}
			G(x,y) \leqslant \frac{1}{\pi} 	\log\left(1+2\, \frac{B}{d}\right).
		\end{equation*}
		Moreover, it holds
		\begin{equation*}
			\delta(x) \leqslant \delta(x)(1+|x|) = (1-|x|)(1+|x|)= 1-x^2 \leqslant 2\, \delta(x).
		\end{equation*}
		Therefore, we have
		\begin{equation*}
			B=\sqrt{(1-x^2)(1-y^2)} \leqslant 2\, \sqrt{\delta(x)\, \delta(y)}.
		\end{equation*}
		Plugging this into \eqref{eq:est-G1}, we obtain the result.
	\end{proof}
	
	As a corollary of Lemma \ref{lm:log-G}, we obtain the following pointwise estimate useful near the boundary.
	\begin{lemma}\label{lm:boundary_G}
		There exists $C>0$ such that, for every $x,y\in(-1,1)$ with $x\neq y$,
		\begin{equation*}
		\frac{G(x,y)}{\sqrt{\delta(x)}} \leqslant \frac{C}{|x-y|^{1/2}}.
		\end{equation*}
	\end{lemma}
	
	\begin{proof}
		We denote
		\begin{equation*}
			\rho\coloneqq \delta(x), \qquad \text{ and } \qquad \sigma\coloneqq \delta(y), \qquad \text{ and }\qquad  d\coloneqq |x-y|.
		\end{equation*}
		Since the distance function is $1$-Lipschitz,
		\begin{equation*}
			\begin{cases} 
			|\delta(x)-\delta(y)| \leqslant |x-y|, \\[2mm]
			\delta(y) \leqslant \delta(x)+|x-y|.
			\end{cases} 
		\end{equation*}
		We distinguish two cases.
		
		\underline{Case 1: $|x-y|\leqslant \delta(x)$.}
		In this case, we have 
		\begin{equation*}
			\delta(y) \leqslant \delta(x) +|x-y| \leqslant 2\, \delta(x).
		\end{equation*}
		Consequently, it holds
		\begin{equation*}
			\sqrt{\delta(x)\, \delta(y)} \leqslant \sqrt{2}\,\delta(x).
		\end{equation*}
		We obtain from Lemma \ref{lm:log-G}
		\begin{equation*}		
			G(x,y) \leqslant \frac{1}{\pi} \log\left(1+4\sqrt{2}\, \frac{\delta(x)}{|x-y|}\right).
		\end{equation*}
		There exists $C>0$ such that for every $t\geqslant 0$, it holds $\log(1+t)\leqslant C\sqrt{t}$. Hence, we have
		\begin{equation*}
			G(x,y) \leqslant C\sqrt{\frac{\delta(x)}{|x-y|}}.
		\end{equation*}

		\underline{Case 2: $|x-y|>\delta(x)$.}
		In this case, it holds
		\begin{equation*}
			\delta(y) \leqslant \delta(x) + |x-y| < 2\, |x-y|.
		\end{equation*}
		It follows that
		\begin{equation*}
		\frac{\sqrt{\delta(x)\, \delta(y)}}{|x-y|} \leqslant \frac{\sqrt{2\, \delta(x)\, |x-y|}}{|x-y|} = \sqrt{\frac{2\,\delta(x)}{|x-y|}}.
		\end{equation*}
		Using $\log(1+t)\leqslant t$ for $t\geqslant 0$, we obtain 
		\begin{equation*}
		G(x,y) \leqslant C\frac{\sqrt{\delta(x)\, \delta(y)}}{|x-y|} \leqslant C\sqrt{\frac{\delta(x)}{ |x-y| }}.\qedhere
		\end{equation*}
	\end{proof}
	
	\subsection{Pointwise estimate of $\dr_x G$.}\label{sec:pointwise_dGreen}
	
	From the exact formula \eqref{eq:Green}, we derive a pointwise estimate on the derivatives of the Green function in \eqref{eq:gx_G}.
	\begin{claim}
		It holds
		\begin{equation}\label{eq:der-log1}
			\dr_x\log\left(1-xy + \sqrt{(1-x^2)(1-y^2)}\right) = \frac{ -y\sqrt{1-x^2} -x \sqrt{1-y^2} }{(1-xy)\sqrt{1-x^2} + (1-x^2)\, \sqrt{1-y^2}}.
		\end{equation}
	\end{claim}
	\begin{proof}
		On the one hand, we have
		\begin{equation}\label{eq:der-log2}
			\dr_x \log|x-y| = \frac{\g_x |x-y|}{|x-y|} = \frac{\mathrm{sign}(x-y)}{|x-y|} = \frac{1}{x-y}. 
		\end{equation}
		On the other hand, one has
		\begin{equation*}
			\dr_x\left(1-xy + \sqrt{(1-x^2)(1-y^2)}\right) = -y + \sqrt{1-y^2}\, \frac{-x}{\sqrt{1-x^2}}.
		\end{equation*}
		Hence, it holds
		\begin{equation*}
			\dr_x\log\left(1-xy + \sqrt{(1-x^2)(1-y^2)}\right) = \frac{ -y + \sqrt{1-y^2}\, \frac{-x}{\sqrt{1-x^2}} }{1-xy + \sqrt{(1-x^2)(1-y^2)}}.
		\end{equation*}
	\end{proof}
	
	It holds $G(x,y) = G(y,x)$ and we have from \eqref{eq:Green}, \eqref{eq:der-log1} and \eqref{eq:der-log2}
	\begin{equation*}
		\begin{aligned} 
			\dr_x G(x,y) & = \frac{1}{\pi}\ \left( \frac{ -y\sqrt{1-x^2} -x \sqrt{1-y^2} }{ (1-xy)\sqrt{1-x^2} + (1-x^2)\, \sqrt{1-y^2}} - \frac{1}{x-y} \right) \\[4mm]
			& = \frac{1}{\pi}\ \frac{ \big( -y\sqrt{1-x^2} -x \sqrt{1-y^2} \big)(x-y) - (1-xy)\sqrt{1-x^2} - (1-x^2)\, \sqrt{1-y^2} }{(x-y) \big( (1-xy)\sqrt{1-x^2} + (1-x^2)\, \sqrt{1-y^2} \big)} \\[4mm]
			& = \frac{1}{\pi}\ \frac{ \sqrt{1-x^2}\, (y(y-x) + xy-1) + \big( -x(x-y) - (1-x^2)\big)\, \sqrt{1-y^2} }{(x-y) \big( (1-xy)\sqrt{1-x^2} + (1-x^2)\, \sqrt{1-y^2} \big)} \\[4mm]
			& = \frac{1}{\pi}\ \frac{ -\sqrt{1-y^2}\, \big( \sqrt{1-x^2}\, \sqrt{1-y^2} + 1- xy\big) }{(x-y)\, \sqrt{1-x^2}\,  \big( (1-xy) +  \sqrt{(1-x^2)\,(1-y^2)} \big)} \\[4mm]
			& = \frac{1}{\pi}\, \frac{-\sqrt{1-y^2}}{(x-y)\, \sqrt{1-x^2}}.
		\end{aligned} 
	\end{equation*}
	Therefore, we obtain the following pointwise estimate
	\begin{equation}\label{eq:gx_G}
		\forall x,y\in(-1,1),\qquad |\dr_x G(x,y)| \aleq \frac{\sqrt{\delta(y)}}{|x-y|\, \sqrt{\delta(x)}}.
	\end{equation}
	
	\subsection{Fractional Laplacian of the Green kernel}\label{sec:Lap_Green}
	In this section, we use the previous estimates to derive pointwise estimates of $(-\lap)^{t/2}_x G(x,y)$ for $t\in(0,1)$. We decompose the estimates in the cases where $x,y\in(-1,1)$ and $x\not\in (-1,1)$ but $y\in(-1,1)$.
	
	\underline{Step 1: for $x,y\in(-1,1)$.} We prove the following estimate.
	\begin{proposition}\label{pr:lapt_G}
		Let $t\in(0,1)$. There exists $C_t>0$ such that, for every distinct $x,y\in (-1,1)$, it holds
		\begin{equation}\label{eq:lapt_G2}
			\left|(-\lap)^{t/2}_x G(x,y)\right| \leqslant
			C_t\left[ \frac{1+|\log|x-y||}{|x-y|^t} + \frac{\Phi_t(\delta(x))}{|x-y|^{1/2}} \right],
		\end{equation}
		where
		\begin{equation*}
		\Phi_t(r) \coloneqq  \begin{cases}
			1, & 0<t<1/2,\\[1mm]
			1+|\log r|, & t=\frac{1}{2},\\[1mm]
			r^{1/2-t}, & \frac{1}{2}<t<1.
		\end{cases}
		\end{equation*}
	\end{proposition}
	
	\begin{proof}
		Let $x,y\in (-1,1)$, with $x\neq y$, and denote
		\begin{equation*}
			a\coloneqq \frac{1}{4}\, \min\{\delta(x), |x-y| \}.
		\end{equation*}
		Since $0<a\leqslant \delta(x)/4$, it holds $(x-a,x+a)\subset (-1,1)$.
		
		Since $G(z,y)=0$ for $z\in\R\setminus (-1,1)$, we have
		\begin{align}
			\left|(-\Delta)^{t/2}_xG(x,y)\right|
			&\leqslant \int_{(x-a,x+a)} \frac{|G(x,y)-G(z,y)|}{|x-z|^{1+t}}\,\dd z \label{eq:T1} \\[2mm]
			&\quad + G(x,y) \int_{\R\setminus (x-a,x+a)} \frac{\dd z}{|x-z|^{1+t}} \label{eq:T2}\\[2mm]
			&\quad + \int_{(-1,1)\setminus (x-a,x+a)} \frac{G(z,y)}{|x-z|^{1+t}}\,\dd z. \label{eq:T3}
		\end{align}

		\underline{Estimate of \eqref{eq:T1}.}
		
		Let $z\in (x-a,x+a)$, and let $\xi$ be a point between $x$ and $z$. Since the distance function to the boundary is $1$-Lipschitz, we have
		\begin{equation*}
		\delta(\xi) \geqslant \delta(x)-|\xi-x| \geqslant \delta(x)-|z-x|.
		\end{equation*}
		Since $|z-x|<a\leqslant \delta(x)/4$, it follows that
		\begin{equation}\label{eq:delta_xi}
			\delta(\xi)\geqslant \frac{3\, \delta(x)}{4}.
		\end{equation}
		On the other hand, we have
		\begin{equation*}
		      |\xi-y| \geqslant |x-y|-|\xi-x| \geqslant |x-y|-|z-x|.
		\end{equation*}
		Since $|z-x|<a\leqslant \frac{|x-y|}{4}$, we obtain the following
		\begin{equation}\label{eq:distance_xi_y}
			|\xi-y|\geqslant \frac{3\, |x-y|}{4}.
		\end{equation}
		Using \eqref{eq:gx_G}, \eqref{eq:delta_xi} and \eqref{eq:distance_xi_y}, we obtain
		\begin{equation*}
		      |\dr_\xi G(\xi,y)| \leqslant C\, \frac{ \sqrt{\delta(y)} }{|x-y|\, \sqrt{\delta(x)}}.
		\end{equation*}
		Thanks to the mean value theorem, we obtain
		\begin{equation*}
		      |G(x,y)-G(z,y)| \leqslant C\, \frac{ \sqrt{\delta(y)} }{|x-y|\, \sqrt{\delta(x)}} |x-z|.
		\end{equation*}
		Consequently, we have 
		\begin{align}
			\int_{(x-a,x+a)} \frac{|G(x,y)-G(z,y)|}{|x-z|^{1+t}}\,\dd z
			&\leqslant C\, \frac{\sqrt{\delta(y)}}{|x-y|\, \sqrt{\delta(x)} } \int_{(x-a,x+a)} \frac{\dd z}{|x-z|^t} \nonumber\\[2mm]
			&\leqslant  C\, \frac{\sqrt{\delta(y)}}{|x-y|\, \sqrt{\delta(x)} } \left(2\int_0^a r^{-t}\,\dd r\right) \nonumber\\[2mm]
			&\leqslant  \frac{2C}{1-t}\, \frac{\sqrt{\delta(y)}}{|x-y|\, \sqrt{\delta(x)}}\, a^{1-t}.
			\label{eq:T11}
		\end{align}
		We distinguish two cases.
		If $|x-y|\leqslant \delta(x)$, then it holds
		\begin{equation*}
		      \delta(y) \leqslant \delta(x)+|x-y| \leqslant 2\, \delta(x).
		\end{equation*}
		Moreover, we have $a=\frac{|x-y|}{4}$. Hence \eqref{eq:T11} implies that
		\begin{equation*}
		      \int_{(x-a,x+a)} \frac{|G(x,y)-G(z,y)|}{|x-z|^{1+t}}\,\dd z
		      \leqslant C_t\, \frac{ \sqrt{\delta(x)} }{|x-y|\, \sqrt{\delta(x)}}\, |x-y|^{1-t} 
		      = \frac{C_t}{|x-y|^t}.
		\end{equation*}
		If $\delta(x)<|x-y|$, then we have
		\begin{equation*}
		\delta(y) \leqslant \delta(x)+|x-y| \leqslant 2\,|x-y|.
		\end{equation*}
		We also have $a=\delta(x)/4$. Therefore, we have
		\begin{equation*}
		\int_{(x-a,x+a)} \frac{|G(x,y)-G(z,y)|}{|x-z|^{1+t}}\,\dd z
		\leqslant C_t\, \frac{\sqrt{|x-y|}}{|x-y|\, \sqrt{\delta(x)} }\, \delta(x)^{1-t}
		= C_t\, \frac{\delta(x)^{1/2-t}}{|x-y|^{1/2}}.
		\end{equation*}
		We obtain
		\begin{equation}\label{eq:T1final}
			\int_{(x-a,x+a)} \frac{|G(x,y)-G(z,y)|}{|x-z|^{1+t}}\,\dd z
			\leqslant C_t \left[ \frac{1}{|x-y|^t} + \frac{\delta(x)^{1/2-t}}{|x-y|^{1/2}} \right].
		\end{equation}
		
		\underline{Estimate of \eqref{eq:T2}.}
		
		We first compute the following integral
		\begin{equation*}
		\int_{\R\setminus (x-a,x+a)} \frac{\dd z}{|x-z|^{1+t}} = 2\int_a^\infty r^{-1-t}\,\dd r = \frac{2}{t\,a^t}.
		\end{equation*}
		Thus, it holds
		\begin{equation}\label{eq:T21}
			G(x,y) \int_{\R\setminus (x-a,x+a)} \frac{\dd z}{|x-z|^{1+t}} \leqslant C_t\,\frac{G(x,y)}{a^t}.
		\end{equation}
		
		If $|x-y|\leqslant \delta(x)$, we have from Lemma \ref{lm:log-G}
		\begin{equation*}
		G(x,y) \leqslant C\left(1+|\log |x-y||\right).
		\end{equation*}
		Since $a=\frac{|x-y|}{4}$, it follows that
		\begin{equation}\label{eq:T22}
			G(x,y) \int_{\R\setminus (x-a,x+a)} \frac{\dd z}{|x-z|^{1+t}}
			\leqslant C_t\, \frac{1+|\log |x-y||}{|x-y|^t}.
		\end{equation}
		
		If $\delta(x)\leqslant |x-y|$, we apply Lemma \ref{lm:boundary_G} and use $a=\delta(x)/4$, we obtain
		\begin{equation}\label{eq:T23}
			G(x,y) \int_{\R\setminus (x-a,x+a)} \frac{\dd z}{|x-z|^{1+t}}
			\leqslant C_t\, \frac{\sqrt{\delta(x)}}{\sqrt{|x-y|}\,\delta(x)^t}
			= C_t\, \frac{\delta(x)^{1/2-t}}{|x-y|^{1/2}}.
		\end{equation}
		
		Combining \eqref{eq:T22} and
		\eqref{eq:T23}, we conclude that
		\begin{equation}\label{eq:T2final}
			G(x,y) \int_{\R\setminus (x-a,x+a)} \frac{\dd z}{|x-z|^{1+t}}
			\leqslant C_t\left[ \frac{1+|\log |x-y||}{|x-y|^t} + \frac{\delta(x)^{1/2-t}}{|x-y|^{1/2}} \right].
		\end{equation}

		\underline{A local integral estimate for $G$.}
		
		For $0<R\leqslant 1$, we obtain from Lemma \ref{lm:log-G} 
		\begin{equation*}
			\int_{(y-R,y+R)} G(z,y)\,\dd z
			\aleq \int_{(y-R,y+R)} \left(1+|\log|z-y||\right)\,\dd z \aleq \int_0^R(1+|\log r|)\,\dd r.
		\end{equation*}
		Since $0<r\leqslant R\leq1$, one has
		\begin{equation*}
		\int_0^R(-\log r)\,\dd r = \left[-r\log r+r\right]_{r=0}^{r=R} = R\, (1+|\log R|).
		\end{equation*}
		Therefore, we have 
		\begin{equation}\label{eq:int-G}
			\int_{(y-R,y+R)} G(z,y)\,\dd z \aleq R\, \left(1+|\log R|\right).
		\end{equation}
		
		\underline{Estimate of \eqref{eq:T3}.}
		
		We split the integration domain into $E_y \coloneqq \left\{z\in (-1,1) :|z-y|<|x-y|/2\right\}$ and its complement. We denote
		\begin{equation*}
			\int_{(-1,1)\setminus (x-a,x+a)} \frac{G(z,y)}{|x-z|^{1+t}}\,\dd z = I + J,
		\end{equation*}
		where
		\begin{align*}
			I &\coloneqq  \int_{E_y\setminus  (x-a,x+a)}
			\frac{G(z,y)}{|x-z|^{1+t}}\,\dd z,\\[2mm]
			J &\coloneqq  \int_{ (-1,1)\setminus \left[ E_y\cup (x-a,x+a)\right]  } \frac{G(z,y)}{|x-z|^{1+t}} \, \dd z.
		\end{align*}
		For $z\in E_y$, it holds
		\begin{equation*}
		|x-z| \geqslant |x-y|-|z-y| > |x-y|/2.
		\end{equation*}
		Hence, by \eqref{eq:int-G} with $R=|x-y|/2$, it holds
        \begin{equation}\label{eq:I}
		\begin{aligned}
			I &\leqslant \left(\frac{2}{|x-y|}\right)^{1+t}
			\int_{\left( y- |x-y|/2, y+ |x-y|/2\right) } G(z,y)\,\dd z
			\\[2mm]
			&\leqslant C_t\,|x-y|^{-1-t} \, |x-y|\, \left(1+|\log |x-y||\right) \\[2mm]
			&\leqslant C_t\,\frac{1+|\log |x-y||}{|x-y|^t}.
		\end{aligned}
		\end{equation}
		We now estimate $J$. Assume first that $|x-y|\leqslant \delta(x)$. For $z\in (-1,1)\setminus \left[E_y\cup (x-a,x+a)\right]$, we have $|z-y|\geqslant |x-y|/2$.
		Since $|z-y|\leqslant 2$, we obtain from Lemma \ref{lm:log-G}
		\begin{equation*}
		      G(z,y) \leqslant C\left(1+|\log|z-y||\right) \leqslant C\left(1+|\log |x-y||\right).
		\end{equation*}
		We also have $a=\frac{|x-y|}{4}$. Therefore, it holds
        \begin{equation}\label{eq:J1}
		\begin{aligned}
			J&\leqslant C\left(1+|\log |x-y||\right) \int_{\left\{ |x-z| \geqslant \frac{|x-y|}{4} \right\}} \frac{\dd z}{|x-z|^{1+t}}  \\[2mm]
			&\leqslant C_t\, \frac{1+|\log |x-y||}{|x-y|^t}. 
		\end{aligned}
        \end{equation} 
		We now consider the case where $\delta(x)<|x-y|$. By Lemma \ref{lm:boundary_G}, it holds
		\begin{equation*}
		G(z,y) \leqslant C\,\frac{\sqrt{\delta(z)}}{\sqrt{|z-y|}}.
		\end{equation*}
		For $z\in (-1,1)\setminus \left[E_y\cup (x-a,x+a)\right]$, we have $|z-y|\geqslant |x-y|/2$ and hence we obtain 
		\begin{equation}\label{eq:Gzy}
			G(z,y) \leqslant C\,\frac{\sqrt{\delta(z)}}{\sqrt{|x-y|}}.
		\end{equation}
		Since the distance function $\delta$ is $1$-Lipschitz, we have
		\begin{equation*}
		      \delta(z) \leqslant \delta(x)+|x-z| .
		\end{equation*}
		We also have $|x-z|\geqslant a=\delta(x)/4$. It follows that $\delta(x) \leqslant 4\, |x-z|$. 
		Consequently, we have
		\begin{equation}\label{eq:delta_z}
			\delta(z)\leqslant 5|x-z|.
		\end{equation}
		Combining \eqref{eq:Gzy} and
		\eqref{eq:delta_z}, we obtain
		\begin{equation*}
		G(z,y) \leqslant C\,\frac{|x-z|^{1/2}}{|x-y|^{1/2}}.
		\end{equation*}
		Coming back to $J$, we obtain
		\begin{align}
			J&\leqslant \frac{C}{|x-y|^{1/2}}
			\int_{(-1,1)\setminus \left(x-\delta(x)/4,x+\delta(x)/4\right)} \frac{\dd z}{|x-z|^{t+1/2}} \nonumber\\[2mm]
			&\leqslant \frac{C}{|x-y|^{1/2}} \left( 2\int_{\delta(x)/4 }^{2} r^{-t-1/2}\,\dd r \right).
			\label{eq:J2}
		\end{align}
		If $0<t<1/2$, then it holds
		\begin{equation*}
		      \int_{\delta(x)/4}^{2} r^{-t-1/2}\,\dd r 
            = \left[ \frac{r^{1/2-t}}{1/2-t} \right]_{\delta(x)/4}^{2}
		      \leqslant C_t.
		\end{equation*}
		If $t=1/2$, then we obtain 
		\begin{equation*}
		      \int_{\delta(x)/4}^{2}\frac{\dd r}{r} 
              = \log\left(\frac{8}{\delta(x)}\right) \leqslant C\left(1+|\log\delta(x)|\right).
		\end{equation*}
		If $\frac{1}{2}<t<1$, then it holds
		\begin{equation*}
			\int_{\delta(x)/4}^{2}r^{-t-1/2}\,\dd r 
            = \left[ \frac{r^{1/2-t}}{1/2-t} \right]_{\delta(x)/4}^{2} 
            = \frac{\left( \delta(x)/4 \right)^{1/2-t}-2^{1/2-t}}{t-1/2} 
            \leqslant C_t\,\delta(x)^{1/2-t}.
		\end{equation*}
		Thus, in all cases, we obtain 
		\begin{equation}\label{eq:J3}
			J \leqslant C_t\, \frac{\Phi_t(\delta(x))}{|x-y|^{1/2}}.
		\end{equation}
		
		Combining \eqref{eq:I}, \eqref{eq:J1}, and
		\eqref{eq:J3}, we obtain
		\begin{equation}\label{eq:T3final}
			\int_{(-1,1)\setminus (x-a,x+a)} \frac{G(z,y)}{|x-z|^{1+t}}\,\dd z
			\leqslant C_t\left[ \frac{1+|\log |x-y||}{|x-y|^t} 
            + \frac{\Phi_t(\delta(x))}{|x-y|^{1/2}} \right].
		\end{equation}
		
		\underline{Conclusion.}
		
		Combining \eqref{eq:T1final}, \eqref{eq:T2final}, and \eqref{eq:T3final} we conclude that
		\begin{equation*}
		      \left|(-\lap)^{t/2}_x G(x,y)\right| 
            \leq C_t\left[ \frac{1+|\log|x-y||}{|x-y|^t} + \frac{\Phi_t(\delta(x))}{|x-y|^{1/2}} \right].
		\end{equation*}
	\end{proof} 
	
	\underline{For $y\in(-1,1)$ and $x\notin(-1,1)$.}	We prove the following result.
	\begin{proposition}\label{pr:lapt_G_outside}
		Let $y\in(-1,1)$, $x\in\R\setminus (-1,1)$, and $t\in(1/2,1)$. 
		Then it holds
		\begin{equation}\label{eq:lapt_G_outside}
			\left|(-\lap)^{t/2}_x G(x,y)\right| \aleq_t
			\begin{cases}
				\displaystyle \frac{1}{(1+|x|)^{1+t}} & \text{ if }|x|\geq2,\\[4mm]
				\displaystyle \frac{1+|\log|x-y||}{|x-y|^t} + \frac{1} {\delta(x)^{t-1/2}|x-y|^{1/2}} & \text{ if } 1<|x|<2.
			\end{cases}
		\end{equation}
	\end{proposition}
	
	\begin{proof}
		It holds
		\begin{equation*}
			(-\lap)^{t/2}_x G(x,y) = C_t\, \int_{\R} \frac{G(x,y)-G(z,y)}{|x-z|^{1+t}}\,\dd z 
            = -C_t\, \int_{-1}^{1} \frac{G(z,y)}{|x-z|^{1+t}}\,\dd z.
		\end{equation*}
		Since $G(z,y)\geqslant 0$, it follows that
		\begin{equation}\label{eq:abs_lapt_G}
			\left|(-\lap)^{t/2}_x G(x,y)\right| 
            = C_t\, \int_{-1}^{1} \frac{G(z,y)}{|x-z|^{1+t}}\,\dd z.
		\end{equation}
		Lemma \ref{lm:log-G} implies that the $\Ll^1$-norm of $G(\cdot,y)$ is bounded uniformly with respect to $y$
		\begin{equation*}
			\int_{-1}^{1} G(z,y)\,\dd z
			\leqslant C\int_{-1}^{1} \left(1+\left|\log|z-y|\right|\right)\,\dd z 
            \leqslant C\int_{-2}^{2} \left(1+|\log|w||\right)\,\dd w 
            \leqslant C.
		\end{equation*}
		
		\underline{Case 1: $|x|\geqslant 2$.}
		
		For $z\in(-1,1)$ and $|x|\geqslant 2$, it holds
		\begin{equation*}
		\frac{1}{|x-z|^{1+t}} \leqslant \frac{C_t}{(1+|x|)^{1+t}}.
		\end{equation*}
		Inserting this into \eqref{eq:abs_lapt_G} yields
		\begin{equation*}
			\left|(-\lap)^{t/2}_x G(x,y)\right|
			\leqslant \frac{C_t}{(1+|x|)^{1+t}} \int_{-1}^{1}G(z,y)\,\dd z 
            \leqslant \frac{C_t}{(1+|x|)^{1+t}}.
		\end{equation*}
		
		\underline{Case 2: $1<|x|<2$.}
		It holds $|x-y| \geqslant \delta(x)$.
        Moreover, because $|x|<2$ and $|y|<1$, we obtain $0<|x-y|<3$.
		We decompose $(-1,1)$ into
		\begin{equation*}
        \begin{aligned} 
			I_1& \coloneqq \left\{ z\in(-1,1):|z-y|\geqslant |x-y|/2 \right\}, \\[2mm]
			I_2 & \coloneqq \left\{ z\in(-1,1):|z-y|<|x-y|/2 \right\}.
        \end{aligned} 
		\end{equation*}
		We obtain
		\begin{align}
			\left|(-\lap)^{\frac t2}_xG(x,y)\right|
			&\leqslant C_t\int_{I_1} \frac{1}{|x-z|^{1+t}}\, \frac{\sqrt{\delta(z)}}{\sqrt{|z-y|}}\, \dd z \label{eq:I1}\\[2mm]
			&\quad+ C_t\int_{I_2} \frac{1+|\log|z-y||}{|x-z|^{1+t}}\,\dd z.  \label{eq:I2}
		\end{align}
		We first estimate \eqref{eq:I1}. For $z\in I_1$ and $x>1$, we have 
		\begin{equation*}
			\delta(z) = 1-z \leqslant (x-1) + (1-z) = x-z=|x-z|.
		\end{equation*}
		If $x<-1$, we have 
		\begin{equation*}
			\delta(z) = 1+z \leqslant (-1-x) + (1+z) = z-x=|x-z|.
		\end{equation*}
		Hence, it holds
		\begin{equation}\label{eq:I1_1}
			\int_{I_1} \frac{1}{|x-z|^{1+t}}\, \frac{\sqrt{\delta(z)}}{\sqrt{|z-y|}}\,\dd z
			\leqslant \frac{\sqrt{2}}{\sqrt{|x-y|}} \int_{-1}^{1} \frac{\dd z}{|x-z|^{1/2+t}}.
		\end{equation}
		If $x>1$, we have
        \begin{equation*} 
		\begin{aligned}
			\int_{-1}^{1}\frac{\dd z}{|x-z|^{1/2+t}}
			&= \int_{x-1}^{x+1}\frac{\dd r}{r^{1/2+t}}\\[2mm]
			&= \frac{1}{t-1/2} \left( (x-1)^{1/2-t} - (x+1)^{1/2-t} \right)\\[2mm]
			&\aleq_t \delta(x)^{1/2-t}.
		\end{aligned}
        \end{equation*}
		A similar computation holds for $x<-1$. Returning to \eqref{eq:I1_1}, we obtain
		\begin{equation}\label{eq:I1_2}
			\int_{I_1} \frac{1}{|x-z|^{1+t}}\, \frac{\sqrt{\delta(z)}}{\sqrt{|z-y|}}\,\dd z 
            \leqslant C_t\, \frac{ \delta(x)^{1/2-t} }{ |x-y|^{1/2} }.
		\end{equation}
		We now estimate \eqref{eq:I2}. If $z\in I_2$, it holds
		\begin{equation*}
			|x-z| \geqslant |x-y|-|z-y| > |x-y|/2.
		\end{equation*}
		Consequently, we have
        \begin{equation}\label{eq:I2_1}
		\begin{aligned}
			\int_{I_2} \frac{1+|\log|z-y||}{|x-z|^{1+t}}\,\dd z
			&\leqslant \left(\frac{2}{|x-y|}\right)^{1+t} \int_{\left\{|z-y|<|x-y|/2\right\}} \left(1+|\log|z-y||\right)\,\dd z  \\[2mm]
			&\leqslant \frac{C_t}{|x-y|^{1+t}} \int_0^{|x-y|/2} (1+|\log r|)\,\dd r. 
		\end{aligned}
        \end{equation}
		To estimate the last integral, we consider the change of variable $r=|x-y|\, u$. We obtain
        \begin{equation*} 
		\begin{aligned}
			\int_0^{|x-y|/2}(1+|\log r|)\,\dd r
			&= |x-y| \int_0^{1/2} \left(1+|\log |x-y|+\log u|\right)\,\dd u\\[2mm]
			&\leqslant C\, |x-y|\, \left(1+|\log |x-y| |\right).
		\end{aligned}
        \end{equation*}
		It follows from \eqref{eq:I2_1} that
		\begin{equation}\label{eq:I2_2}
			\int_{I_2} \frac{1+|\log|z-y||}{|x-z|^{1+t}}\,\dd z 
            \leqslant C_t\, \frac{1+|\log |x-y||}{|x-y|^t}.
		\end{equation}
		Combining \eqref{eq:I1_1}, \eqref{eq:I2_1}, \eqref{eq:I1_2} and \eqref{eq:I2_2}, we conclude that
		\begin{equation*}
		\left|(-\lap)^{t/2}_x G(x,y)\right| \aleq_t \frac{1+|\log |x-y||}{|x-y|^t} + \frac{ \delta(x)^{1/2-t} }{ |x-y|^{1/2} }.
		\end{equation*}
	\end{proof}
	
	\subsection{Global regularity}\label{sec:global_reg}
	
	In this section, we study the regularity up to the boundary for Dirichlet problem of the $\frac{1}{2}$-Laplacian. We obtain first the following lemma concerning the solution itself (i.e.\ not its derivatives).
	
	\begin{lemma}\label{lm:boundary_u}
		We consider $f\in \Ll^p(-1,1)$ for some $2>p>1$ and $u\colon \R\to \R$ be the solution to 
		\begin{equation*}
			\begin{cases}
				\displaystyle (-\lap)^{1/2}u = f & \text{ in }(-1,1) \\[2mm]
				u = 0 & \text{ in } \R\setminus (-1,1).
			\end{cases}
		\end{equation*}
		Then it holds
		\begin{equation*}
			\|u\|_{\Ll^{\infty}(-1,1)} + \left\| \frac{u}{\delta^{1/2}} \right\|_{\Ll^{\frac{2p}{2-p}}(-1,1)} \aleq_{s,t,\eps} \|f\|_{\Ll^p(-1,1)}.
		\end{equation*}
	\end{lemma}
	
	\begin{proof}
		By Lemma \ref{lm:log-G}, it holds
		\begin{equation*}
			|u(x)|  \aleq \int_{-1}^1 \Big(1 + |\log|x-y||\Big)\, |f(y)|\, \dd y.
		\end{equation*}
		Since $|\log|x||\in \Ll^{p'}$, we have by Young inequality for the convolution that 
		\begin{equation*}
			\|u\|_{\Ll^{\infty}(-1,1)} \aleq \|f\|_{\Ll^p(-1,1)}.
		\end{equation*}
		Moreover, we deduce from Lemma \ref{lm:boundary_G} that 
		\begin{equation*}
			\frac{|u(x)|}{\delta(x)^{1/2}} \aleq \int_{-1}^1 \frac{|f(y)|}{|x-y|^{1/2}}\, \dd y.
		\end{equation*}
		We conclude by Riesz potentials estimates \cite[Theorem 6.1.3]{Gr2014}.
	\end{proof}
	
	In order to obtain higher-order estimates, we first derive a pointwise decomposition of the fractional Laplacian of the solution.
	
	\begin{lemma}\label{lm:pointwise_lapt_u}
		Let $f\in \Ll^p(-1,1)$ for some $p>1$, and let $u\colon \R \to \R$ be the solution of
		\begin{equation*}
		\begin{cases}
			(-\lap)^{1/2} u=f & \text{in }(-1,1),\\[1mm]
			u=0 & \text{in }\R\setminus(-1,1).
		\end{cases}
		\end{equation*}
		Let $1/2<t<1$. Then, for almost every $x\in(-1,1)$,
		\begin{equation}
			\left|(-\lap)^{t/2} u(x)\right| \aleq_t
			\int_{-1}^{1} \frac{1+|\log|x-y||}{|x-y|^t} |f(y)|\,|\dd y|+ \frac{1}{\delta(x)^{t-1/2}} \int_{-1}^{1} \frac{|f(y)|}{|x-y|^{1/2}}\,\dd y.
			\label{eq:pointwise_lapt_u}
		\end{equation}
	\end{lemma}
	
	\begin{proof}
		The Green representation formula gives
		\begin{equation}\label{eq:rep_formula}
			\forall x\in\R,\qquad u(x) = \int_{-1}^{1} G(x,y)f(y)\,\dd y.
		\end{equation}
		By Proposition \ref{pr:lapt_G}, for $x,y\in(-1,1)$, it holds
		\begin{equation*}
		      \left|(-\lap)^{t/2}_x G(x,y)\right| 
              \aleq_t \frac{1+|\log|x-y||}{|x-y|^t} + \frac{1}{ \delta(x)^{t-1/2}|x-y|^{1/2} }.
		\end{equation*}
		Applying $(-\lap)^{t/2}$ to \eqref{eq:rep_formula}, we obtain for almost every $x\in(-1,1)$
		\begin{equation*}
		      (-\lap)^{t/2} u(x) = \int_{-1}^{1} (-\lap)^{t/2}_x G(x,y)f(y)\,\dd y.
		\end{equation*}
		We obtain the following estimate
		\begin{align*}
			\left|(-\lap)^{t/2} u(x)\right|
			&\leqslant \int_{-1}^{1} \left|(-\lap)^{t/2}_x G(x,y)\right|\, |f(y)|\,\dd y\\[2mm]
			&\aleq_t \int_{-1}^{1} \frac{1+|\log|x-y||}{|x-y|^t}\, |f(y)|\,\dd y + \frac{1}{\delta(x)^{t-1/2}} \int_{-1}^{1} \frac{|f(y)|}{|x-y|^{1/2}}\,\dd y.
		\end{align*}
		This proves \eqref{eq:pointwise_lapt_u}.
	\end{proof}

	We can now gather the previous result to obtain global estimates on $\R$.
	
	\begin{theorem}\label{th:global_reg_Lpdata}
		Let $1/2<t<1$, $0<\eps<1-t$, and $1<p<2$.
		We define $q_{\eps}>1$ such that 
		\begin{equation}\label{eq:def_q_epsilon}
			\frac{1}{q_{\eps}} \coloneqq  \frac{1}{p}-(1-t-\eps) = \frac{1}{p} + t + \eps -1.
		\end{equation}
        We consider $f\in \Ll^p(-1,1)$ and $u\colon \R\to \R$ be the solution to 
		\begin{equation*}
			\begin{cases}
				\displaystyle (-\lap)^{1/2}u = f & \text{ in }(-1,1) \\[2mm]
				u = 0 & \text{ in } \R\setminus (-1,1).
			\end{cases}
		\end{equation*}
		Then it holds
		\begin{equation}\label{eq:global_reg_Lpdata}
			\left\|(-\lap)^{t/2} u \right\|_{\Ll^{q_{\eps}}(\R)} \aleq_{p,t,\eps} \|f\|_{\Ll^p(-1,1)}.
		\end{equation}
	\end{theorem}

	\begin{proof}
		We extend $f$ by zero outside $(-1,1)$. For $0<\alpha<1$, we denote the Riesz potential
		\begin{equation*}
		      I_\alpha h(x) \coloneqq  \pv \int_{\R} \frac{h(y)}{|x-y|^{1-\alpha}}\,\dd y.
		\end{equation*}
		We denote $q_0$ the number obtained by the formula \eqref{eq:def_q_epsilon} for $\eps=0$.
		We have
		\begin{equation}\label{eq:qepsilon_qzero}
			\frac{1}{q_{\eps}} = \frac{1}{q_0}+\eps.
		\end{equation}
		In particular, it holds $q_{\eps}<q_0$. We divide the proof into three steps.
		
		\underline{Step 1: estimate on $(-2,2)$.}
		
		By Lemma \ref{lm:pointwise_lapt_u} in $(-1,1)$, and \eqref{eq:lapt_G_outside} with \eqref{eq:rep_formula} in $(-2,2)\setminus[-1,1]$, we have for $x\in(-2,2)$
		\begin{equation}\label{eq:lapt_u_two_terms}
			\left|(-\lap)^{t/2} u(x)\right| \aleq_t A(x) + B(x),
		\end{equation}
		where
		\begin{align}
			A(x)& \coloneqq \int_{-1}^{1} \frac{1+|\log|x-y||}{|x-y|^t}\, |f(y)|\,\dd y, \label{eq:A}\\[2mm]
			B(x)& \coloneqq \delta(x)^{1/2-t} \int_{-1}^{1} \frac{|f(y)|}{|x-y|^{1/2}}\,\dd y.\label{eq:B}
		\end{align}
		
		\underline{Estimate of \eqref{eq:B}.} We have $B(x) = \delta(x)^{1/2-t}I_{1/2}|f|(x)$.
		Since $1<p<2$, we obtain by Riesz potential estimates \cite[Theorem 6.1.3]{Gr2014},
		\begin{equation}\label{eq:HLS_one_half_rewritten}
			\left\|I_{1/2}|f|\right\|_{\Ll^{\frac{2p}{2-p}}(\R)} \aleq_p \|f\|_{\Ll^p(-1,1)},
		\end{equation}
		Since $t>1/2$, we have 
		\begin{equation}\label{eq:rho_weak_space}
			\delta^{1/2-t} \in L^{\left(\frac{1}{t-1/2},\infty\right)}(-2,2).
		\end{equation}
		Indeed, for $\lambda\geqslant 1$, it holds
		\begin{equation*}
		\left| \left\{ x\in(-2,2): \delta(x)^{1/2-t}>\lambda \right\} \right| 
		= \left| \left\{ x\in(-2,2): \delta(x)<\lambda^{-\frac{1}{t-1/2}} \right\}\right|
		\leqslant C\lambda^{-\frac{1}{t-1/2}}.
		\end{equation*}
		This is exactly \eqref{eq:rho_weak_space}.
		
		Moreover, it holds
		\begin{equation*}
			t-\frac{1}{2} + \frac{2-p}{2p} = \frac{1}{p} + t-1 = \frac{1}{q_0}.
		\end{equation*}
		By Hölder inequality, we obtain
		\begin{equation}\label{eq:B_weak_qzero}
			\|B\|_{\Ll^{(q_0,\infty)}(-2,2)} \aleq_{p,t} \left\| \delta^{1/2-t} \right\|_{\Ll^{\left(\frac{1}{t-1/2},\infty\right)}(-2,2)} \,  \left\|I_{1/2}|f|\right\|_{\Ll^{\frac{2p}{2-p} }(\R)}  \aleq_{p,t} \|f\|_{\Ll^p(-1,1)},
		\end{equation}
		By \eqref{eq:qepsilon_qzero}, we have $q_\eps<q_0$. Hence \eqref{eq:B_weak_qzero} implies that
		\begin{equation}\label{eq:B_strong_qepsilon}
			\|B\|_{\Ll^{q_\eps}(-2,2)} \aleq_{p,t,\eps} \|f\|_{\Ll^p(-1,1)}.
		\end{equation}
		
		\underline{Estimate of \eqref{eq:A}.}
		
		Let $0<\eps<1-t$. For every $0 < r \leqslant 3$, we have
		\begin{equation}\label{eq:log_absorption}
			1+|\log r| \leqslant C_\eps\, r^{-\eps}.
		\end{equation}
		We obtain, with $0<1-t-\eps<1/2$, that
		\begin{equation*}
			A(x) \aleq_{\eps} \int_{-1}^{1} \frac{|f(y)|}{|x-y|^{t+\eps}}\,\dd y = I_{1-t-\eps}|f|(x).
		\end{equation*}
		Since $1<p<2$, we have
		\begin{equation*}
			p<2<\frac{1}{1-t-\eps}.
		\end{equation*}
		By Riesz potential estimates \cite[Theorem 6.1.3]{Gr2014}, we have
		\begin{equation}\label{eq:A_qepsilon}
			\|A\|_{\Ll^{q_\eps}(\R)} \aleq_{p,t,\eps} \|f\|_{\Ll^p(-1,1)},
		\end{equation}		
		Combining \eqref{eq:lapt_u_two_terms}, \eqref{eq:B_strong_qepsilon}, and \eqref{eq:A_qepsilon}, we conclude that
		\begin{equation}\label{eq:lapt_u_1}
			\left\|(-\lap)^{t/2} u\right\|_{\Ll^{q_\eps}(-2,2)} \aleq_{p,t,\eps} \|f\|_{\Ll^p(-1,1)}.
		\end{equation}
		
		\underline{Step 2: estimate on $\R\setminus(-2,2)$.}
		
		For $|x|\geq2$, using Proposition \ref{pr:lapt_G_outside} and \eqref{eq:rep_formula}, we obtain
		\begin{equation*}
			\left|(-\lap)^{t/2} u(x)\right|
			\leqslant \int_{-1}^{1} \left|(-\lap)^{t/2}_x G(x,y)\right|\, |f(y)|\,\dd y 
            \aleq_t \frac{1}{(1+|x|)^{1+t}} \int_{-1}^{1} |f(y)|\,\dd y 
            \aleq_{t,p} \frac{\|f\|_{\Ll^p(-1,1)}}{(1+|x|)^{1+t}}.
		\end{equation*}
		We have $(1+|x|)^{-1-t} \in \Ll^{q_\eps} (\R\setminus(-2,2))$ since $(1+t)q_{\eps}>1$.
		we deduce that
		\begin{equation}\label{eq:lapt_u_2}
			\left\|(-\lap)^{t/2} u\right\|_{\Ll^{q_\eps}(\R\setminus(-2,2))} \aleq_{p,t,\eps} \|f\|_{\Ll^p(-1,1)}.
		\end{equation}
		
		\underline{Step 3: conclusion.}
		
		Combining \eqref{eq:lapt_u_1} and \eqref{eq:lapt_u_2}, we complete the proof.
	\end{proof}
	
	As an application of Theorem \ref{th:global_reg_Lpdata} and Lemma \ref{lm:boundary_u}, we obtain the following estimates for $p=2$.
	\begin{theorem}\label{th:global_reg_L2data}
		Let $f\in \Ll^2(-1,1)$ and $u\colon\R\to \R$ be the solution to 
		\begin{equation}\label{eq:equation-u}
			\begin{cases}
				\displaystyle (-\lap)^{1/2} u = f & \text{ in }(-1,1),\\[2mm]
				u = 0 & \text{ in }\R\setminus (-1,1).
			\end{cases}
		\end{equation} 
		Then for any $p\in[1,\infty)$, we have $\frac{u}{\delta^{1/2}}\in \Ll^p(-2,2)$ and $u\in \W^{1/2,p}(-2,2)$, with the estimate
		\begin{equation*}
			\left\| \frac{u}{\delta^{1/2}} \right\|_{\Ll^p(-2,2)} + \|u\|_{\W^{1/2,p}(-2,2)} \aleq_p \|f\|_{\Ll^2(-1,1)}.
		\end{equation*}
		In particular, it holds $u\in \Cd^{0,\alpha}(\R)$ for any $\alpha\in(0,1/2)$.
	\end{theorem}
	The inclusion $\W^{1/2,p}\hookrightarrow \Cd^{0,1/2 - 1/p}$ can be found in \cite[Theorem 8.2]{DiNPV2012}.
	
	\begin{lemma}\label{lm:C0half}
		\label{lem:half_laplacian_p_above_two}
		Let $p>2$, let $f\in \Ll^p(-1,1)$, and let $u$ be given by a solution to \eqref{eq:equation-u}. Then there exists $C=C(p)>0$ such that
		\begin{equation}\label{eq:Chalf_estimate}
			\|u\|_{\Cd^{0,1/2}(\R)} \leqslant C\|f\|_{\Ll^p(-1,1)}.
		\end{equation}
		Moreover, it holds
		\begin{align}
			\lim_{r\to 0^+}\frac{u(1-r)}{\sqrt r} &= \frac{\sqrt{2}}{\pi} \int_{-1}^1 \sqrt{\frac{1+y}{1-y}}\,f(y)\,\dd y,
			\label{eq:right}\\
			\lim_{r\to 0^+}\frac{u(-1+r)}{\sqrt r} &= \frac{\sqrt2}{\pi} \int_{-1}^1 \sqrt{\frac{1-y}{1+y}}\,f(y)\,\dd y. \label{eq:left}
		\end{align}
	\end{lemma}
	
	\begin{proof}
		Since $p>2$, we have
		\begin{equation*}
			\frac{1}{p'}\coloneqq \frac{p-1}{p}>\frac{1}{2} \qquad\text{and}\qquad 1<p'<2 .
		\end{equation*}
		The $\Ll^\infty$-estimate follows from Lemma \ref{lm:boundary_u}.
		
		Let $x,z\in(-1,1)$, and denote $h\coloneqq |x-z|$.
		By \eqref{eq:rep_formula} and Hölder's inequality, we have
		\begin{equation}\label{eq:diff-u}
			|u(x)-u(z)| \leqslant \|G(x,\cdot)-G(z,\cdot)\|_{\Ll^{p'}(-1,1)}\, \|f\|_{\Ll^p(-1,1)}.
		\end{equation}
		
		We first consider the case where
		\begin{equation*}
			\min\{\delta(x),\delta(z)\} \leqslant 4h.
		\end{equation*}
		Since $\delta$ is $1$-Lipschitz, we obtain
		\begin{equation*}
		\max\{\delta(x),\delta(z)\} \leqslant \min\{\delta(x),\delta(z)\}+h \leqslant 5 h.
		\end{equation*}
		By Lemma \ref{lm:boundary_G} and using $p'<2$, we obtain
        \begin{equation*}
		\begin{aligned}
			\|G(x,\cdot)-G(z,\cdot)\|_{\Ll^{p'}(-1,1)} &\leq
			\|G(x,\cdot)\|_{\Ll^{p'}(-1,1)} + \|G(z,\cdot)\|_{\Ll^{p'}(-1,1)}\\[2mm]
			&\aleq \sqrt{\delta(x)}+\sqrt{\delta(z)}\\[2mm]
			&\aleq h^{1/2}.
		\end{aligned}
        \end{equation*}
 		We now consider the case where 
		\begin{equation}\label{eq:m}
			m\coloneqq \min\{\delta(x),\delta(z)\}>4h.
		\end{equation}
		Let
		\begin{equation*}
		      E_0 \coloneqq  \{y\in(-1,1):|y-x|<2h\} \cup \{y\in(-1,1):|y-z|<2h\}.
		\end{equation*}
		By Lemma \ref{lm:log-G}, and using $|E_0| \aleq h$, we obtain 
		\begin{equation}\label{eq:E0}
			\|G(x,\cdot)-G(z,\cdot)\|_{\Ll^{p'}(E_0)} \aleq h^{1/p'}\big(1+|\log h|\big) \aleq h^{1/2}.
		\end{equation}
		It remains to estimate the difference on $E_1\coloneqq (-1,1)\setminus E_0$. Let $I_{x,z}$ be the segment joining $x$ and $z$, and let
		\begin{equation*}
			d(y)\coloneqq \dist(y,I_{x,z}).
		\end{equation*}
		For $y\in E_1$, we have $d(y)\geqslant 2h$. Moreover, for every $\xi\in I_{x,z}$, we have
		\begin{equation*}
		      \delta(\xi) \geqslant m-h \geqslant \frac{3m}{4}.
		\end{equation*}
		If $\eta\in I_{x,z}$ satisfies $|y-\eta|=d(y)$, then it holds
		\begin{equation*}
			d(y) \leqslant \delta(\eta)+d(y) \aleq m+d(y).
		\end{equation*}
		Using \eqref{eq:gx_G} and the mean-value theorem, we therefore obtain
		\begin{align*}
			|G(x,y)-G(z,y)|
			&\aleq h \sup_{\xi\in I_{x,z}} \frac{\sqrt{\delta(y)}}{|\xi-y|\, \sqrt{\delta(\xi)}} \aleq h\, \frac{\sqrt{m+d(y)}}{d(y)\, \sqrt m}\\[2mm]
			&\aleq h\left( \frac{1}{d(y)} + \frac{1}{\sqrt{m\, d(y)}} \right).
		\end{align*}
		Taking the $\Ll^{p'}$-norm, we obtain by using $1<p'<2$
		\begin{equation*}
        \begin{aligned} 
			\|G(x,\cdot)-G(z,\cdot)\|_{\Ll^{p'}(E_1)}
			& \aleq_p h\,\left( \int_{2h}^{2}\frac{\dd r}{r^{p'}} \right)^{1/p'} 
            + \frac{h}{\sqrt m} \left( \int_{2h}^{2}\frac{\dd r}{r^{p'/2 }} \right)^{1/p'} \\[2mm]
            & \aleq_p h^{1/p'} + \frac{h}{\sqrt{m}}.
        \end{aligned} 
		\end{equation*}
		Since $m>4h$ by \eqref{eq:m} and $1/p'>1/2$, we end up with 
		\begin{equation}\label{eq:E1}
			\|G(x,\cdot)-G(z,\cdot)\|_{\Ll^{p'}(E_1)}  \aleq_p h^{1/2}.
		\end{equation}
		Combining \eqref{eq:E0} and \eqref{eq:E1}, we obtain
		\begin{equation*}
		      \|G(x,\cdot)-G(z,\cdot)\|_{\Ll^{p'}(-1,1)} \aleq_p |x-z|^{1/2}.
		\end{equation*}
        Returning to \eqref{eq:diff-u}, we conclude that
		\begin{equation}\label{eq:holder}
			|u(x)-u(z)| \aleq_p |x-z|^{1/2}\, \|f\|_{\Ll^p(-1,1)}.
		\end{equation}
		Since $u=0$ on $\R\setminus(-1,1)$, it remains only to compare an interior point with an exterior point. If $x\in(-1,1)$ and $z\notin(-1,1)$, then $\delta(x)\leqslant |x-z|$.
		By Lemma \ref{lm:boundary_G}, we obtain 
		\begin{align*}
			|u(x)-u(z)|
			&= |u(x)| \leqslant \|G(x,\cdot)\|_{\Ll^{p'}(-1,1)}\, \|f\|_{\Ll^p(-1,1)}\\[2mm]
			&\aleq_p \sqrt{\delta(x)}\, \|f\|_{\Ll^p(-1,1)}\\[2mm]
			&\aleq_p |x-z|^{1/2}\, \|f\|_{\Ll^p(-1,1)}.
		\end{align*}
		This proves \eqref{eq:Chalf_estimate}.
		
		We now prove the boundary asymptotic expansions. For $x=1-r$, we denote
		\begin{equation*}
			K_r(y)\coloneqq \frac{G(1-r,y)}{\sqrt r}.
		\end{equation*}
		Using \eqref{eq:Green}, we have
		\begin{equation*}
		K_r(y) = \frac{1}{\pi\sqrt{r}} \log\left( \frac{ 1-y+ry+\sqrt{(2r-r^2)(1-y^2)} }{ |1-y-r| } \right).
		\end{equation*}
		For every fixed $y\in(-1,1)$ and $r$ small enough, we have
		\begin{equation*}
			K_r(y) = \frac{1}{\pi\sqrt{r}} \log\left( 1+ \frac{ r(1+y)+\sqrt{(2r-r^2)(1-y^2)} }{ 1-y-r } \right) \xrightarrow[r\to 0]{} K_+(y) \coloneqq \frac{\sqrt{2}}{\pi} \sqrt{\frac{1+y}{1-y}}.
		\end{equation*}
		It remains to justify the convergence under the integral. We have $K_r\to K_+$ pointwise almost everywhere and by Lemma \ref{lm:boundary_G}, there exists $s\in(p',2)$ such that
		\begin{equation*}
		      \sup_{0<r<1}\|K_r\|_{\Ll^s(-1,1)}<+\infty.
		\end{equation*}
		Fatou's lemma gives $K_+\in \Ll^s(-1,1)$. Consequently, the family $(|K_r-K_+|^q)_{0<r<1}$ is uniformly integrable. Vitali's convergence theorem therefore gives the strong convergence in $\Ll^{p'}(-1,1)$.
		It follows from Hölder's inequality that
        \begin{equation*}
		\begin{aligned}
			\lim_{r\to 0^+} \frac{u(1-r)}{\sqrt{r}}
			&= \lim_{r\to 0^+} \int_{-1}^{1} K_r(y)\, f(y)\,\dd y\\[2mm]
			&= \int_{-1}^{1} K_+(y)\, f(y)\,\dd y\\[2mm]
			&= \frac{\sqrt{2}}{\pi} \int_{-1}^{1} \sqrt{\frac{1+y}{1-y}}\,f(y)\,\dd y.
		\end{aligned}
        \end{equation*} 
		This proves \eqref{eq:right}.
		
		Finally, using the symmetry $G(-x,-y)=G(x,y)$
		we similarly obtain
		\begin{equation*}
		\lim_{r\to 0^+}\frac{u(-1+r)}{\sqrt{r}} = \frac{\sqrt{2}}{\pi} \int_{-1}^{1} \sqrt{\frac{1-y}{1+y}}\,f(y)\,\dd y.
		\end{equation*}
		This proves \eqref{eq:left}.
	\end{proof}

	\bibliography{mobius.bib}
	\bibliographystyle{plain}
	
\end{document}